\documentclass[3p,sort&compress]{elsarticle}

\usepackage{amsmath,amsfonts,amssymb}
\usepackage{subfigure,epsfig,psfrag}
\usepackage{epstopdf}
\usepackage{color,url}
\usepackage{algorithm,algorithmic}

\newcommand{\mJ}{\mathcal{J}}

\newcommand{\bs}{{\bf s}}
\newcommand{\bg}{{\bf g}}
\newcommand{\bv}{{\bf v}}
\newcommand{\RR}{\mathbb{R}}

\def\bnabla{\boldsymbol{\nabla}}

\newcommand\norm[1]{\left\Vert{} #1 \right\Vert}
\newcommand\lnorm[1]{\left\vert{} #1 \right\vert}
\newcommand\seminorm[1]{\left[ #1 \right]}
\newcommand\D[2]{\frac{\partial{} #1}{\partial{} #2}}

\newcommand{\Ind}{\mathbf{1}}

\def\J{\mathcal{J}}

\let\l\ell

\def\ivspace{B_2^{3/2 + 2\alpha}}

\def\scontrolspace{W_2^2}
\def\gcontrolspace{B_2^{1/2 + \alpha}}
\def\fcontrolspace{B_{2, x, t}^{1, 1/4 + \alpha}}

\def\muspace{B_2^{1/4}}
\def\adjointsolnspace{W_2^{2, 1}}
\def\chigammaspace{B_{2,x,t}^{3/2 + 2\alpha^*,3/4 + \alpha^*}}

\newcommand{\argmin}{\operatorname{argmin}}

\newtheorem{thm}{Theorem}

\newtheorem{cor}[thm]{Corollary}
\newdefinition{rmk}{Remark}
\newproof{pf}{Proof}
\newproof{pot}{Proof of Theorem \ref{thm2}}
\newtheorem{definition}{Definition}[section]
\DeclareMathOperator*{\esssup}{ess~sup}

\begin{document}

\begin{frontmatter}
\title{Gradient Method in Hilbert-Besov Spaces for the Optimal Control of Parabolic Free Boundary Problems}
\author[fit]{Ugur G.~Abdulla\corref{cor1}}
\ead{abdulla@fit.edu}
\cortext[cor1]{Corresponding author}
\author[fit]{Vladislav Bukshtynov}
\ead{vbukshtynov@fit.edu}
\author[fit]{Ali Hagverdiyev}
\ead{ahaqverdiyev2011@my.fit.edu}
\address[fit]{Department of Mathematical Sciences,
Florida Institute of Technology, Melbourne, FL 32901, USA}

\begin{abstract}
This paper presents computational analysis of the inverse Stefan type free boundary problem, where information
on the boundary heat flux is missing and must be
found along with the temperature and the free boundary.
We pursue optimal control framework introduced in {\it U.G. Abdulla, Inverse Problems and Imaging, 7, 2(2013), 307-340;\ 10, 4(2016), 869--898}, where boundary heat flux and free boundary are components of the control vector,
and optimality criteria consist of the minimization of the quadratic declinations from the available measurements of the temperature distribution at the final moment, phase transition temperature on the free boundary, and  the final position of the free boundary. We develop gradient descent algorithm based on Frechet differentiability in Hilbert-Besov spaces complemented with preconditioning or increase of regularity of the Frechet gradient through implementation of the Riesz representation theorem. Three model examples with various levels of complexity are considered.
Extensive comparative analysis through implementation of preconditioning and Tikhonov regularization, calibration of preconditioning and regularization parameters, effect of noisy data, comparison of simultaneous identification of control parameters vs. nested optimization is pursued.

\end{abstract}


\end{frontmatter}

\section{Introduction and Motivation}
\label{sec:intro}
The goal of this paper is to implement and analyze gradient method in Besov spaces framework for the numerical solution of the optimal control problem introduced recently as a variational formulation of the {\it inverse Stefan problem} (ISP) in~\cite{abdulla13,abdulla15}.
Consider the general one-phase Stefan problem:
\begin{gather}
  Lu := (a(x,t) u_x)_x + b(x,t) u_x + c(x,t) u - u_{t} = f(x,t),~\mathrm{in}~\Omega\label{eq:pde-1}
\\
  u(x,0) = \phi (x),~0 \leq x \leq s(0)=s_0\label{eq:pde-init}
\\
  a(0,t) u_x (0,t) = g(t),~0 \leq t \leq T\label{eq:pde-bound}
\\
  a(s(t),t) u_x (s(t),t) + \gamma (s(t),t) s'(t) = \chi (s(t),t),~0 \leq t \leq T\label{eq:pde-stefan}
\\
  u(s(t),t) = \mu (t),~0 \leq t \leq T,\label{eq:pde-freebound}
\end{gather}
where
\begin{equation}
  \Omega = \{(x,t): 0 < x < s(t),~0 < t \leq T\}\label{eq:pde-domain-defn}
\end{equation}
and $a, b, c, f, \phi, g, \gamma, \chi, \mu$ are given functions.
Assume now that some of the data is not available, or involves some measurement error.  For example, assume that the heat flux $g(t)$ on the fixed boundary $x=0$ is not known and must be found along with the temperature $u(x,t)$ and the free boundary $s(t)$.  In order to do that, some additional information is needed.  Assume that we are able to measure the temperature on our domain and the position of the free boundary at the final moment $T$.
\begin{equation}
  u(x,T) = w(x),~ 0 \leq x \leq s(T) = s_*.\label{eq:pde-finaltemp}
\end{equation}
Under these conditions, we are required to solve an \emph{inverse} Stefan problem (ISP): find a triple
  $\{u,s,g\}$
that satisfies conditions~\eqref{eq:pde-1}--\eqref{eq:pde-finaltemp}.

The motivation for this type of inverse problem arose, in particular, from the modeling of bioengineering problems on the laser ablation of biological tissues through a Stefan problem~\eqref{eq:pde-1}--\eqref{eq:pde-finaltemp}, where $s(t)$ is the ablation depth at the moment $t$.
The boundary temperature measurement $u(0,t)$ contains an error, which makes it impossible to get reliable measurement of the boundary heat flux $g(t)$.
Lab experiments pursued on laser ablation of biological tissues allow for the measure of final temperature distribution and final ablation depth; the ISP must be solved for the identification of $g$.
Our approach allows us to regularize an error contained in the final moment temperature measurement $w(x)$ and final moment ablation depth $s_*$.
Another advantage of this approach is that the condition~\eqref{eq:pde-freebound} can be treated as a measurement of the temperature on the ablation front, allowing us to regularize the error contained in temperature measurement $\mu(t)$ on the ablation front.
Still another important motivation arises from the optimal control of the laser ablation process.
A typical control problem arises when an unknown control parameter, such as the heat flux $g$ on the known boundary must be chosen with the purpose of achieving a desired ablation depth and temperature distribution at the end of the time interval.

ISP is not well posed in the sense of Hadamard: the solution may not exist; if it exists, it may not be unique, and in general it does not exhibit continuous dependence on the data.
The goal of this paper is to pursue numerical analysis of the gradient method in Besov-Sobolev spaces based on the Fr\'echet differential and necessary condition for optimality (\cite{abdulla16,abdulla17}) in the optimal control problem introduced recently as a variational formulation of the \emph{inverse Stefan problem} (ISP) in~\cite{abdulla13,abdulla15}.

The inverse Stefan problem first appeared in~\cite{cannon67}; the problem
discussed was the determination of a heat flux on the fixed boundary for which
the solution of the Stefan problem has a desired free boundary.
The variational approach for solving this ill-posed inverse Stefan problem was developed in~\cite{budak72,budak73,budak74}.
In~\cite{vasilev69}, the problem of finding the optimal value for the external
temperature in order to achieve a given measurement of temperature at the final
moment was considered, and existence was proven.
In~\cite{yurii80}, the Fr\'echet differentiability and
convergence of difference schemes was proven for the same problem, and Tikhonov
regularization was suggested.

Later development of the inverse Stefan problem proceeded along two lines:
inverse Stefan problems with given phase boundaries
in~\cite{bell81,budak74,cannon64,carasso82,ewing79,ewing79a,goldman97,hoffman81,sherman71},
and inverse problems with unknown phase boundaries in~\cite{baumeister80,fasano77,goldman97,hoffman82,hoffman86,jochum80,jochum80a,knabner83,ladyzhenskaya68,lurye75,niezgodka79,nochetto87,primicerio82,sagues82,talenti82}.
We refer to the monograph~\cite{goldman97} for a complete list of references for both types of inverse Stefan problem, both for linear and quasilinear parabolic equations.

The established variational methods in earlier works fail in general to address two issues:
\begin{itemize}
    \item The solution of ISP does not depend continuously on the phase transition temperature.
    A small perturbation of the phase transition temperature may imply significant change of the solution to the ISP.\@
    \item In the existing formulation, at each step of the iterative method a
    Stefan problem must be solved which incurs a high computational cost.
\end{itemize}

A new method developed in~\cite{abdulla13,abdulla15} addresses both issues with a new variational formulation.
Existence of the optimal control and the convergence of the sequence of discrete
optimal control problems to the continuous optimal control problem was proved
in~\cite{abdulla13,abdulla15}.
In~\cite{abdulla16}, the Fr\'echet differentiability and necessary optimality condition in Besov spaces
were established under minimal conditions on the data, when control parameters are chosen as a free boundary $s$, the heat flux $g$,
and the density of sources $f$;
In~\cite{abdulla17} the results are extended to the case when the control vector includes the coefficients
$a,b,c$. A new method for solving optimal control of multiphase Stefan problem is presented in a recent paper \cite{abdullapoggi}.

The structure of the remainder of the paper is as follows: in Section~\ref{sec:notations} we define all
the functional spaces.
Section~\ref{sec:opt_ctrl_problem} formulates optimal control problem. In Section~\ref{sec:discr_conv} we introduce  discrete optimal control problem. Theorem~\ref{thm:convergence} formulates the result on the convergence of the sequence of discrete optimal control problem to the continuous optimal control problem. In Section~\ref{sec:frechet} we introduce the adjoined PDE problem and present the Fr\'echet differentiability result in Theorem~\ref{thm:gradient-result2}. Corollary~\ref{optimalitycondition} presents the necessary condition for the optimal control in the form of the variational inequality. In Section~\ref{sec:grad_besov} we describe the numerical algorithm based on the gradient method in Besov spaces. Section~\ref{sec:numerical} presents the numerical results. Finally, conclusions are presented in Section~\ref{sec:conclusions}

\section{Notations}
\label{sec:notations}

We will use the notation
\[
  \Ind_{I}(x)
  := \begin{cases}
  1,&~ x \in I
  \\
  0,&~ x \not\in I
  \end{cases}
\]
for the indicator function of the set $I$, and $[r]$ for the integer part of the real number $r$.
We will require the notions of Sobolev-Slobodeckij or Besov spaces~\cite{besov79,besov79a,kufner77,nikolskii75}.
In this section, assume $U$ is a domain in ${\mathbb R}$ and denote by
\[
  Q_T := (0,1)\times (0, T].
\]
\begin{itemize}
  \item For $\l \in Z_+$, $W_p^{\l}(U)$ is the Banach space of measurable functions with finite norm
\[
    \norm{u}_{W_p^{\l}(U)}
      := \left(\int_U \sum_{k=0}^{\l} \lnorm{\frac{d^k u}{dx^k}}^p\,dx\right)^{1/p}\nonumber
\]
  \item For $\l >0$, $B_p^{\l}(U)$ is the Banach space of measurable functions with finite norm

\[
    \norm{u}_{B_p^{\l}(U)}
      :=
    \norm{u}_{W_p^{[\l]}(U)}
    + \seminorm{u}_{B_p^{\l}(U)}.
    \]
    If $\l \in Z_+$, the seminorm $\seminorm{u}_{B_p^{\l}(U)}$ is given by
    \[
    \seminorm{u}_{B_p^\l(U)}^p
      :=\int_{U} \int_{U} \frac{\lnorm{
    \D{^{[\l]} u(x)}{x^{[\l]}}
    - \D{^{[\l]} u(y)}{x^{[\l]}}}^p
    }{\lnorm{x-y}^{1 + p(\l - [\l])}} \,dx \,dy,\nonumber
\]
  while if $\l \in Z_+$, $\seminorm{u}_{B_p^{\l}(U)}$ is given by
  \[
    \seminorm{u}_{B_p^\l(U)}^p
      :=
    \int_{-\infty}^\infty \int_{-\infty}^\infty \frac{\lnorm{
    \D{^{\l-1} u(x)}{x^{\l-1}}
    - 2 \D{^{\l-1} u\left(\frac{x+y}{2}\right)}{x^{\l-1}}
    + \D{^{\l-1} u(y)}{x^{\l-1}}
    }^p}{\lnorm{x-y}^{1+p}} \,dy \,dx\nonumber
\]\cite[thm.\ 5, p. 72]{solonnikov64}.
By~\cite[\S 18, thm.\ 9]{besov79a}, it follows that for $p=2$ and $\l \in Z_+$, the $B_p^{\l}(U)$ norm is equivalent to the $W_p^{\l}(U)$ norm (i.e.\ the two spaces coincide.) Common notation $H^\l$ is used instead of $B_2^\l$ or $W_2^\l$ if $\l \in Z_+$.

%
  \item Let $1 \leq p < \infty$, $\l_1,\l_2>0$. The Besov space $B_{p,x,t}^{\l_1, \l_2}(Q_T)$ is defined as the closure of the set of smooth functions under the norm
\[
    \norm{u}_{B_{p,x,t}^{\l_1, \l_2}(Q_T)}
    := \left(\int_0^T \norm{u(x,t)}_{B_p^{\l_1}(0,1)}^p \,dt\right)^{1/p}
\]
\[
    + \left(\int_0^1 \norm{u(x,t)}_{B_p^{\l_2}(0,T)}^p \,dx\right)^{1/p}.
\]
  When $p=2$, if either $\l_1$ or $\l_2$ is an integer, the Besov seminorm may be replaced with the corresponding Sobolev seminorm (and the corresponding space denoted by $W_2^{\l_1, \l_2}$ due to equivalence of the norms.
  \item The H\"older space $C_{x,t}^{\alpha,\alpha/2}(Q_T)$ is the set of continuous functions with $[\alpha]$ $x$-derivatives and $[\alpha/2]$ $t$-derivatives, and for which the highest order $x$- and $t$-derivatives satisfy H\"older conditions of order $\alpha-[\alpha]$ and $\alpha/2-[\alpha/2]$, respectively.
  
   \item $V_{2}(\Omega)$ is the subspace of $B_{2}^{1,0}(\Omega)$ for which the norm
  \[
    \norm{u}_{V_{2}(\Omega)}^2
    = \esssup_{0\leq t \leq T} \norm{u(\cdot, t)}_{L_{2}\big(0,s(t)
    \big)}^{2}
    + \lnorm{\D{u}{x}}_{L_{2}(\Omega)}^{2}
  \]

  \item $V_{2}^{1,0}(\Omega)$ is the completion of $B_{2}^{1,1}(\Omega)$ in the $V_{2}(\Omega)$ norm.
  For $u \in V_{2}^{1,0}(\Omega)$, the function
  \[
    \phi(t)=\norm{u(\cdot,t)}_{L_{2}(0,s(t))}
  \]
  varies continuously.
  $V_{2}^{1,0}(\Omega)$ is a Banach space with norm
  \[
    \norm{u}_{V_{2}^{1,0}(\Omega)}^2
    = \max_{0\leq t \leq T} \norm{u(\cdot, t)}_{L_{2}\big(0,s(t)\big)}^{2}
    + \lnorm{\D{u}{x}}_{L_{2}(\Omega)}^{2}
  \]

\end{itemize}

\section{Optimal Control Problem}
\label{sec:opt_ctrl_problem}

Consider a minimization of the cost functional
\begin{gather}
\mathcal{J}(v)=\beta_{0}\Vert u(x,T)-w(x)\Vert_{L_{2}[0,s(T)]}^{2}+\beta_{1}\Vert u(s(t),t)-\mu(t)\Vert_{L_{2}[0,T]}^2 \nonumber\\
+\beta_2 |s(T)-s_*|^2\label{Eq:W:1:8}
\end{gather}
on the control set
\begin{gather*}
V_{R}=\{ v=(s,g) \in B_{2}^{2}[0,T]\times B_{2}^{1}[0,T]: \delta \leq s(t)\leq l, s(0)=s_0, s'(0)=0, \nonumber\\
\max (~\Vert s\Vert_{B_{2}^{2}}; ~\Vert g\Vert_{B_{2}^{1}} \leq R\}
\end{gather*}
where $\delta, l,R, \beta_0, \beta_1$ are given positive numbers, and $u=u(x,t;v)$ be a solution of the Neumann problem \eqref{eq:pde-1}--\eqref{eq:pde-stefan}.

\begin{definition}
The function $u \in W_{2}^{1,1}(\Omega)$ is called a weak solution of the problem \eqref{eq:pde-1}--\eqref{eq:pde-stefan} if $u(x,0)=\phi(x) \in W_{2}^{1}[0,s_0]$ and
\begin{gather}
0=\int_{0}^{T}\int_{0}^{s(t)}[ a u_{x}\Phi_{x}-bu_{x}\Phi - c u \Phi + u_{t} \Phi+f\Phi] \,dx\,dt \nonumber\\
 +\int_{0}^{T}[ \gamma(s(t),t)s'(t)-\chi(s(t),t)]\Phi(s(t),t)\, dt
+\int_{0}^{T}g(t)\Phi(0,t)\, dt\label{Eq:W:1:9}
\end{gather}
for arbitrary $\Phi \in W_{2}^{1,1}(\Omega)$
\end{definition}

Furthermore, formulated optimal control problem will be called Problem $I$.

\subsection{Discretization and convergence}
\label{sec:discr_conv}
Let
\[\omega_{\tau}=\{ t_{j}=j \cdot \tau,~j=0,1,\ldots,n\} \]
be a grid on $[0,T]$ and $\tau=\frac{T}{n}$. Consider a discretized control set
\begin{equation*}
V^n_{R}=\{ [v]_{n}=([s]_n,[g]_n) \in {\mathbb R}^{2n+2}:~0<\delta\leq s_{k} \leq l,~ \max(\Vert [s]_{n}\Vert_{w_{2}^{2}}^2; ~\Vert [g]_{n}\Vert_{w_{2}^{1}}^2) \leq R^2\}
\end{equation*}
where,
\[ [s]_n=(s_0,s_1,...,s_n) \in {\mathbb R}^{n+1}, \ [g]_n=(g_0,g_1,...,g_n) \in {\mathbb R}^{n+1} \]
\[
\Vert [s]_{n}\Vert_{w_{2}^{2}}^2= \sum\limits_{k=0}^{n-1}\tau s_k^2+\sum\limits_{k=1}^{n}\tau s_{\overline{t},k}^2+\sum\limits_{k=1}^{n-1}\tau s_{\overline{t}t,k}^2, \  \Vert [g]_{n}\Vert_{w_{2}^{1}}^2= \sum\limits_{k=0}^{n-1}\tau g_k^2+\sum\limits_{k=1}^{n}\tau g_{\overline{t},k}^2.
\]
under the standard notation for the finite differences:
\[ s_{\overline{t},k}=\frac{s_k-s_{k-1}}{\tau}, \ s_{t,k}=\frac{s_{k+1}-s_{k}}{\tau}, \ s_{\overline{t}t,k}^2=\frac{s_{k+1}-2s_k+s_{k-1}}{\tau^2}. \]
Introduce two mappings $\mathcal{Q}_n$ and $\mathcal{P}_n$ between continuous and discrete control sets:
\[ \mathcal{Q}_n(v)=[v]_n=([s]_n,[g]_n), \quad \text{for}~ v\in V_R \]
where $s_k=s(t_k), g_k=g(t_k), k=0,1,...,n$.
\[ \mathcal{P}_n([v]_n)=v^n=(s^n,g^n)\in W_2^2[0,T]\times W_2^1[0,T] \quad \text{for}~ [v]_n \in V_R^n, \]
where
\begin{equation}\label{Eq:W:1:11}
s^n(t)=
\left\{
\begin{array}{l}
s_0+\frac{t^2}{2\tau} s_{\overline{t},1} \ \ 0\le t \le \tau,\\
s_{k-1}+(t-t_{k-1}-\frac{\tau}{2})s_{\overline{t},k-1}+\frac{1}{2}(t-t_{k-1})^2 s_{\overline{t}t,k-1} \ \ t_{k-1}\le t \le t_k, k=\overline{2,n}.
\end{array}\right.
\end{equation}
\begin{equation*}
g^n(t)=g_{k-1}+\frac{g_k-g_{k-1}}{\tau}(t-t_{k-1}), \ \ t_{k-1} \le t \le t_k, k=\overline{1,n}.
\end{equation*}
Let us now introduce a spatial grid. Let $[v]_n\in V_R^n$, let $(p_0,p_1,\cdots,p_n)$ be a permutation of $(0,1,\cdots,n)$ according to order
\[ s_{p_0}\le s_{p_1}\le \cdots \le s_{p_n} \]
In particular, according to this permutation for arbitrary $k$ there exists a unique $j_k$ such that
\begin{equation}\label{Eq:W:1:11a}
s_k=s_{p_{j_k}}
\end{equation}
Furthermore, unless it is necessary in the context, we are going to write simply $j$ instead of subscript $j_k$. Let
\[\omega_{p_0}=\{ x_{i}: x_i=i \cdot h,~i=0,1,\ldots,m_0^{(n)}\} \]
be a grid on $[0,s_{p_0}]$ and $h=\frac{s_{p_0}}{m_0^{(n)}}$. Furthermore, we always assume that
\begin{equation}\label{htau}
h = O(\sqrt{\tau}), \quad \text{as}~ \tau \rightarrow 0.
\end{equation}
We continue construction of the spatial grid by induction. Having constructed $\omega_{p_{k-1}}$ on $[0,s_{p_{k-1}}]$ we construct
\[ \omega_{p_k}=\{ x_i:~i=0,1,\cdots, m_k^{(n)} \} \]
on $[0,s_{p_{k}}]$, where $m_k^{(n)}\ge m_{k-1}^{n}$, and inequality is strict if and only if $s_{p_{k}}>s_{p_{k-1}}$; for $i\le m_{k-1}^{(n)}$ points $x_i$ are the same as in grid $\omega_{p_{k-1}}$.
Finally, if $s_{p_{n}}<l$, then we introduce a grid on $[s_{p_n},l]$
\[ \overline{\omega}=\{x_i: x_i=s_{p_n}+(i-m_n^{(n)}) \overline{h}, ~i=m_n^{(n)},\cdots, N \} \]
of stepsize order $h$, i.e. $\overline{h}=O(h)$ as $h \rightarrow 0$.
Furthermore we simplify the notation and write $m_k^{(n)}\equiv m_k$. Let
\[ h_i=x_{i+1}-x_i, \ i=0,1,\cdots,N-1; \]
and assume that
\begin{equation*}
m_k \rightarrow +\infty, \quad \text{as}~ n\rightarrow \infty.
\end{equation*}
Introduce Steklov averages
\[ d_{k}(x)=\frac{1}{\tau}\int_{t_{k-1}}^{t_{k}}d(x,t)\,dt, \ h_{k}=\frac{1}{\tau}\int_{t_{k-1}}^{t_{k}}h(t)\,dt, \ d_{ik}=\frac{1}{h_i \tau} \int_{x_i}^{x_{i+1}}\int_{t_{k-1}}^{t_k} d(x,t)\,dt\,dx, \]
where $i=0,1,\cdots,N-1; \ k=1,\cdots,n;$ $d$ stands for any of the functions $a$, $b$, $c$, $f$, and $h$ stands for any of the functions $\nu$, $\mu$, $g$ or $g^n$. Given $v=(s,g) \in V_R$ we define Steklov averages of traces
\begin{equation}\label{Eq:W:1:12}
\chi^{k}_s=\frac{1}{\tau} \int_{t_{k-1}}^{t_{k}}\chi(s(t),t) \,dt, \
(\gamma_s s')^k=\frac{1}{\tau} \int_{t_{k-1}}^{t_{k}}\gamma(s(t),t)s'(t) \,dt.
\end{equation}
Given $[v]_n=([s]_n,[g]_n) \in V_R^n$ we define Steklov averages $\chi^{k}_{s^n}$ and $(\gamma_{s^n} (s^n)')^k$ through \eqref{Eq:W:1:12} with  $s$ replaced by $s^n$ from \eqref{Eq:W:1:11}.

Next we define a discrete state vector through discretization of the integral identity \eqref{Eq:W:1:9}
\begin{definition}\label{discretestatevector}
Given discrete control vector $[v]_n$, the vector function
\[ [u([v]_n)]_n=(u(0),u(1),...,u(n)), \ u(k)\in {\mathbb R}^{N+1}, \ k=0,\cdots,n  \]
is called a discrete state vector if
\begin{description}
\item{\bf(a)} First $m_0+1$ components of the vector $u(0)\in {\mathbb R}^{N+1}$ satisfy
\[ u_i(0)=\phi_i := \phi(x_i), \ i=0,1,\cdots,m_0; \]
\item{\bf(b)} Recalling \eqref{Eq:W:1:11a}, for arbitrary $k=1,\cdots,n$, the first $m_j+1$ components of the vector $u(k)\in {\mathbb R}^{N+1}$ solve the following system of $m_j+1$ linear algebraic equations:
\begin{gather}
\Big [ a_{0k}+hb_{0k}-h^2c_{0k}+\frac{h^2}{\tau} \Big ] u_0(k) - \Big [ a_{0k}+hb_{0k} \Big ] u_1(k)=\frac{h^2}{\tau}u_0(k-1)-h^2f_{0k}-hg^n_{k}, \nonumber\\
-a_{i-1,k}h_iu_{i-1}(k)+\Big [ a_{i-1,k}h_i+a_{ik}h_{i-1}+b_{ik}h_ih_{i-1}-c_{ik}h_i^2h_{i-1}+\frac{h_i^2h_{i-1}}{\tau} \Big ] u_i(k)- \nonumber\\
\Big [ a_{ik}h_{i-1}+b_{ik}h_ih_{i-1} \Big ] u_{i+1}(k) = -h_i^2h_{i-1}f_{ik}+\frac{h_i^2h_{i-1}}{\tau}u_i(k-1), \ i=1,\cdots,m_j-1  \nonumber\\
-a_{m_j-1,k} u_{m_j-1}(k)+a_{m_j-1,k} u_{m_j}(k)=-h_{m_j-1} \Big [  (\gamma_{s^n} (s^n)')^k-\chi^{k}_{s^n} \Big ].\label{alma}
\end{gather}
\item{\bf(c)} For arbitrary $k=0,1,...,n$, the remaining components of $u(k)\in {\mathbb R}^{N+1}$ are calculated as
\[ u_i(k)= \hat{u}(x_i;k), \ m_j\le i \le N  \]
where $\hat{u}(x;k) \in W_2^1[0,l]$ is a piecewise linear interpolation of $\{u_i(k): \ i=0,\cdots,m_j \}$, that is to say
\[  \hat{u}(x;k)=u_i(k)+\frac{u_{i+1}(k)-u_i(k)}{h_i} (x-x_i), \ x_i\le x\le x_{i+1}, i=0,\cdots,m_j-1,  \]
iteratively continued to $[0,l]$ as
\begin{equation}\label{Eq:W:1:14}
\hat{u}(x;k)=\hat{u}(2^ns_k-x;k), \ 2^{n-1}s_k \le x \le 2^ns_k, n=\overline{1,n_k}, \ n_k\le n_*=1+\log_2\Big [ \frac{l}{\delta}\Big ]
\end{equation}
where $[r]$ means integer part of the real number $r$.
\end{description}
\end{definition}

It should be mentioned that for any $k=1,2,\cdots,n$, system \eqref{alma} is equivalent to the following summation identity
\begin{gather}
\sum_{i=0}^{m_j-1}h_i \Big [ a_{ik}u_{ix}(k)\eta_{ix}-b_{ik}u_{ix}(k)\eta_i-c_{ik}u_i(k)\eta_i+
f_{ik}\eta_i+u_{i\overline{t}}(k)\eta_i \Big ] + \nonumber\\ \Big [ (\gamma_{s^n} (s^n)')^k-\chi^{k}_{s^n} \Big ]\eta_{m_j}+g^n_k \eta_0=0,\label{Eq:W:1:19}
\end{gather}
for arbitrary numbers $\eta_i, i=0,1,\cdots,m_j$.

Consider a discrete optimal control problem of minimization of the cost functional
\begin{equation}\label{Eq:W:1:15}
 \mathcal{I}_n([v]_n)
  =\beta_{0}\sum_{i=0}^{m_n-1} h_i \Big( u_i(n) - w_i \Big)^2
  +\beta_{1}\tau\sum_{k=1}^n \Big( u_{m_k}(k)-\mu_k \Big)^2
  +\beta_2 \lnorm{s_n - s_*}^2
\end{equation}
on a set $V_R^n$ subject to the state vector defined in Definition 1.3. Furthermore, formulated discrete optimal control problem will be called Problem $I_n$.

Throughout, we use piecewise constant and piecewise linear interpolations of the discrete state vector:
given discrete state vector $[u([v]_n)]_n=(u(0),u(1),...,u(n))$, let
\[ u^\tau(x,t)=\hat{u}(x;k), \quad \text{if}~ t_{k-1}<t\le t_k, \ 0\le x \le l, \ k=\overline{0,n}, \]
\[ \hat{u}^\tau(x,t)=\hat{u}(x;k-1)+\hat{u}_{\overline{t}}(x;k)(t-t_{k-1}), \quad \text{if}~ t_{k-1}<t\le t_k, \ 0\le x \le l, \ k=\overline{1,n}, \]
\[ \hat{u}^\tau(x,t)= \hat{u}(x;n), \quad \text{if}~ t\ge T, \ 0\le x \le l. \]
\[ \tilde{u}^\tau(x,t)=u_i(k), \quad \text{if}~ t_{k-1}< t \le t_k, \ x_i \le x < x_{i+1}, \ k=\overline{1,n}, \ i=\overline{0,N-1}.\]
Obviously, we have
\[ u^\tau \in V_2(D), \ \ \hat{u}^\tau \in W_2^{1,1}(D), \ \ \tilde{u}^\tau \in L_2(D).  \]
As before, we employ standard notations for difference quotients of the discrete state vector:
\[ u_{ix}(k)=\frac{u_{i+1}(k)-u_i(k)}{h_i}, \ u_{i\overline{t}}=\frac{u_i(k)-u_i(k-1)}{\tau}, \ \quad \text{etc.}  \]

Assume that the following assumptions are satisfied:
\begin{gather}
  a,~\D{a}{x},b,c \in L_\infty(D); \ a\geq a_0>0
  ~\text{a.e.}~\text{in}~D;~\int_{0}^T \esssup_{0 \leq
    x \leq \l} \lnorm{\D{a}{t}} \,dt < +\infty\nonumber\\
  w \in L_2(0,\l), f \in L_2(D),
  ~\chi,\gamma \in B_2^{1,1}(D),
  ~\phi \in B_2^1(0,s_0),
  ~\mu \in L_2(0,T)\nonumber
\end{gather}

The following results characterize the convergence of the sequence of discrete optimal control problems to the continuous optimal control problem.

\begin{thm}\label{thm:convergence}\cite{abdulla15}
The sequence of discrete optimal control problems $I_n$ approximates the optimal control problem $I$ with respect to the functional, i.e.
\begin{equation}\label{Eq:W:1:18}
\lim\limits_{n\to +\infty} \mathcal{I}_{n_*}=\mathcal{J}_*,
\end{equation}
where
\[ \mathcal{I}_{n_*}=\inf\limits_{V_R^n} \mathcal{I}_n([v]_n), \ n=1,2,... \]
If $[v]_{n_\epsilon}\in V_R^n$ is chosen such that
\[ \mathcal{I}_{n_*} \le \mathcal{I}_n([v]_{n_\epsilon})\le \mathcal{I}_{n_*}+\epsilon_n, \ \epsilon_n \downarrow 0, \]
then the sequence $v_n=(s_n,g_n)=\mathcal{P}_n([v]_{n_\epsilon})$ converges to some element $v_*=(s_*,g_*) \in V_*$ weakly
in $W_2^2[0,T] \times W_2^1[0,T]$, and strongly in $W_2^1[0,T] \times L_2[0,T]$. In particular $s_n$ converges to $s_*$ uniformly on $[0,T]$. For any $\delta>0$, define
  \[
    \Omega_*^{\prime} = \Omega_* \cap \{x<s_*(t)-\delta,~0<t<T\}
  \] Then the piecewise linear interpolation $\hat{u}^\tau$ of the
discrete state vector $[u[v]_{n_\epsilon}]_n$ converges to the solution $u(x,t;v_*) \in W_2^{1,1}(\Omega_*)$ of the Neumann
problem \eqref{eq:pde-1}--\eqref{eq:pde-stefan} weakly in $W_2^{1,1}(\Omega_*^{\prime})$.
\end{thm}
\begin{rmk}\label{comparisonwithIPI}
The only difference between Problems $I$, $I_n$  and the corresponding optimal control problems in \cite{abdulla15} is that the cost functionals \eqref{Eq:W:1:8} and \eqref{Eq:W:1:15} are replaced respectively with
\begin{gather}
\mathcal{J}(v)=\beta_{0}\Vert u(0,t)-\nu(t)\Vert_{L_{2}[0,T]}^{2}+\beta_{1}\Vert u(s(t),t)-\mu(t)\Vert_{L_{2}[0,T]}^2 \nonumber\\
\mathcal{I}_n([v]_n)=\beta_{0}\tau\sum\limits_{k=1}^n \Big ( u_0(k)-\nu_k \Big )^2+\beta_{1}\tau\sum\limits_{k=1}^n \Big ( u_{m_k}(k)-\mu_k \Big )^2.\nonumber
\end{gather}
\end{rmk}
The proof of Theorem~\ref{thm:convergence} is almost identical to the proof of the corresponding convergence theorem of \cite{abdulla15}.

\subsection{Fr\'echet differentiability in Besov spaces and optimality condition}
\label{sec:frechet}
Fr\'echet differentiability of the cost functional $\mathcal{J}(v)$ is true under slightly higher regularity assumptions on the data. Let $\alpha>0$ be fixed, $H := \scontrolspace(0, T)\times \gcontrolspace(0, T)$ and
\begin{align}
  V^1_R=\Big\{ v=(s, g) \in H: s(0) = s_0,
      ~s'(0) = 0,
  ~g(0) = a(0, 0) \phi'(0),\nonumber\\
  ~0 < \delta \leq s(t), \norm{v}_H := \max\left( \norm{s}_{\scontrolspace(0,T)}, \norm{g}_{\gcontrolspace(0,T)}, \right)\leq R
  \Big\},\label{eq:control-set}
\end{align}
In addition to the assumptions formulated in Section~\ref{sec:discr_conv} we assume that
\begin{gather*}
  a,a_x,b,c \in C_{x, t}^{1/2 + 2\alpha^*, 1/4 + \alpha^*}(D),\label{eq:datacond-coeff},~w \in W_2^1(0, \l),\quad \phi \in \ivspace(0, s_{0}),
  \\
  \chi,\gamma \in \chigammaspace(D), \mu \in \muspace(0,T), f \in \fcontrolspace(D)
\end{gather*}
where $\alpha^* > \alpha$ is arbitrary, and $\chi$, $\phi$ satisfy the compatibility condition
\begin{equation*}
  \chi(s_{0},0) = \phi'(s_{0})a(s_{0},0).
\end{equation*}
Given a control vector $v \in V_R$, under this
conditions there
exists a unique pointwise a.e.\ solution $u \in W_2^{2,1}(\Omega)$ of the
Neumann problem~\eqref{eq:pde-1}--\eqref{eq:pde-stefan} (\cite{ladyzhenskaya68,solonnikov64}).
\begin{definition}\label{defn:adjoint}
  For given $v$ and $u = u(x, t; v)$, $\psi \in \adjointsolnspace(\Omega)$ is a solution to the adjoint problem if
  \begin{gather}
    L^* \psi := {\big(a\psi_x\big)}_x - {(b\psi)}_x + c\psi + \psi_t = 0,\quad\mathrm{in}~\Omega\label{eq:adj-pde}
    \\
    \psi(x, T) = 2\beta_0(u(x, T) - w(x)),~0 \leq x \leq s(T)\label{eq:adj-finalmoment}
    \\
    a(0, t)\psi_x(0, t) - b(0, t)\psi(0, t)=0,~0 \leq t \leq T\label{eq:adj-robin-fixed}
    \\
    {\Big[a\psi_x - (b + s'(t))\psi\Big]}_{x=s(t)} = 2\beta_1(u(s(t), t) - \mu(t)), ~0 \leq t \leq T\label{eq:adj-robin-free}
  \end{gather}
\end{definition}
Given a control vector $v \in V_R$ and the corresponding state vector $u \in
W_2^{2,1}(\Omega)$, there exists a unique pointwise a.e.\ solution $\psi \in
W_2^{2,1}(\Omega)$ of the adjoint problem~\eqref{eq:adj-pde}--\eqref{eq:adj-robin-free} \cite{ladyzhenskaya68,solonnikov64}.

The following theorem formulates the Fr\'echet differentiability of the cost functional $J(v)$(\cite{abdulla16}):
\begin{thm}[Fr\'echet Differentiability]\label{thm:gradient-result2}\cite{abdulla16}
  The functional $\J(v)$ is differentiable in the sense of Fr\'echet, and the Fr\'echet differential is
  \begin{gather}
    \left\langle{}\J'(v),{\delta v}\right\rangle_H
      = - \int_{0}^T \psi(0,t) {\delta g}(t)\,dt-\int_0^T \big[\gamma \psi\big]_{x=s(t)}{\delta s}'(t)\,dt \nonumber
    \\
      + \int_0^T \left[2 \beta_1 (u-\mu) u_x
    + \psi \left(\chi_x - \gamma_x s' -\big(a u_x\big)_x\right)\right]_{x=s(t)}{\delta s}(t)\, dt  \nonumber
       \\
      + \left(\beta_0\lnorm{u(s(T),T) - w(s(T))}^2 + 2\beta_2(s(T)-s_*)\right){\delta s}(T),\label{eq:gradient-full}
  \end{gather}
  where $\J'(v)\in H'$ is the Fr\'echet derivative, $\langle{} \cdot,\cdot\rangle{}_H$ is a pairing
  between $H$ and its dual $H'$,
$\psi$ is a solution to the adjoint problem in the sense of Definition~\ref{defn:adjoint}, and ${\delta v} = ({\delta s},{\delta g})$ is a variation of the control vector $v \in V^1_R$ such that $v + \delta v \in V^1_R$.
\end{thm}
\begin{cor}[Optimality Condition]\label{optimalitycondition}
  If $\mathbf{v} = (\mathbf{s}, \mathbf{g})$ is an optimal control,
  then the following variational inequality is satisfied:
  \begin{equation}
    \left\langle \J'(\mathbf{v}), v-\mathbf{v} \right\rangle_H\geq 0\label{eq:optimality-condition}
  \end{equation}
  for arbitrary $v = (s, g) \in V_R$.
\end{cor}

\subsection{Gradient method in Besov spaces}
\label{sec:grad_besov}
Fr\'echet differentiability result of Theorem~\ref{thm:gradient-result2} and the formula~\eqref{eq:gradient-full}
for the Fr\'echet differential suggest the following algorithm based on the projective gradient method:
\begin{description}
\item[{\bf Step 1.}] Set $k=0$ and choose initial vector function $v_0=(s_0, g_0) \in V_R$.

\item[{\bf Step 2.}] Solve the Neumann problem~\eqref{eq:pde-1}--\eqref{eq:pde-stefan} to find $u_k=u(x,t;v_k)$ and $\J(v_k)$.

\item[{\bf Step 3.}] If $k=0$, move to Step 4.
  Otherwise, check the following criteria:
  \begin{equation}
    \left| \frac{\J(v_{k})-\J(v_{k-1})}{\J(v_{k-1})} \right| <\epsilon, \quad
   \frac{ \norm{v_{k}-v_{k-1}}}{\norm{v_{k-1}}} < \epsilon, \label{convergencecriteria}
  \end{equation}
  where $\epsilon$ is the required accuracy.
  If the criteria are satisfied, then terminate the iteration. Otherwise, move to Step 4.
\item[{\bf Step 4.}] Having $u_k$, solve the adjoined PDE problem~\eqref{eq:adj-pde}--\eqref{eq:adj-robin-free} to find $\psi_k=\psi(x,t;v_k)$.
\item[{\bf Step 5.}] Choose stepsize parameter $\alpha_k>0$ and compute new control vector
  $v_{k+1}=(s_{k+1}, g_{k+1})\in H$ as follows:
\begin{gather}
  g_{k+1}(t)=g_k(t)+\alpha_k
  \psi_k(0,t),\label{gradientupdate_g}
  \\
  s_{k+1}(t)=s_k(t)-\alpha_k \Big[2 \beta_1 (u_k-\mu) u_{kx}
   \nonumber
  \\
  +\psi_k \left(\chi_x - \gamma_x s_k' -\big(a
  u_{kx}\big)_x\right)\Big]_{x=s_k(t)},\label{gradientupdate_s}
  \\
  s'_{k+1}(t)=s'_k(t)+\alpha_k \big[\gamma
    \psi_k\big]_{x=s_k(t)},\label{gradientupdate_s'}
  \\
  s_{k+1}(T)=s_k(T)-\alpha_k [\beta_0\lnorm{u_k(s_k(T),T) -
    w(s_k(T))}^2 + \nonumber
  \\
  2\beta_2(s_k(T)-s_*)].
  \label{gradientupdate_s(T)}
\end{gather}
\item[{\bf Step 6.}]
Replace $v_{k+1}$ with $\mathcal{P}_{V_R}(v_{k+1})$, where $\mathcal{P}_{V_R}:H \to V^1_R$ is the projection operator to the closed and convex subset $V^1_R$.
Then replace $k$ with $k+1$ and move to Step 2.
\end{description}
Note that the construction of the component $s_{k+1}\in W_2^2(0,T)$ is achieved through interpolation using the values of $s_{k+1}$ and $s_{k+1}'$ at grid points $t=t_k$ according to the formulae~\eqref{gradientupdate_s} and \eqref{gradientupdate_s'}. Moreover, $s_{k+1}(T)$ is updated according to~\eqref{gradientupdate_s(T)}. In practical applications, the fact that $s$ and $s'$ are updated independently causes some inconvenience, and an alternative algorithm where only $s$ is updated would be preferred. By slight increase of the regularity assumption on $\gamma$ (precisely $\gamma \in B^{2,1}_{2,x,t}(D)$), one can transform \eqref{eq:gradient-full} to the alternative form:
\begin{gather}
    \left\langle{}\J'(v),{\delta v}\right\rangle_H
      = - \int_{0}^T \psi(0,t) {\delta g}(t)\,dt\nonumber
    \\
    + \int_0^T \left[2 \beta_1 (u-\mu) u_x
    + \psi \left(\chi_x+\gamma_t\right) + \gamma \psi_x s' +\psi_t \gamma-\psi \big(a u_x\big)_x\right]_{x=s(t)}{\delta s}(t)\, dt  \nonumber
       \\
      + \left(\beta_0\lnorm{u(s(T),T) - w(s(T))}^2 + 2\beta_2(s(T)-s_*)-\gamma \psi |_{(s(T),T)}\right){\delta s}(T).\label{eq:gradient-full-alternative}
  \end{gather}
This suggests a modification of the described above algorithm where \eqref{gradientupdate_s}--\eqref{gradientupdate_s(T)}
are replaced with
\begin{gather}
s_{k+1}(t)=s_k(t)-\alpha_k \Big[2 \beta_1 (u_k-\mu) u_{kx}+\psi_k \left(\chi_x+\gamma_t\right)
   \nonumber
  \\
  +\gamma \psi_{kx} s_k' +\psi_{kt} \gamma-\psi_k\left(a
  u_{kx} \right)_x \Big]_{x=s_k(t)},\label{gradientupdate_s_modified}
  \\
  s_{k+1}(T)=s_k(T)-\alpha_k [\beta_0\lnorm{u_k(s_k(T),T) -
    w(s_k(T))}^2\nonumber\\
     + 2\beta_2(s_k(T)-s_*)-\gamma \psi_k |_{(s_k(T),T)}].
     \label{gradientupdate_s(T)_modified}
\end{gather}
\begin{rmk}
\label{rmk:2}
From \eqref{eq:gradient-full-alternative} it follows that the Fr\'echet gradient with respect to $s$ is
\begin{gather}
\J_s'(v)= \left[2 \beta_1 (u-\mu) u_x
    + \psi \left(\chi_x+\gamma_t\right) + \gamma \psi_x s' +\psi_t \gamma-\psi \big(a u_x\big)_x\right]_{x=s(t)}\nonumber\\
    + \left(\beta_0\lnorm{u(s(T),T) - w(s(T))}^2 + 2\beta_2(s(T)-s_*)-\gamma \psi |_{(s(T),T)}\right) \delta_T,
    \label{s-gradient}
    \end{gather}
    where $\delta_T$ is a Dirac measure on $[0,T]$ with support at $t=T$. Fr\'echet gradient with respect to 
    $g$ is
    \begin{equation}\label{g-gradient}
    \J_g'(v)=-\psi(0,t).
    \end{equation}
\end{rmk}
\begin{rmk}
\label{rmk:3}
We will implement Tikhonov regularization by replacing the cost functional \eqref{Eq:W:1:8} with
\begin{gather}
\mathcal{J}(v)=\beta_{0}\Vert u(x,T)-w(x)\Vert_{L_{2}[0,s(T)]}^{2}+\beta_{1}\Vert u(s(t),t)-\mu(t)\Vert_{L_{2}[0,T]}^2 \nonumber\\
+\beta_2 |s(T)-s_*|^2+\beta\Vert s-\bar{s}\Vert_{L_2[0,T]}^2\label{Eq:W:1:8reg}
\end{gather}
where $\beta>0$ is a regularization parameter. In this case instead of \eqref{eq:gradient-full-alternative} one can derive the following expression for the Fr\'echet differential:
\begin{gather}
    \left\langle{}\J'(v),{\delta v}\right\rangle_H
      = - \int_{0}^T \psi(0,t) {\delta g}(t)\,dt
    + \int_0^T [2 \beta_1 (u-\mu) u_x
    + \psi \left(\chi_x+\gamma_t\right)\nonumber\\ + \gamma \psi_x s' +\psi_t \gamma-\psi \big(a u_x\big)_x+2\beta(s(t)-\bar{s}(t))]\delta s(t)\, dt  \nonumber
       \\
      + \left(\beta_0\lnorm{u(s(T),T) - w(s(T))}^2 + 2\beta_2(s(T)-s_*)-\gamma \psi |_{(s(T),T)}\right){\delta s}(T).\label{eq:gradient-full-alternative-reg}
  \end{gather}
 Therefore, the Fr\'echet gradient with respect to $s$ is
\begin{gather}
\J_s'(v)= \left[2 \beta_1 (u-\mu) u_x
    + \psi \left(\chi_x+\gamma_t\right) + \gamma \psi_x s' +\psi_t \gamma-\psi \big(a u_x\big)_x\right]_{x=s(t)}+2\beta(s-\bar{s})\nonumber\\
    + \left(\beta_0\lnorm{u(s(T),T) - w(s(T))}^2 + 2\beta_2(s(T)-s_*)-\gamma \psi |_{(s(T),T)}\right) \delta_T.
    \label{s-gradient-reg}
    \end{gather}
      \end{rmk}

\section{Numerical Results}
\label{sec:numerical}

In this section we provide the computational results obtained to solve the inverse Stefan problem
\eqref{eq:pde-1}--\eqref{eq:pde-finaltemp} by finding an optimal control vector
$\mathbf{v} = (\mathbf{s}, \mathbf{g})$ based on the algorithm described in detail in
Section~\ref{sec:grad_besov}. First, we briefly discuss the numerical approaches used for
discretizing the problem both in space and time, as well as the numerical optimization techniques added
to our computational algorithm to improve its performance. Then we describe the models
chosen to represent various levels of complexity and, finally, we show the outcomes of applying the proposed
computational algorithm to these models.

\subsection{Numerical optimization for discretized models}
\label{sec:approaches}

Our computational approach to solve the inverse Stefan problem \eqref{eq:pde-1}--\eqref{eq:pde-finaltemp}
is formulated in the ``optimize--then--discretize'' framework. Following this paradigm, we formulate
this problem as an optimization problem which in its turn is ultimately discretized for the purpose of
a numerical solution. On the other hand, our optimality conditions and the cost functional gradients
are derived in the continuous, i.e.~PDE setting. As a consequence, the main constituents of the proposed
approach are left independent of the specific discretization used for space and time.

Note that the Frechet gradient $ \J'(v)$ is an element of the dual space $H^\prime$:
\[  \J'(v) = (\J_s'(v), \J_g'(v)) \in H^\prime\]
According to formulae \eqref{s-gradient} (or \eqref{s-gradient-reg}), \eqref{g-gradient}, the $s$-gradient is the sum of elements of $L_2(0,T)$ and a constant multiple of the Dirac measure $\delta_T$, while $g$-gradient is an element of $L_2(0,T)$. Due to the lack of a satisfactory regularity gradient formula \eqref{s-gradient} (or \eqref{s-gradient-reg}), \eqref{g-gradient} may not be suitable for the reconstruction of $v = (s,g)$ \cite{Bukshtynov11,Bukshtynov13}. 
Therefore, for the numerical implementation of the gradient method, we are going to derive an equivalent formula for the gradient with higher regularity. The idea is based on Riesz representation theorem \cite{Berger77}, which expresses the isometrical isomorphism between a Hilbert space and its dual space if the underlying field is the real numbers. Our aim
is to derive an equivalent formula for the Frechet gradient which is the element of the real Hilbert space $H^1(0,T)\times H^1(0,T)$. Moreover, we assume that instead of standard norm, the Hilbert space $H^1(0,T)$ is equipped with the equivalent inner product and norm
\[ (u,w)_{H^1}=\int_0^T uw+\ell^2\frac{du}{dt}\frac{dw}{dt} dt, \ \norm{u}_{H^1}=(u,u)_{H^1}^{\frac{1}{2}} \] 
where $\ell \in \RR^+$ is a ``time-scale" parameter with the purpose to improve the convergence of the gradient method.
To pursue this idea we introduce a notation 
\[ \bnabla^{\mathcal{X}}_v \J=(\bnabla^{\mathcal{X}}_s \J, \bnabla^{\mathcal{X}}_g \J)\]
to represent the Frechet gradient of the functional $\J$ in Hilbert space $\mathcal{X}$. With slight abuse of notation, we are going to use the same notation for the finite-dimensional vector obtained through discretization of the Frechet gradient exclusively for our numerical computations. 

With the refined notation in hand, we can rewrite \eqref{s-gradient}, \eqref{g-gradient} as follows:
\begin{equation}\label{eq:grads0}
\J_s'(v)=\bnabla^{L^2}_s \J(v) + \bnabla^{H^1}_{s(T)} \J(v), \ \  \J_g'(v)=\bnabla^{L^2}_s \J(v)
\end{equation}
where
\begin{equation}
       \begin{aligned}
        \bnabla^{L^2}_{s} \mJ(v) &=
       \left[ 2 \beta_1 (u-\mu) u_{x} + \psi\left( \chi_x + \gamma_t \right)
        +\gamma \psi_{x} s' + \psi_{t} \gamma - \psi \left( a u_{x} \right)_x \right]_{x=s(t)},\\
        \bnabla^{H^1}_{s(T)} \mJ(v) &=
       \left(\beta_0 \lnorm{u(s(T),T) - w(s(T))}^2 + 2\beta_2 (s(T)-s_*) - \gamma \psi |_{(s(T),T)}\right) \delta_T,\\
        \bnabla^{L^2}_g \mJ(v) &= - \psi(0,t)
        \label{eq:grads}
       \end{aligned}
      \end{equation}
From the Riesz representation theorem \cite{Berger77} it follows that the Frechet gradient $\bnabla_v^{H^1} \mJ \in H^1$ satisfies the identity         
\begin{equation}
   \left\langle{}\J'(v),{\delta v}\right\rangle_H =
  \Big(\bnabla_v^{H^1} \mJ, \delta v \Big) _{H^1\times H^1},
  \label{eq:riesz_a}
\end{equation}
for arbitrary $\delta v \in H^1 \times H^1$. Separating the $s$-component we have
\begin{gather}
  \int_0^T \bnabla_s^{L^2} \mJ \, \delta s \, dt +  \left(\beta_0\lnorm{u(s(T),T) - w(s(T))}^2 + 2\beta_2(s(T)-s_*)-\gamma \psi |_{(s(T),T)}\right){\delta s}(T)=\nonumber\\
  \int_0^T \left[ \bnabla_s^{H^1} \mJ \, \delta s +
  \ell^2_s \, \frac{d \bnabla_s^{H^1}\mJ}{dt} \, \frac{d\delta s}{dt} \right] \, dt
  \label{eq:riesz_b}
\end{gather}
for arbitrary $\delta s \in H^1$. Therefore, $ \bnabla_s^{H^1} \mJ \in H^1$ is a weak solution of the following boundary--value problem with measure right-hand side:
\begin{equation}
  \begin{aligned}
    \bnabla_s^{H^1} \mJ - \ell^2_s \, \frac{d^2}{dt^2} \, \bnabla_s^{H^1} \mJ & = \bnabla_s^{L_2} \mJ + \bnabla^{H^1}_{s(T)} \mJ(v)
    \qquad && \textrm{on} \ (0, \, T), \\
    \frac{d}{dt} \bnabla_s^{H^1} \mJ & = 0 && \textrm{for} \ t = 0, \, T.
  \end{aligned}
  \label{eq:helm}
\end{equation}
Similarly, from \eqref{eq:riesz_a} it follows that the $g$-componenet of the Frechet gradient $ \bnabla_g^{H^1} \mJ \in H^1$ is a weak solution of the boundary--value problem
\begin{equation}
  \begin{aligned}
    \bnabla_g^{H^1} \mJ - \ell^2_g \, \frac{d^2}{dt^2} \, \bnabla_g^{H^1} \mJ & = \bnabla_g^{L_2} \mJ    \qquad && \textrm{on} \ (0, \, T), \\
    \frac{d}{dt} \bnabla_s^{H^1} \mJ & = 0 && \textrm{for} \ t = 0, \, T.
  \end{aligned}
  \label{eq:helm1}
\end{equation}
We recall that by changing the value of parameters $\ell_s$ and $\ell_g$ we can control the smoothness of the gradient
$\bnabla_v^{H^1} \mJ$, and therefore also the relative smoothness of the resulting reconstruction of the
control vector $v$, and hence also the regularity of $v$.  More specifically, as was shown in
\cite{ProtasBewleyHagen04}, extracting cost functional gradients in the Sobolev spaces $H^p$, $p>0$,
is equivalent to applying a low--pass filter to $L_2$ gradients with the quantity $\ell$ representing
the ``cut-off'' scale.

It should be mentioned that the described procedure, also known as {\it preconditioning}, may be
considered as an alternative form to perform the projection $H \to V^1_R$ scheduled as Step 6 in the
iterative algorithm of Section~\ref{sec:grad_besov}. Based on this algorithm, we could finally conclude
that the iterative reconstruction of control vector $v = (s,g)$ involves computations summarized below
in Algorithm~\ref{alg:gen_opt}.
\begin{algorithm}[!htb]
  \begin{algorithmic}
    \STATE $k \leftarrow 0$
    \STATE $v_0 \leftarrow $ initial guess $(s_{\rm ini}, g_{\rm ini})$
    \REPEAT
    \STATE given estimate of $v_k$, solve forward (Neumann) problem \eqref{eq:pde-1}--\eqref{eq:pde-stefan}
      to find $u_k=u(x,t;v_k)$
    \STATE evaluate $\J(v_k)$ using \eqref{Eq:W:1:8}
    \STATE given $u_k$ and $v_k$, solve adjoined PDE problem \eqref{eq:adj-pde}--\eqref{eq:adj-robin-free}
      to find $\psi_k=\psi(x,t;v_k)$
    \STATE  compute cost functional gradient components $\bnabla^{L^2}_s \mJ(v_k), \bnabla^{H^1}_{s(T)} \mJ(v_k), \bnabla^{L^2}_g \mJ(v_k)$ according to \eqref{eq:grads} with $v, u$ and $\psi$ replaced with $v_k, u_k$ and $\psi_k$ respectively.
         \STATE given $\bnabla^{L_2}_s \mJ(v_k)$, $\bnabla^{H^1}_{s(T)} \mJ(v_k)$ and $\bnabla^{L_2}_g \mJ(v_k)$, solve \eqref{eq:helm} and/or \eqref{eq:helm1}
      compute preconditioned Sobolev gradients $\bnabla^{H^1}_s \mJ(v_k)$ and/or $\bnabla^{H^1}_g \mJ(v_k)$
      for chosen value of parameter $\ell_v$
    \STATE find optimal values for stepsize parameters $\alpha^s_k$ and $\alpha^g_k$ by solving two-dimensional
      optimization problem
      \begin{equation}
        \alpha_k = (\alpha^s_k, \alpha^g_k) =
        \underset{\alpha > 0}{\argmin}\,\,\mJ \left( v_k - \alpha \bnabla^{H^1}_v \mJ(v_k) \right)
        \label{eq:opt_alpha}
      \end{equation}
    \STATE update $v_{k+1} = (s_{k+1}, g_{k+1})$ with a descent direction derived from gradients
      $\bnabla^{H^1}_s \mJ(v_k)$ and $\bnabla^{H^1}_g \mJ(v_k)$
      \begin{subequations}
        \label{eq:opt_SD}
        \begin{alignat}{2}
          s_{k+1} &= s_k - \alpha^s_k \bnabla^{H^1}_s \mJ(v_k), \label{eq:opt_SD_s}\\
          g_{k+1} &= g_k - \alpha^g_k \bnabla^{H^1}_g \mJ(v_k) \label{eq:opt_SD_g}
        \end{alignat}
      \end{subequations}
    \STATE $k \leftarrow k + 1$
    \UNTIL the termination criteria \eqref{convergencecriteria} are satisfied to a given tolerance $\epsilon$
  \end{algorithmic}
  \caption{Adjoint gradient workflow for solving the inverse Stefan problem}
  \label{alg:gen_opt}
\end{algorithm}

\begin{figure}[htb!]
  \begin{center}
  \mbox{
  \subfigure[]{\includegraphics[width=0.5\textwidth]{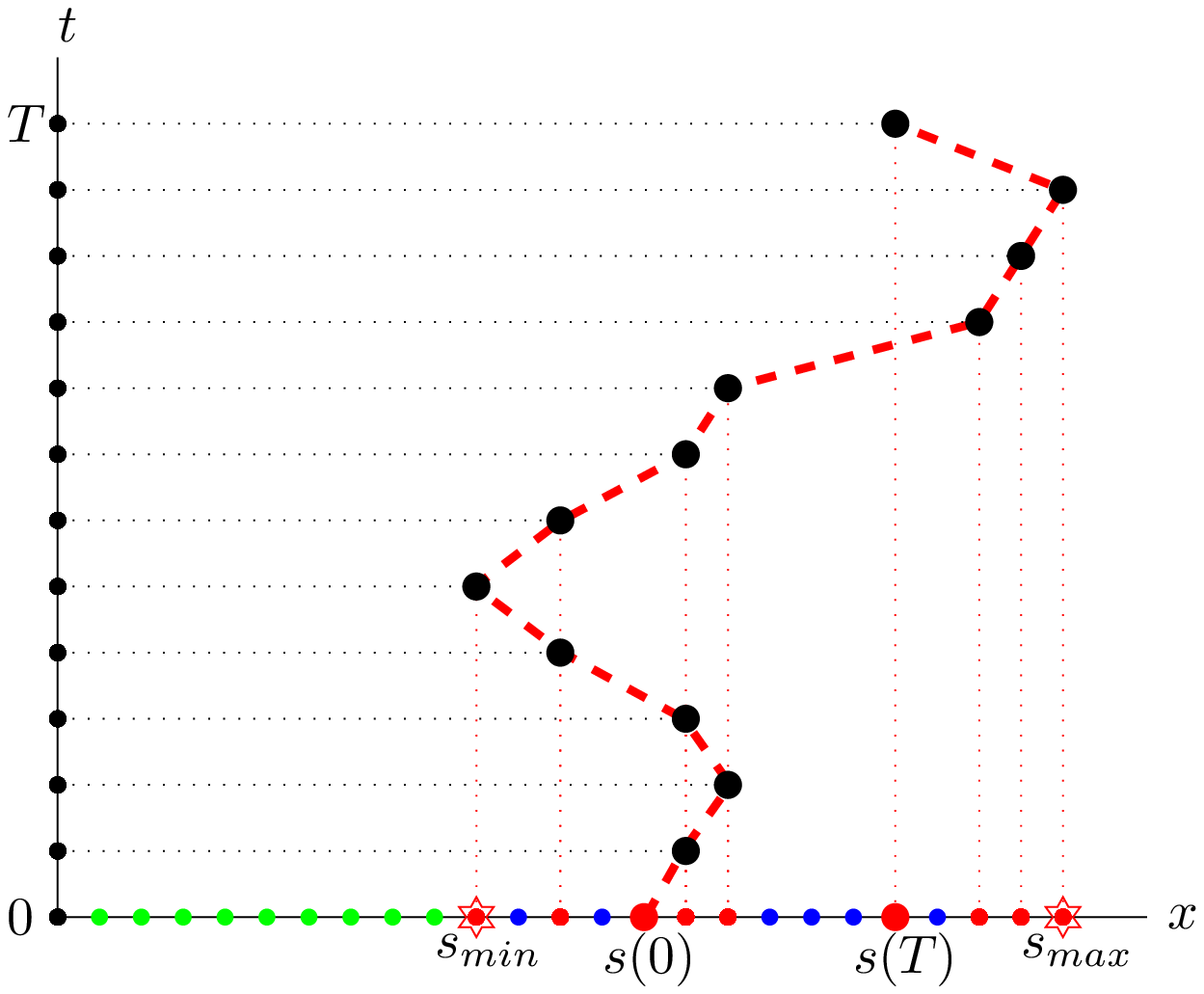}}
  \subfigure[]{\includegraphics[width=0.5\textwidth]{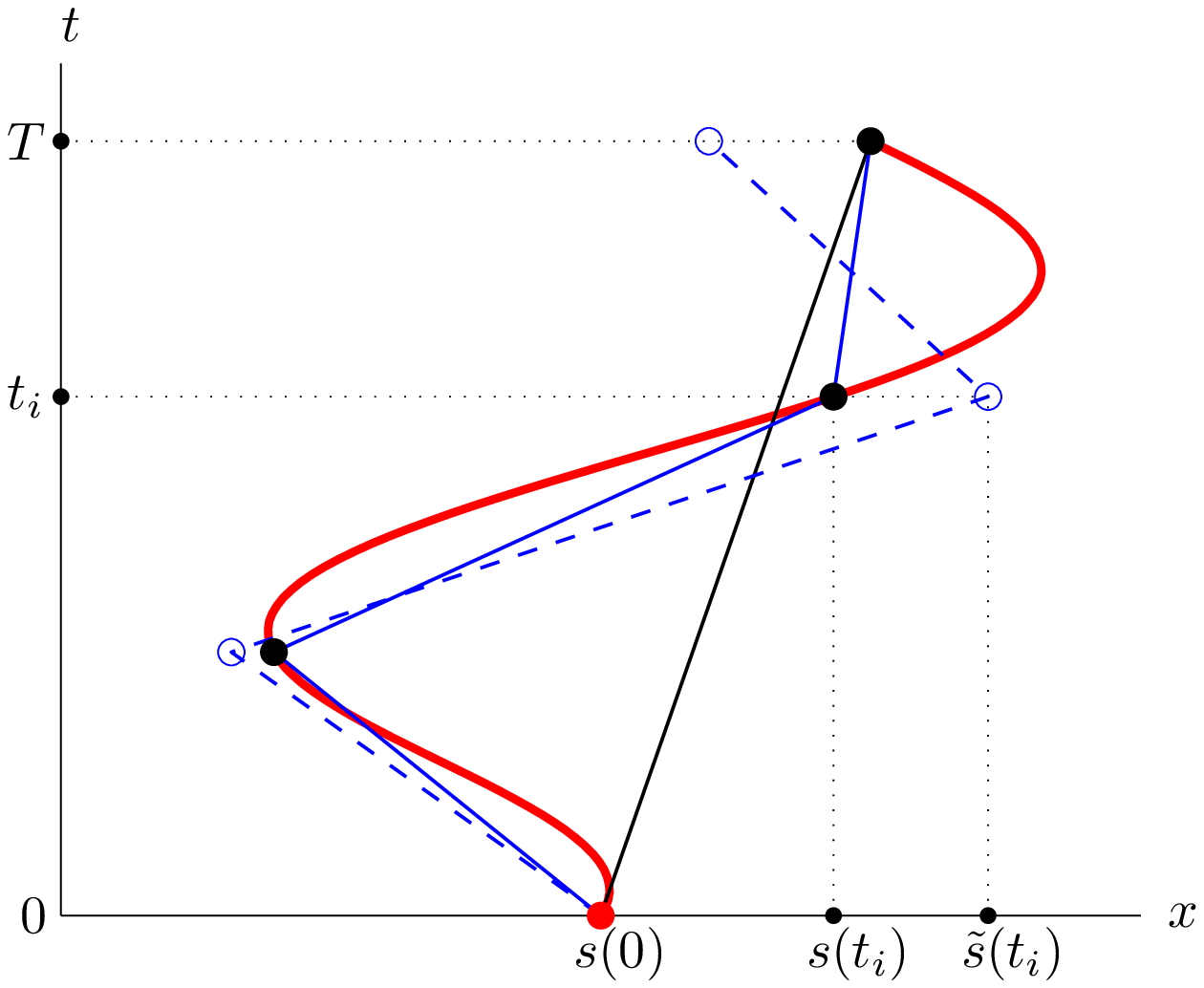}}}
  \end{center}
  \caption{(a)~Schematic showing uniform discretization of $t$-domain (black dots along $t$-axis)
    and non-uniform discretization of $x$-domain (red, green and blue dots along $x$-axis). Black filled
    circles connected by red dashed lines represent free boundary $s(t_j)$ for $j = 1, \ldots,
    n$. Red hexagons are used as endpoints of interval $[s_{\rm{min}}, s_{\rm{max}}]$ where
    discretization of $x$-domain is non-uniform. (b)~Schematic showing of (red line) true
    free boundary $s_{\rm true}(t)$ and additional measurements for $s(t)$ (black filled circles) without
    noise $\{s(t_i) \}_{i=1}^M$ and (blue empty circles) with noise $\{\tilde s^{\eta}_i \}_{i=1}^M$.
    Piecewise linear approximations of $s(t)$ are shown without noise for (black solid line) $M = 1$
    and (blue solid line) $M = 3$, and (blue dashed line) with noise for $M = 3$.}
  \label{fig:discr_tx}
\end{figure}

To validate the accuracy and performance of the proposed approach to find a triple $\{u(x,t),s(t),g(t)\}$ we solve the inverse
Stefan problem \eqref{eq:pde-1}--\eqref{eq:pde-finaltemp} on $\Omega = \{(x,t): 0 < x < s(t),~0 < t \leq T\}$
discretized for $t$ and $x$ according to procedure described in Section~\ref{sec:discr_conv} and demonstrated in Figure~\ref{fig:discr_tx}(a). The $t$-domain is discretized uniformly
(black dots along $t$-axis) using $n$ partial time intervals each of size $\tau = T/n$. 
Spatial discretization for every time instance $t_j = j \tau$, $j = 0, \ldots, n$ is non-uniform. First, the interval from $x = 0$ to $s_{\rm{min}} = \underset{j}{\min}\,s(t_j)$ is discretized
uniformly with step $h \leq h_x$ (green dots along $x$-axis). The rest of the $x$-domain from $x = s_{\rm{min}}$
to $s_{\rm{max}} = \underset{j}{\max}\,s(t_j)$ is discretized by using all available values $s(t_j)$ for
$j = 0, \ldots, n$ sorted in ascending order (red dots). The step size constraint $0 < \epsilon_x \leq h_i \leq h_x$
is enforced to maintain overall accuracy and to prevent grid points to create unreasonably small intervals. Due to the latter
some points are removed, while the former requires to add additional points (blue dots) to the grid. Such discretization
technique is used to solve both forward \eqref{eq:pde-1}--\eqref{eq:pde-stefan} and adjoined PDE \eqref{eq:adj-pde}--\eqref{eq:adj-robin-free} problems to reduce the error accumulated due to interpolating the solutions
near free boundary $x = s(t_j)$. Both parameters $h_x > 0$ and $\epsilon_x > 0$ have to be sufficiently small, but
allow reasonable time to compute the solution for both PDE problems.

In parallel with discretization of both PDE problems, we also discretize continuous measurements
$\mu(t)$ and $w(x)$, given respectively by \eqref{eq:pde-freebound} and \eqref{eq:pde-finaltemp}, using
discretized (pointwise) measurement data $u(x,t)$ which are typically available in actual experiments.
To mimic an actual experimental procedure, model true functions $s_{\rm true}(t)$ and $g_{\rm true}(t)$
are used in combination with PDE system \eqref{eq:pde-1}--\eqref{eq:pde-stefan} to obtain discretized
measurements $\mu(t_j) = u(s(t_j),t_j)$ and $w(x) = u(x,T)$. Functions $s_{\rm true}(t)$ and $g_{\rm true}(t)$
are then ``forgotten'' and reconstructed using the proposed gradient--based framework. While in the calculations
to validate our basic formulation, presented in Sections~\ref{sec:grad_valid}, \ref{sec:single_rec} and
\ref{sec:multiple_rec}, no noise is present in the measurements, its effect is addressed in Section~\ref{sec:noise}.

In terms of the initial guess for every computational model described in the next section we take a
constant approximation $g_{\rm ini} = \dfrac{1}{n+1}\sum_{j=0}^{n} g(t_j)$ to $g(t)$. Unless stated otherwise,
a line segment to connect points $(s(0),0)$ and $(s(T),T)$ is chosen as reasonable initial approximation $s_{\rm ini}(t)$
to $s(t)$ as shown by black solid line in Figure~\ref{fig:discr_tx}(b). We also refer to
Section~\ref{sec:models} for more details.

Our code for solving forward problem \eqref{eq:pde-1}--\eqref{eq:pde-stefan} and adjoined problem
\eqref{eq:adj-pde}--\eqref{eq:adj-robin-free} has been implemented using {\tt FreeFem++} \cite{FreeFem12},
an open--source, high--level integrated development environment for the numerical solution of PDEs based on 
the Finite Element Method (FEM). To solve numerically both problems spatial discretization of domain
\eqref{eq:pde-domain-defn} is carried out using 1D finite elements and the P1 piecewise linear (continuous)
representations for all spatially distributed quantities. The system of algebraic equations obtained after
such discretization at every time step is solved with {\tt UMFPACK}, a solver for nonsymmetric sparse linear
systems \cite{umfpack}. The same meshes for both $t$- and $x$-domains and discussed above P1 representations
are then used to construct the gradients \eqref{eq:grads} and, if required, to project these gradients
onto the Hilbert-Sobolev-Besov space $H^1$ by solving \eqref{eq:helm} and \eqref{eq:helm1} to perform the iterative optimization procedure as
described in Algorithm~\ref{alg:gen_opt}. As seen at \eqref{eq:opt_SD}, this procedure is utilizing the
Steepest Descent (SD) approach \cite{Nocedal06} with stepsize parameters $\alpha_k^s$ and $\alpha_k^g$
obtained by applying line minimization search \cite{NumericalRecipes07} to solve optimization problem
\eqref{eq:opt_alpha}. Unless otherwise stated, for reconstructing the entire control vector $(s,g)$
the same value of the stepsize, i.e.~$\alpha_k^s = \alpha_k^g$, is used for both functions $s(t)$ and $g(t)$.
All cost functional weighting coefficients in \eqref{Eq:W:1:8} are set to be equal,
i.e.~$\beta_0 = \beta_1 = \beta_2 = 1$. The termination condition used is
$\left| \frac{\J(v_{k})-\J(v_{k-1})}{\J(v_{k-1})} \right| < 10^{-5}$.

\subsection{Description of models}
\label{sec:models}

\begin{figure}[htb!]
  \begin{center}
  \mbox{
  \subfigure[model \#1: $u(x,t)$, $s(t)$]{\includegraphics[width=0.33\textwidth]{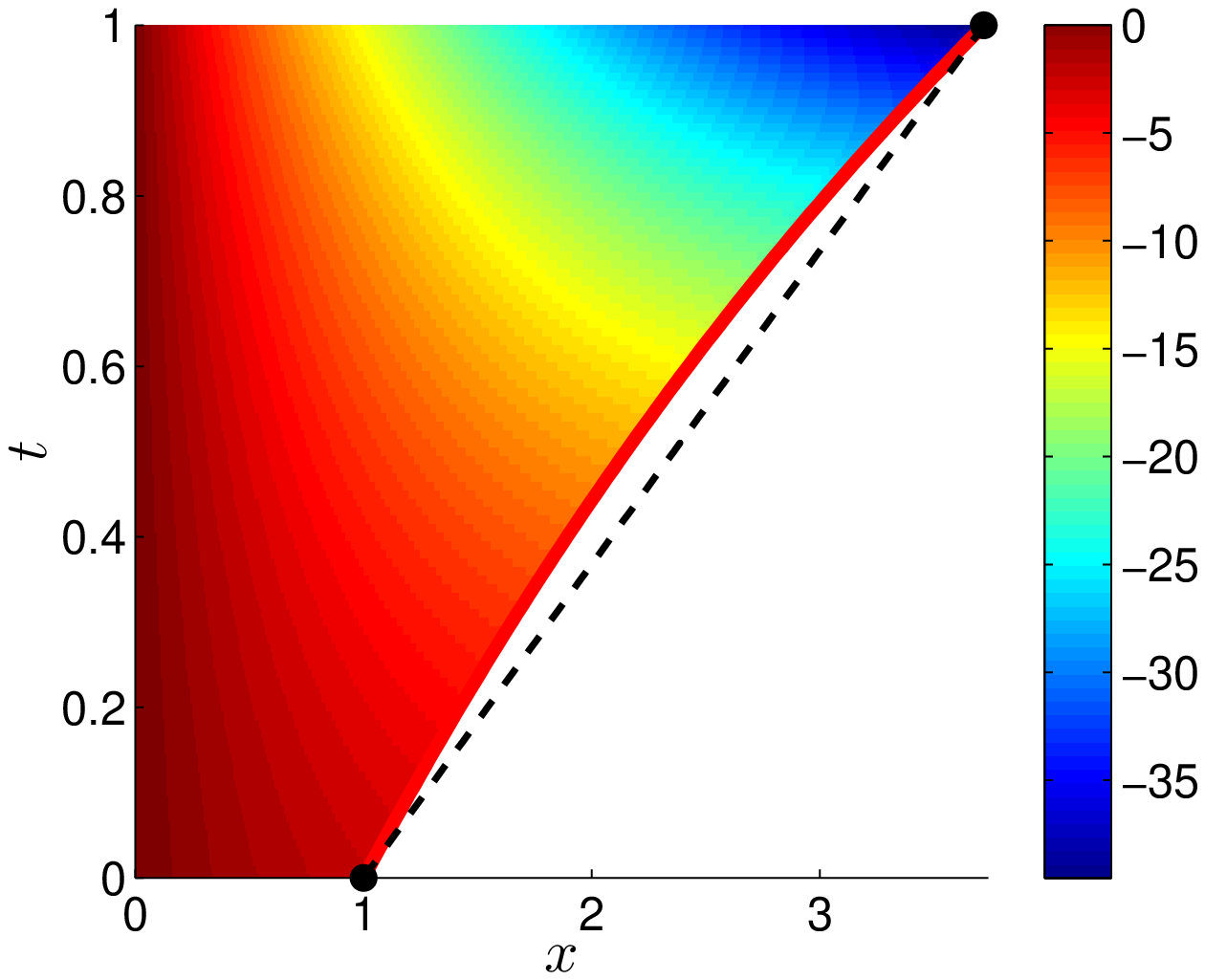}}
  \subfigure[model \#2: $u(x,t)$, $s(t)$]{\includegraphics[width=0.33\textwidth]{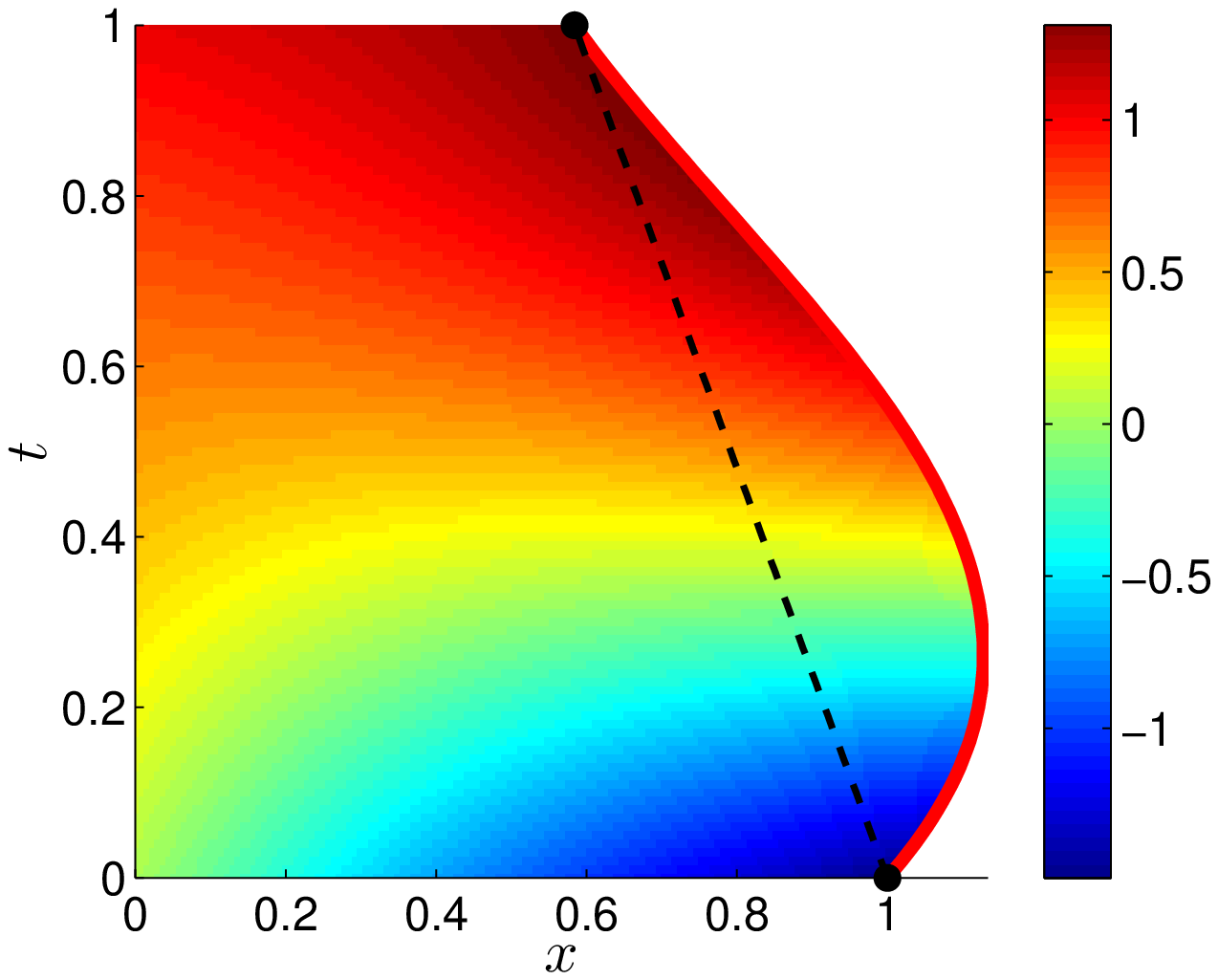}}
  \subfigure[model \#3: $u(x,t)$, $s(t)$]{\includegraphics[width=0.33\textwidth]{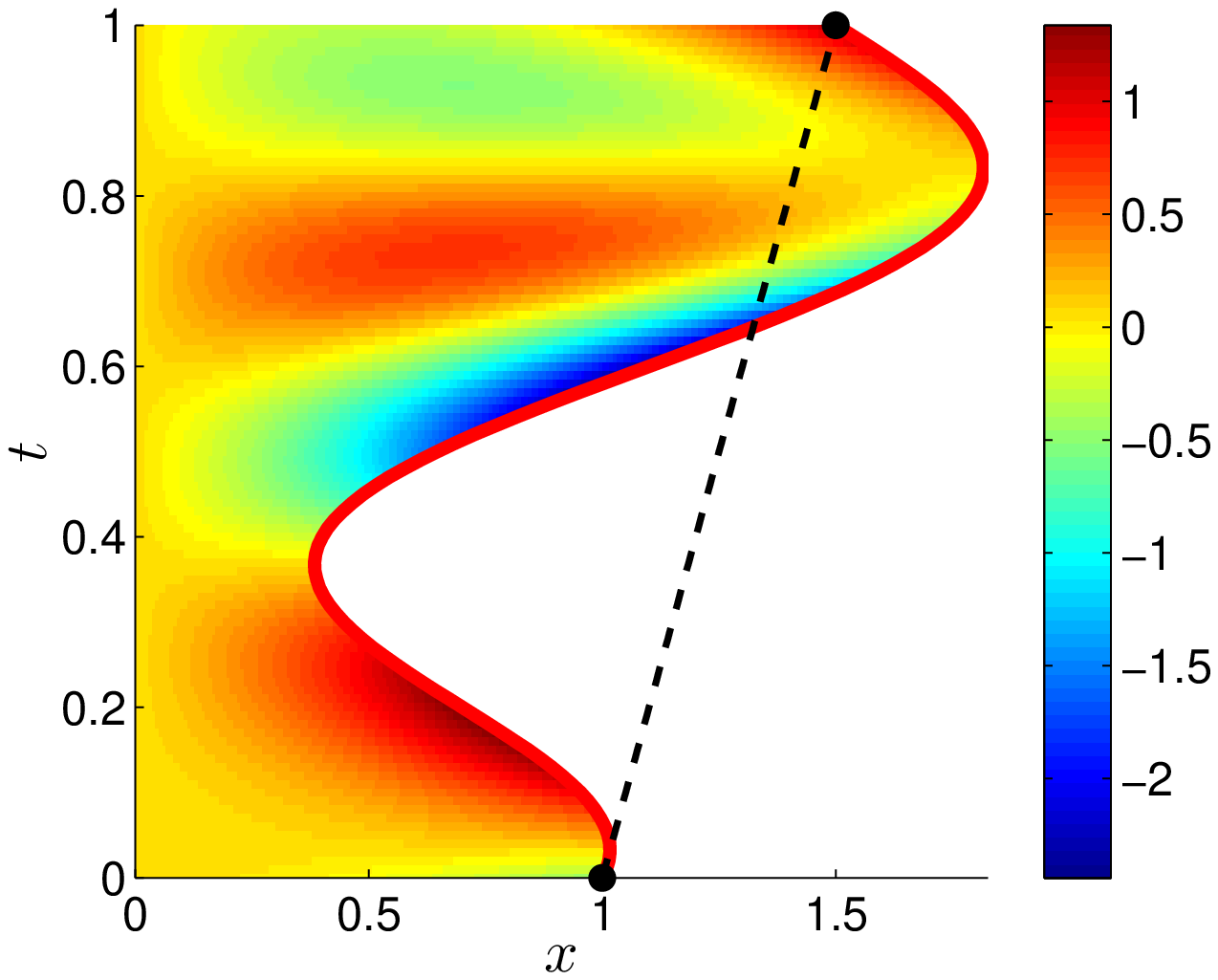}}}
  \mbox{
  \subfigure[model \#1: $g(t)$]{\includegraphics[width=0.33\textwidth]{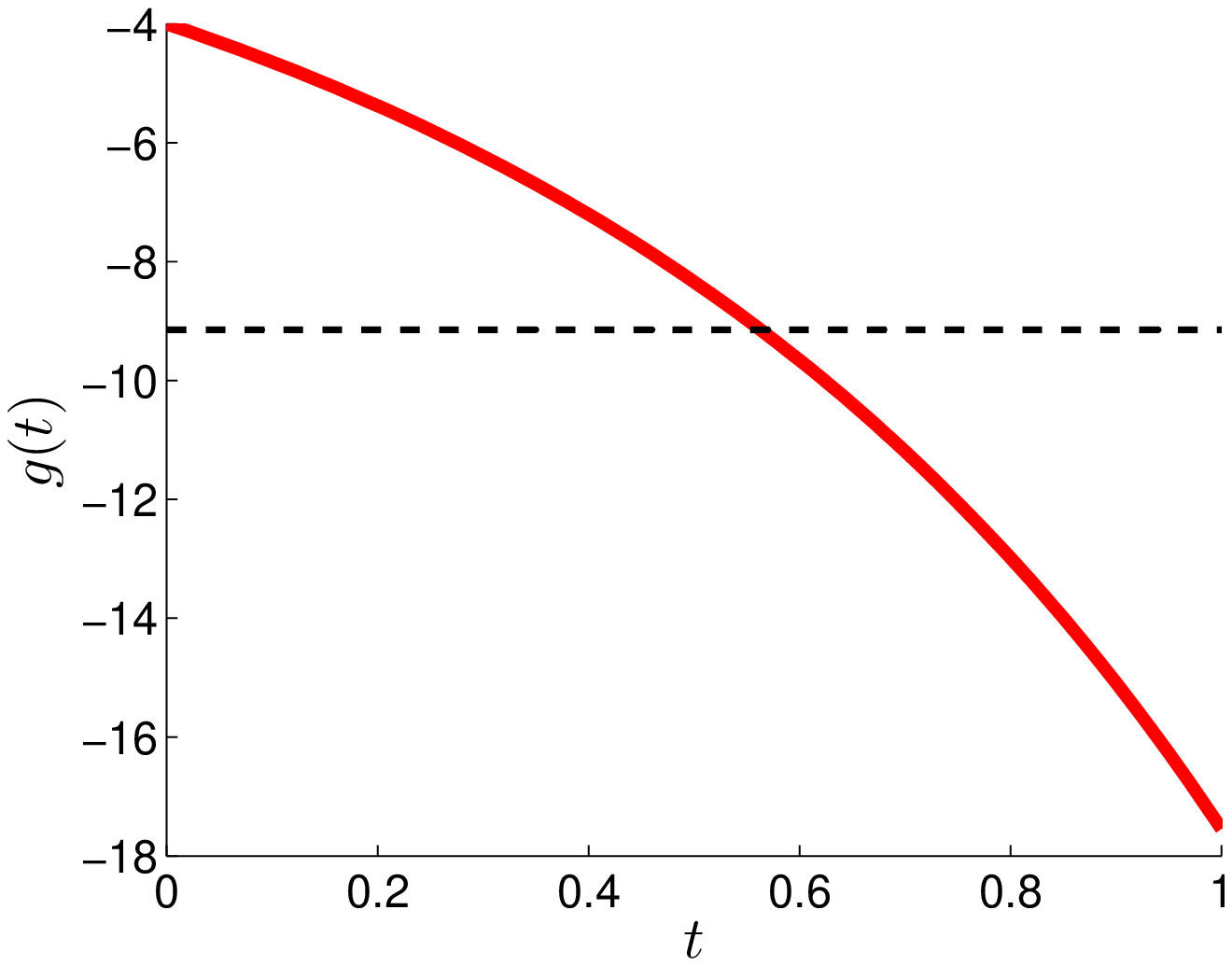}}
  \subfigure[model \#2: $g(t)$]{\includegraphics[width=0.33\textwidth]{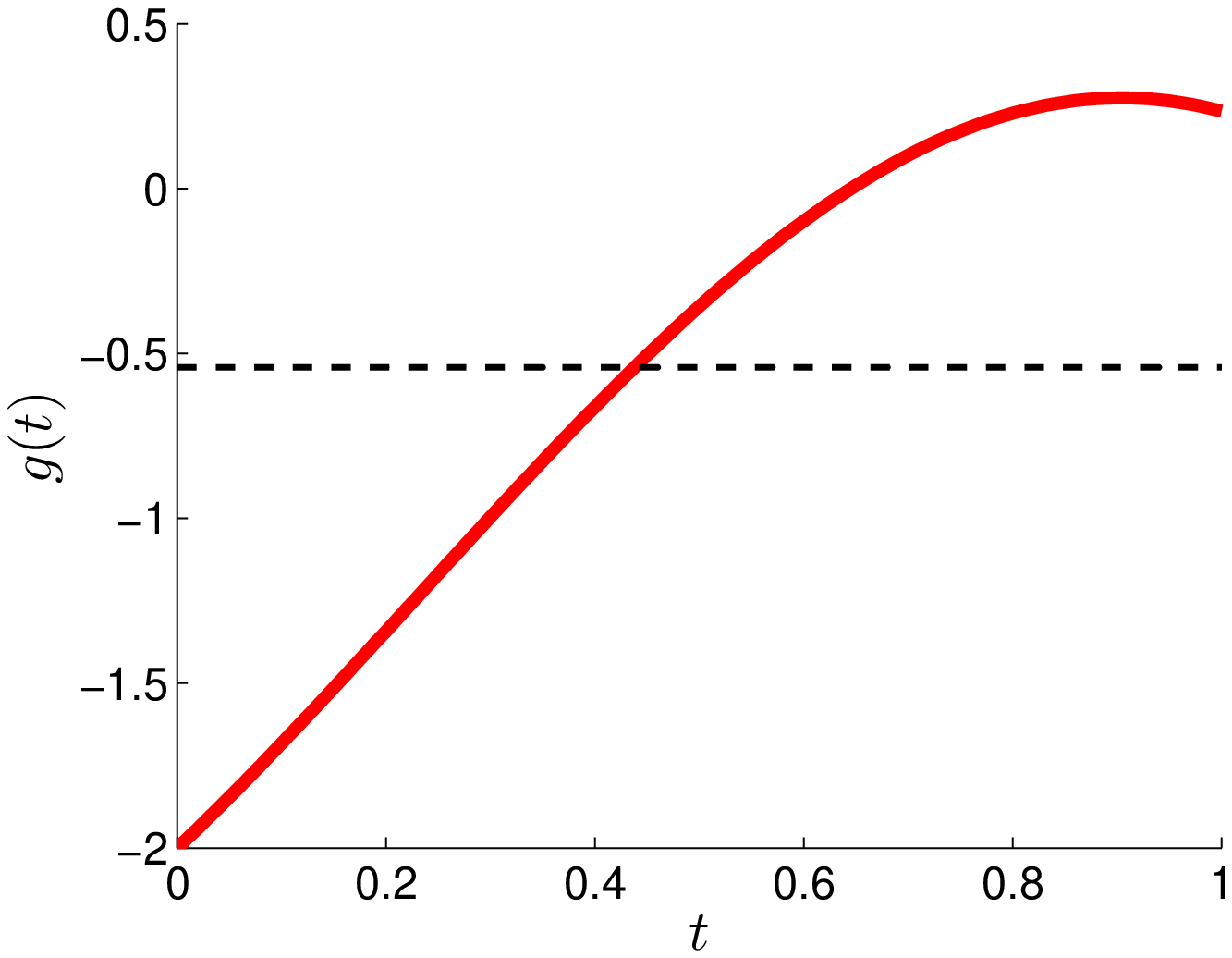}}
  \subfigure[model \#3: $g(t)$]{\includegraphics[width=0.33\textwidth]{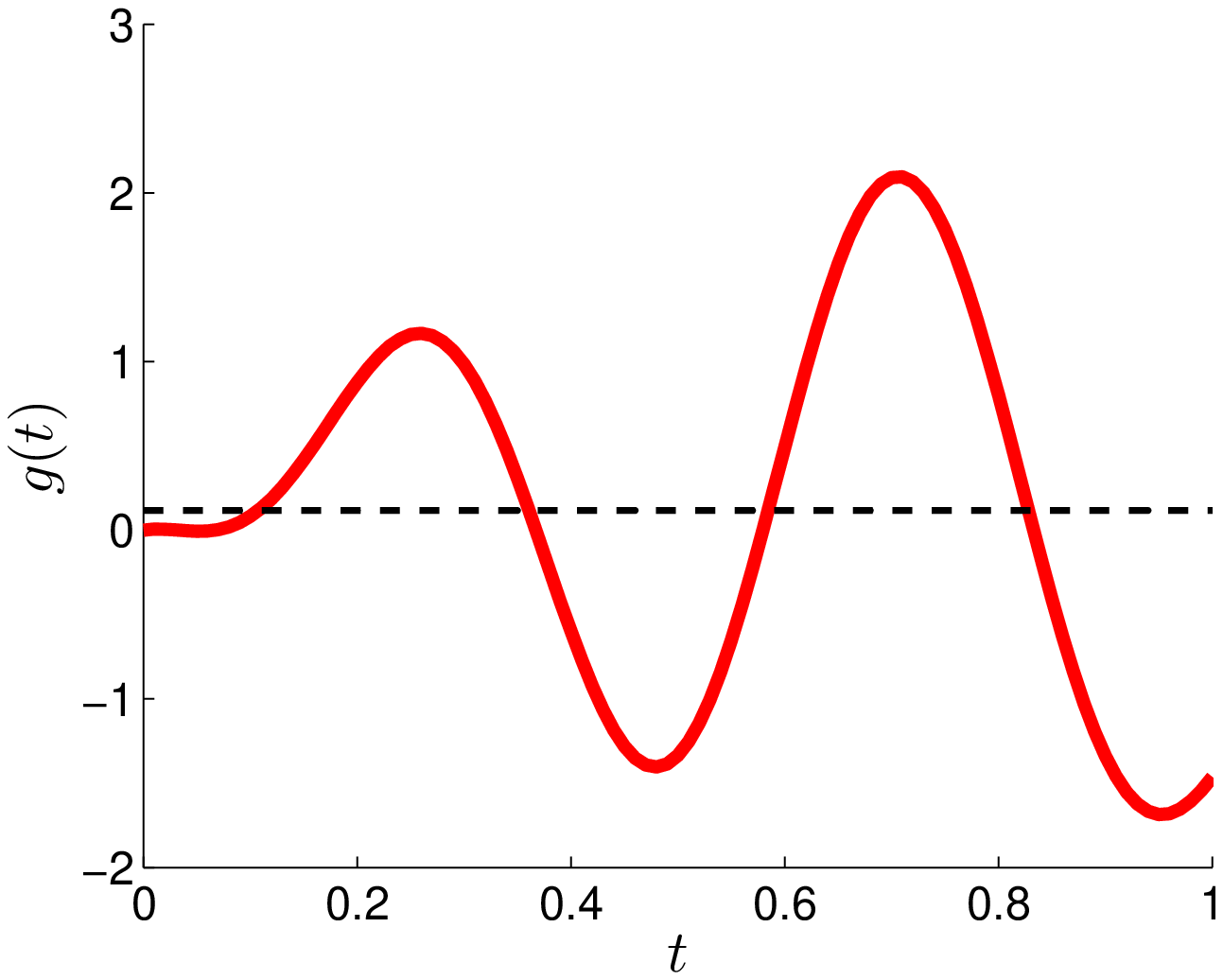}}}
  \end{center}
  \caption{Functions (a-c) $u(x,t)$, $s(t)$ and (d-f) $g(t)$ shown for (a,d) model \#1,
    (b,e) model \#2 and (c,f) model \#3. In (a-c) vertical color bars represent the values for $u(x,t)$;
    solid red lines show the true shape of free boundary $s_{\rm true}(t)$. In (d-f) solid red lines
    represent $g_{\rm true}(t)$ functions. Dashed black lines show initial guesses (a-c)~$s_{\rm ini}(t)$ and
    (d-f)~$g_{\rm ini}(t)$.}
  \label{fig:models}
\end{figure}

To perform computations as described in Section~\ref{sec:approaches} we define three model examples each of
different complexity as described in detail below. For simplicity, but without loss of generality,
the following functions in all models are set to constant values:
\begin{equation}
  a(x,t) = 1, \qquad b(x,t) = 0, \qquad \gamma(s(t),t) = 1, \qquad \chi(s(t),t) = 0.
\end{equation}
%
%
%
Below are shown the analytical expressions for three models which are numbered (from \#1 to \#3) in
the order of increasing complexity:
\begin{itemize}
   \item {\bf Model \#1:}
   \begin{equation}
     \begin{aligned}
       u(x,t) &= - (1 + e^t) \left[ x (t + e^t + 1) - \frac{1}{2} x^2 \right], \qquad c(x,t) = x + t,\\
       s(t) &= t + e^t, \qquad g(t) = - (1 + e^t) (t + e^t + 1).
     \end{aligned}
     \label{eq:model_1}
   \end{equation}

   \item {\bf Model \#2:}
   \begin{equation}
     \begin{aligned}
       u(x,t) &= \frac{1}{2} x^2 + x \left( -\cos 2t - t + 2\sin 2t - 1 \right) + t, \qquad c(x,t) = x + t,\\
       s(t) &= \cos 2t + t, \qquad g(t) = -\cos 2t - t + 2\sin 2t - 1.
     \end{aligned}
     \label{eq:model_2}
   \end{equation}

   \item {\bf Model \#3:}
   \begin{equation}
     \begin{aligned}
       u(x,t) &= \frac{5\pi}{8} x^2 \sin \frac{5\pi}{2} t - \frac{5\pi}{16} x \sin 5\pi t
       - \frac{5\pi}{4} \left( t - \frac{1}{2} \right) x \sin \frac{5\pi}{2} t - \frac{x^2}{2}
       + \frac{x}{2} \cos \frac{5\pi}{2} t + tx - \frac{1}{2} x, \\
       s(t) &= \frac{1}{2} \cos \frac{5\pi}{2} t + t + \frac{1}{2}, \qquad
       c(x,t) = x + t, \\
       g(t) &= -\frac{5\pi}{16} \sin 5\pi t - \frac{5\pi}{4} \left( t - \frac{1}{2} \right) \sin \frac{5\pi}{2} t
       +\frac{1}{2} \cos \frac{5\pi}{2} t + t - \frac{1}{2}.
     \end{aligned}
     \label{eq:model_3}
   \end{equation}
\end{itemize}
Figure~\ref{fig:models} shows functions $u(x,t)$ (in color), $s(t)$ and $g(t)$ (red lines) for all models.
As described in Section~\ref{sec:approaches}, initial guesses to reconstruct $s(t)$ for all models are chosen
as line segments to connect $(s(0),0)$ and $(s(T),T)$ shown by black dots. Figure~\ref{fig:models} also shows
initial guesses $s_{\rm ini}(t)$ and $g_{\rm ini}(t)$ respectively for both $s(t)$ and $g(t)$ as dashed black lines.

In our computations for all three models, we used the following parameters for time and space discretization,
described previously in Section~\ref{sec:approaches}: $T = 1$, $n= 100$, $h_x = 0.01$, $\epsilon_x = 10^{-8}$.
This choice is motivated by finding optimal balance between reasonable computational time and appropriate quality
of cost functional gradients $\bnabla_s^{L_2} \J(v)$ and $\bnabla_g^{L_2} \J(v)$. The discussion on the latter
could be found in the next section.

\subsection{Validation of gradients}
\label{sec:grad_valid}

\begin{figure}[htb!]
  \begin{center}
  \mbox{
  \subfigure[]{\includegraphics[width=0.5\textwidth]{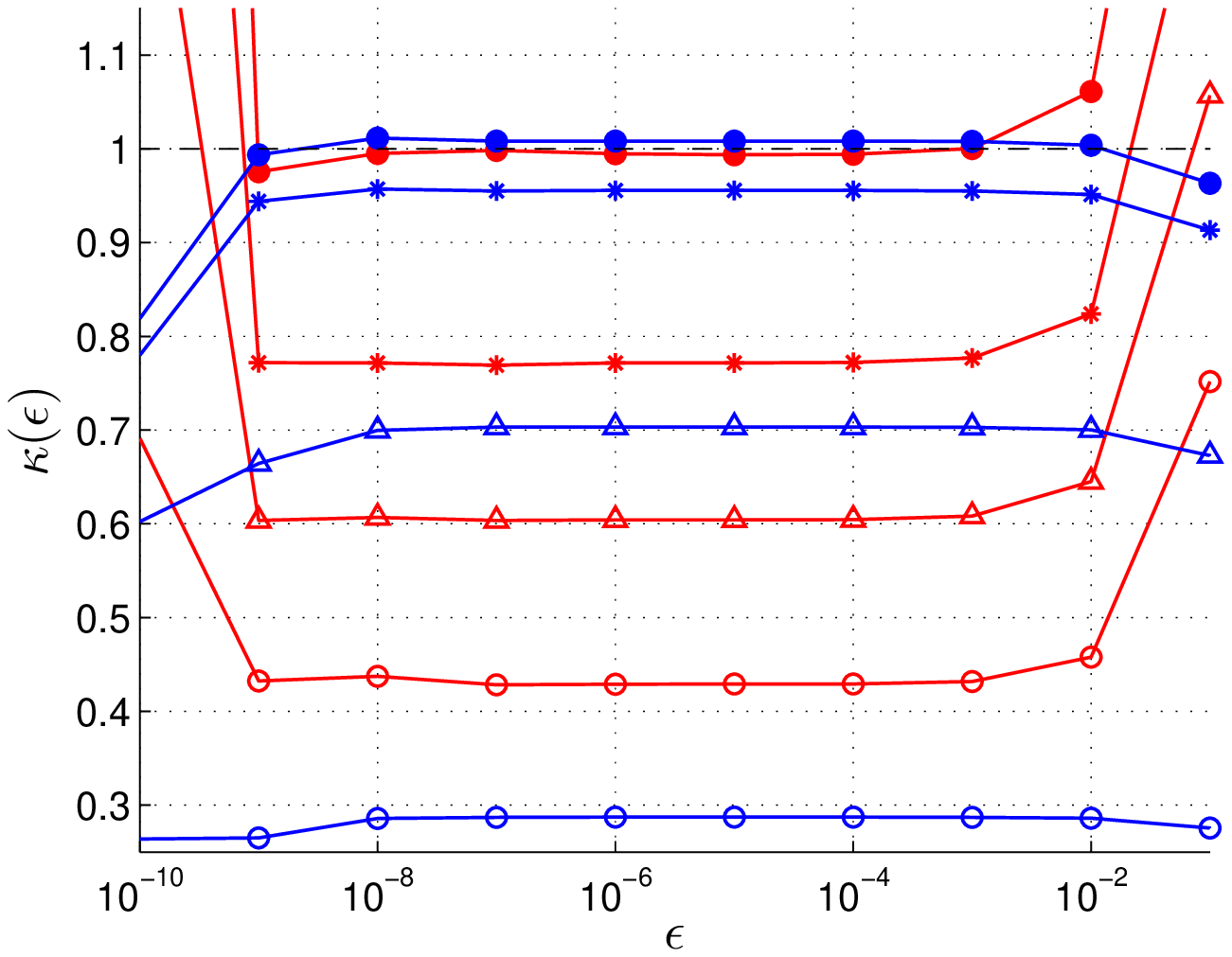}}
  \subfigure[]{\includegraphics[width=0.5\textwidth]{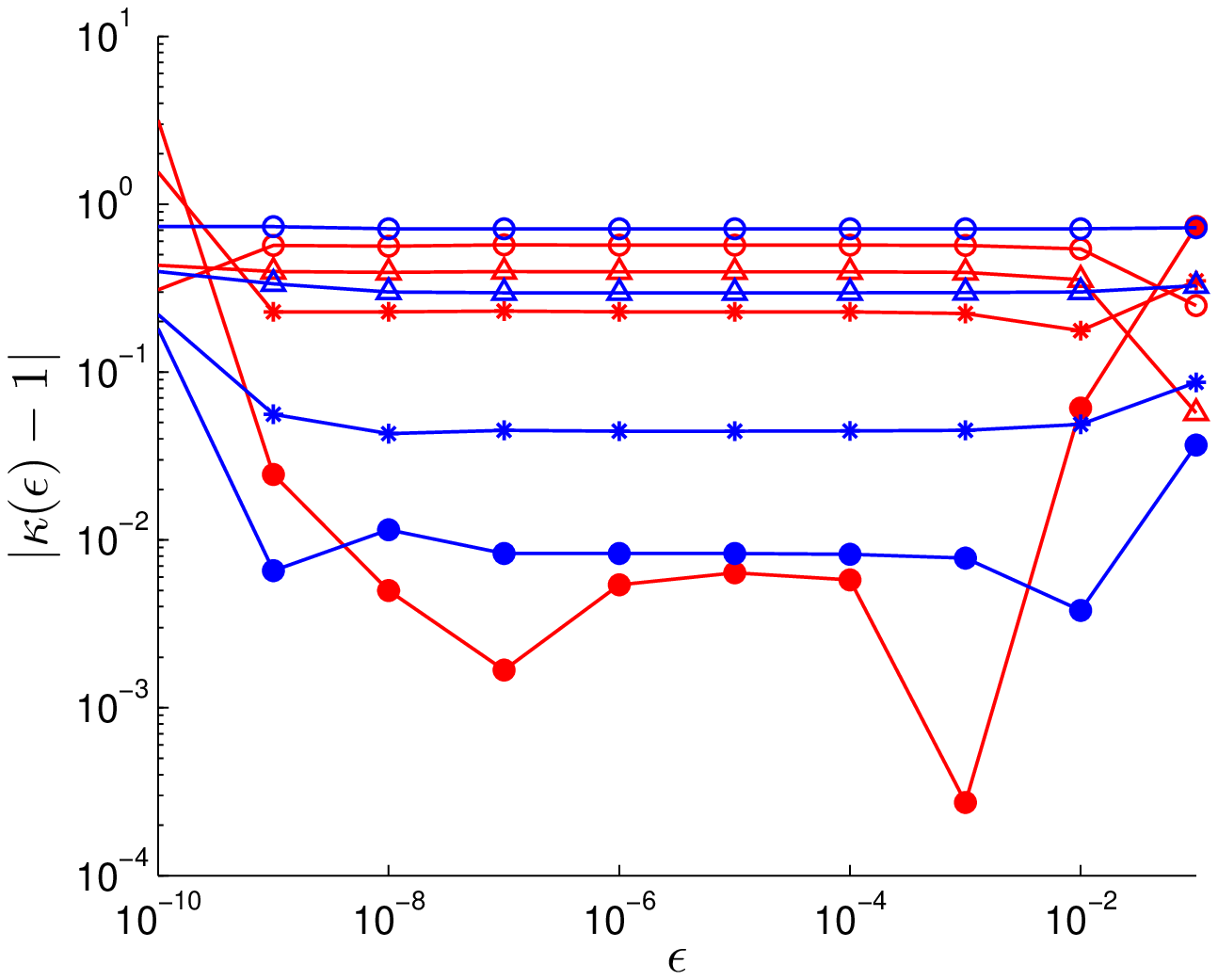}}}
  \end{center}
  \caption{The behavior of (a) $\kappa (\epsilon)$ and (b) $\log_{10} |\kappa(\epsilon) -1 |$ as a function of $\epsilon$
    for both controls (red) $s(t)$ and (blue) $g(t)$ with fixed space discretization parameter $h_x = 0.01$.
    Time steps used in discretizing the $t$-domain for computing $\bnabla_s^{L_2} \J(v)$ are (empty circles)
    $\Delta t = 5 \cdot 10^{-3}$, (triangles) $\Delta t = 2 \cdot 10^{-3}$, (asterisks) $\Delta t = 1 \cdot 10^{-3}$,
    (filled circles) $\Delta t = 1 \cdot 10^{-4}$. Time steps for computing $\bnabla_g^{L_2} \J(v)$ are (empty circles)
    $\Delta t = 2 \cdot 10^{-2}$, (triangles) $\Delta t = 1 \cdot 10^{-2}$, (asterisks) $\Delta t = 5 \cdot 10^{-3}$,
    (filled circles) $\Delta t = 4 \cdot 10^{-3}$. All tests are performed for model \#3.}
  \label{fig:kappa_test}
\end{figure}

In this section we present results demonstrating the consistency of the cost functional gradients obtained with
the approach described in Section~\ref{sec:opt_ctrl_problem} and Algorithm~\ref{alg:gen_opt}.
Figure~\ref{fig:kappa_test} shows the results of a diagnostic test commonly employed to verify the correctness
of the cost functional gradients (see, e.g., \cite{Bukshtynov11,Bukshtynov13}) computed for model \#3. It consists
in computing the Fr\'echet differential $d\mJ(v; \delta v) = \left\langle{}\J'(v),{\delta v}\right\rangle_H$ for
some selected variations (perturbations) $\delta v$ in two different ways, namely, using a finite--difference
approximation and using \eqref{eq:grads} which is based on the adjoint field, and then examining the ratio
of the two quantities, i.e.,
\begin{equation}
  \kappa (\epsilon) = \dfrac{\frac{1}{\epsilon} \left[ \J(v + \epsilon \, \delta v) - \J(v) \right]}
  { \left\langle{}\J'(v),{\delta v}\right\rangle_H}
  \label{eq:kappa}
\end{equation}
for a range of values of $\epsilon$. As the sensitivity of the cost functional $\mJ(v)$ with respect to $v$
may vary significantly for the different contributions of $s(t)$ and $g(t)$, it is reasonable to perform this
test separately for different parts of the gradient, namely $\bnabla_s^{L^2} \J(v)$, $\bnabla_{s(T)}^{H^1} \J(v)$
and $\bnabla_g^{L^2} \J(v)$. If these gradients are computed correctly, then for intermediate values of
$\epsilon$, $\kappa(\epsilon)$ will be close to the unity. Remarkably, this behavior can be observed in
Figure~\ref{fig:kappa_test} over a range of $\epsilon$ spanning about 6 orders of magnitude for both controls
$s(t)$ and $g(t)$. Furthermore, we also emphasize that refining time step $\Delta t$ in discretizing
the $t$-domain while solving both forward \eqref{eq:pde-1}--\eqref{eq:pde-stefan} and adjoined 
\eqref{eq:adj-pde}--\eqref{eq:adj-robin-free} PDE problems yields values of $\kappa (\epsilon)$ closer to the unity.
The reason is that in the ``optimize--then--discretize'' paradigm adopted here such refinement of discretization
leads to a better approximation of the continuous gradient \cite{ProtasBewleyHagen04}. The quality of this
approximation may be further improved by refining parameter $h_x$ of the $x$-domain discretization. However, our
non-uniform $x$-discretization described previously in Section~\ref{sec:approaches} makes the systematic validation
rather complicated, thus it is not considered here. We add that the quantity $\log_{10} | \kappa (\epsilon) - 1 |$
plotted in Figure~\ref{fig:kappa_test}b shows how many significant digits of accuracy are captured in a given
gradient evaluation. As can be expected, the quantity $\kappa(\epsilon)$ deviates from the unity for very small
values of $\epsilon$, which is due to the subtractive cancellation (round--off) errors, and also for large values
of $\epsilon$, which is due to the truncation errors, both of which are well--known effects.

\subsection{Identification of the free boundary}
\label{sec:single_rec}

In this section we present results demonstrating the performance of the proposed numerical approach to identify
free boundary $s(t)$ only. At this point, Algorithm~\ref{alg:gen_opt} is used to find (local) optimal solution
$\bs (t)$ iteratively starting from initial guess $s = s_{\rm ini}(t)$ and setting $g_k(t)$ for every
$k = 0, 1, 2, \ldots$ to the true expressions defined analytically in \eqref{eq:model_1}--\eqref{eq:model_3}.
Thus, Algorithm~\ref{alg:gen_opt} is modified appropriately by skipping computing the corresponding part of
the gradient, namely $\bnabla_g \mJ$ in \eqref{eq:grads}, and setting $\alpha_k^g$ in \eqref{eq:opt_SD} to zero.

\begin{figure}[htb!]
  \begin{center}
  \mbox{
  \subfigure[$k = 1$]{\includegraphics[width=0.5\textwidth]{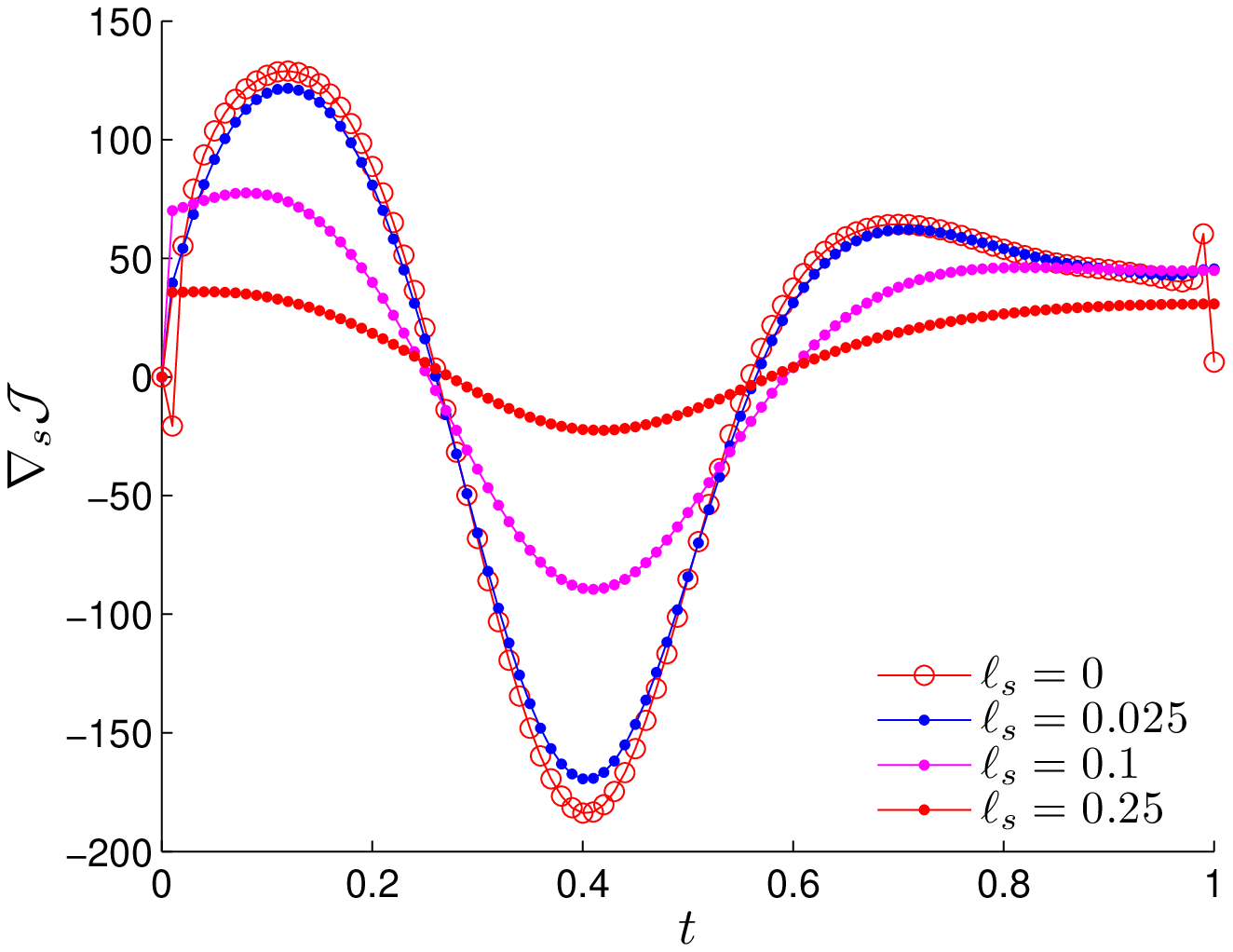}}
  \subfigure[$k = 38$]{\includegraphics[width=0.5\textwidth]{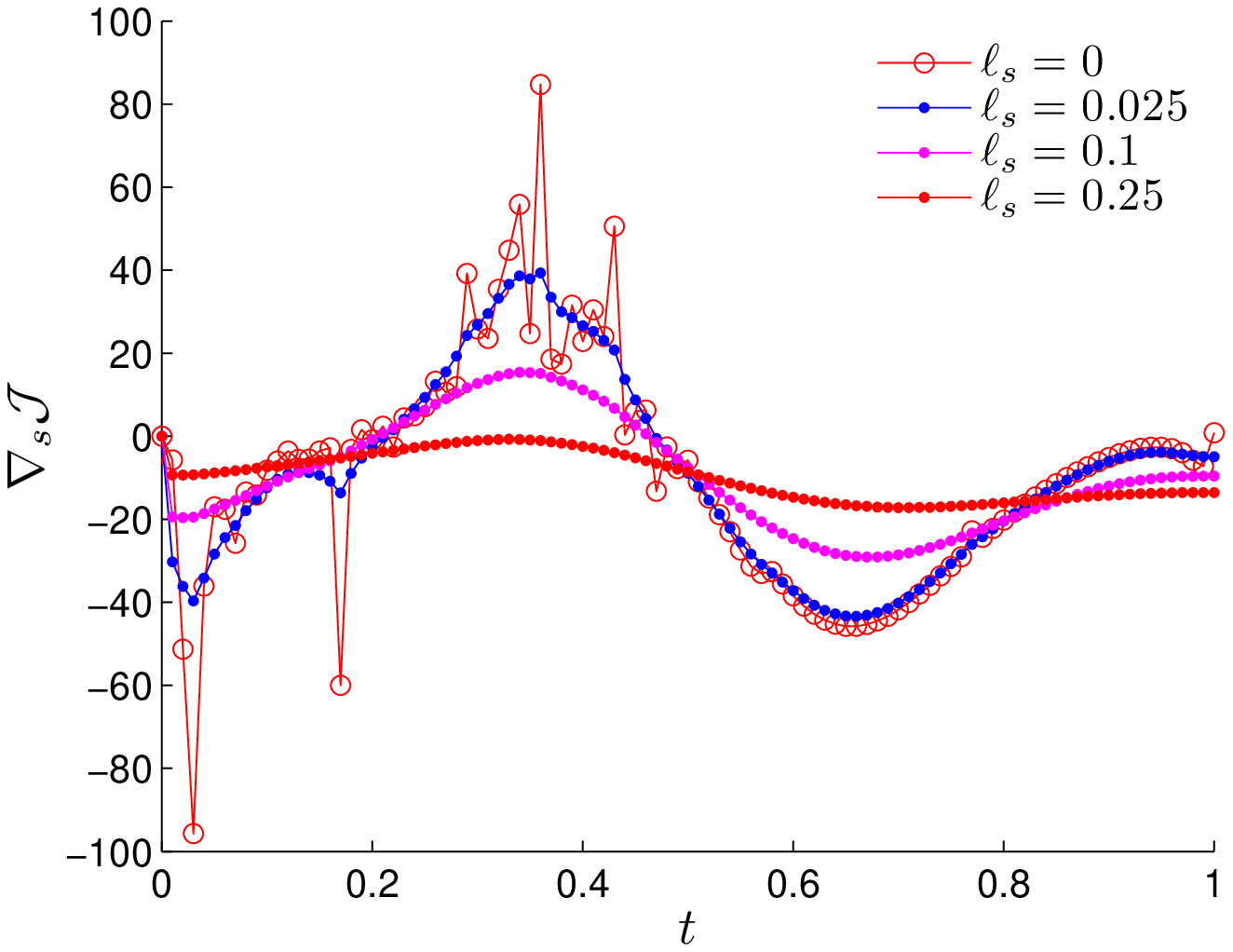}}}
  \end{center}
  \caption{Comparison of (red circles) the $L_2$ gradients $\bnabla^{L_2}_s \mJ(s_k)$ and the Sobolev
    gradients $\bnabla^{H^1}_s \mJ(s_k)$ defined in \eqref{eq:helm} for different values of smoothing
    coefficient (blue dots) $\ell_s = 0.025$, (purple dots) $\ell_s = 0.1$, and (red dots) $\ell_s = 0.25$
    obtained for model \#3 when (a) $k = 1$ and (b) $k = 38$.}
  \label{fig:grad_precond}
\end{figure}

In Figure~\ref{fig:grad_precond} we present the original gradient $\J_s'(v)$ (red circles) calculated according to \eqref{eq:grads0}, \eqref{eq:grads}, and $H^1$ (blue, purple and red dots)
gradient $\bnabla^{H^1}_{s} \J(v)$ which solves \eqref{eq:helm}, obtained for model \#3 at the first iteration, $k = 1$, and when the termination condition is reached,
$k = 38$. In the first place, we observe that the gradient $\J_s'(s_1,g)$ exhibits a smooth shape except
the small parts which are close to the endpoints $t = 0$ and $t = T$. It is explained by the fact that the initial
guess is a smooth (linear) function, see Figure~\ref{fig:models}(c), but $\J_s'(s_1,g)|_{t=0}$ has to be set
to zero as $s(0)$ is fixed and $\bnabla_{s(T)}^{H^1} \mJ(v)$ is computed by approximating Dirac measure $\delta_T$ with time grid function equal to $\tau^{-1}$ at $T$, and zero in all other grid points. The irregularity,
seen initially at the endpoints only, then tends to propagate deeper into $t$-domain and, as was anticipated in
Section~\ref{sec:approaches}, gradient $\J_s'(v)$ loose necessary smoothness which makes them unsuitable to step
forward in the optimization process. On the other hand, the gradients extracted in the Hilbert-Besov space $H^1$ are
characterized by the required smoothness which will be used for the major part of our computations accompanied
by the further analysis on the proper choice of preconditioning parameter $\ell_s$.

\begin{figure}[htb!]
  \begin{center}
  \mbox{
  \subfigure[model \#1]{\includegraphics[width=0.33\textwidth]{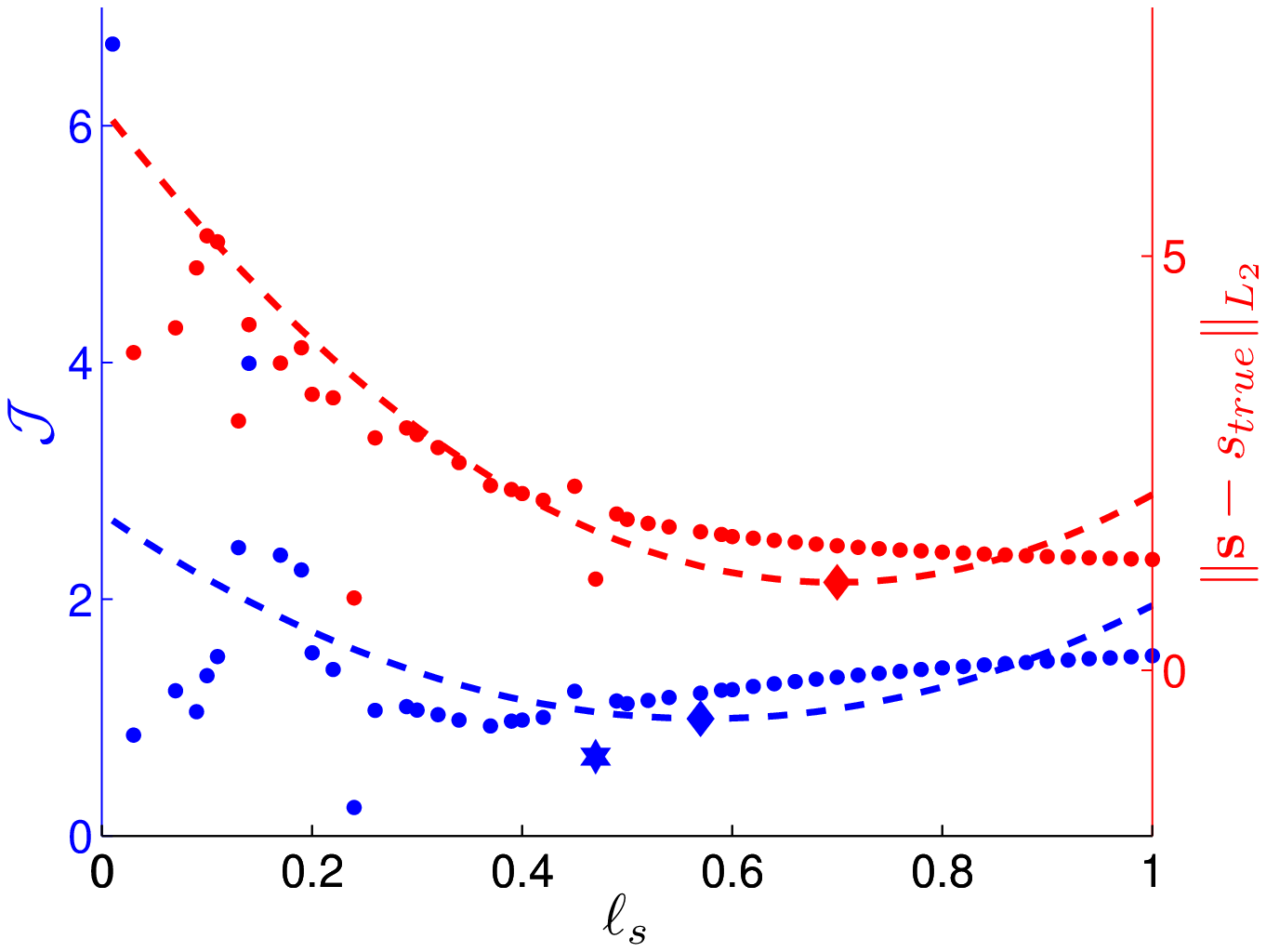}}
  \subfigure[model \#2]{\includegraphics[width=0.33\textwidth]{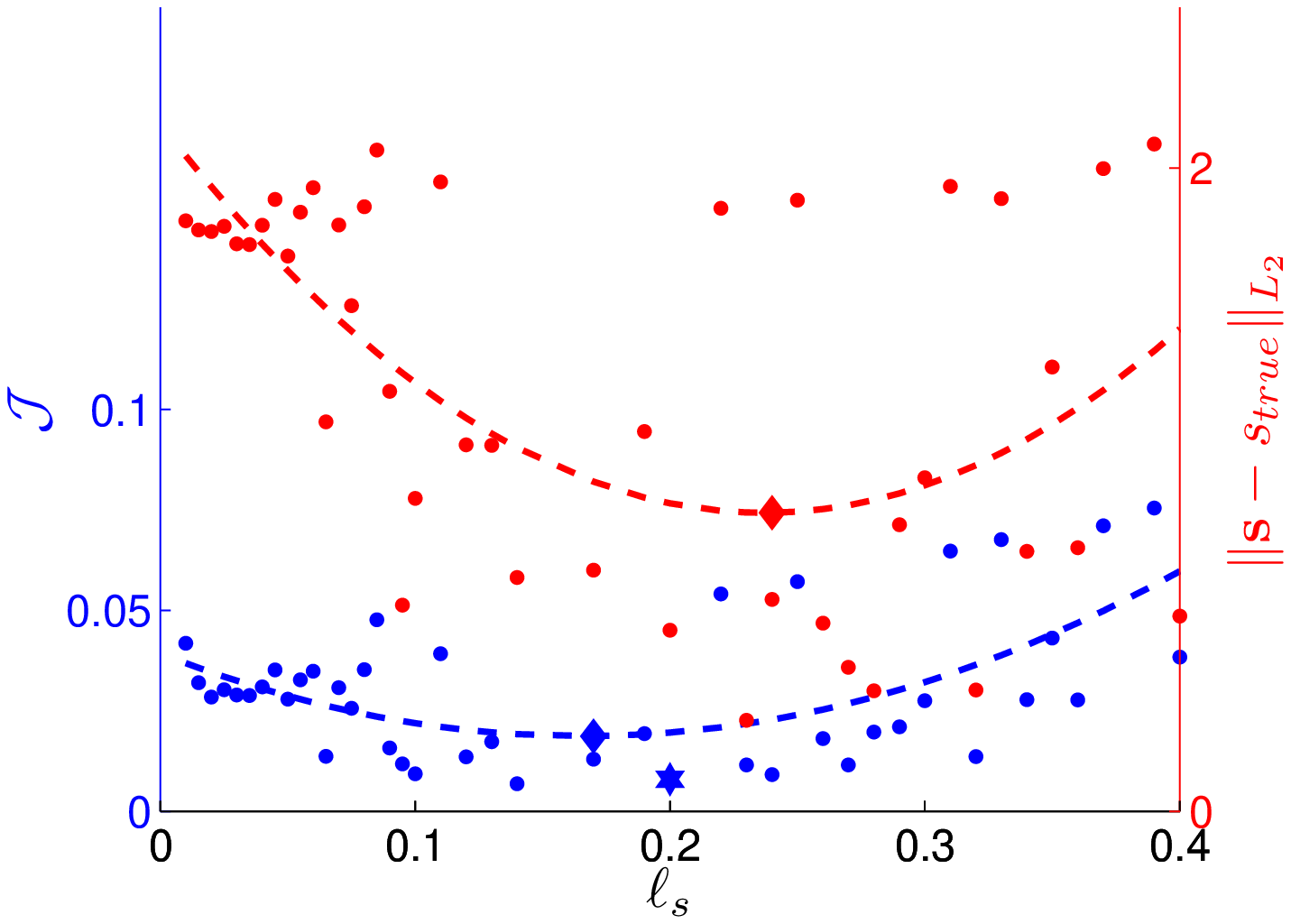}}
  \subfigure[model \#3]{\includegraphics[width=0.33\textwidth]{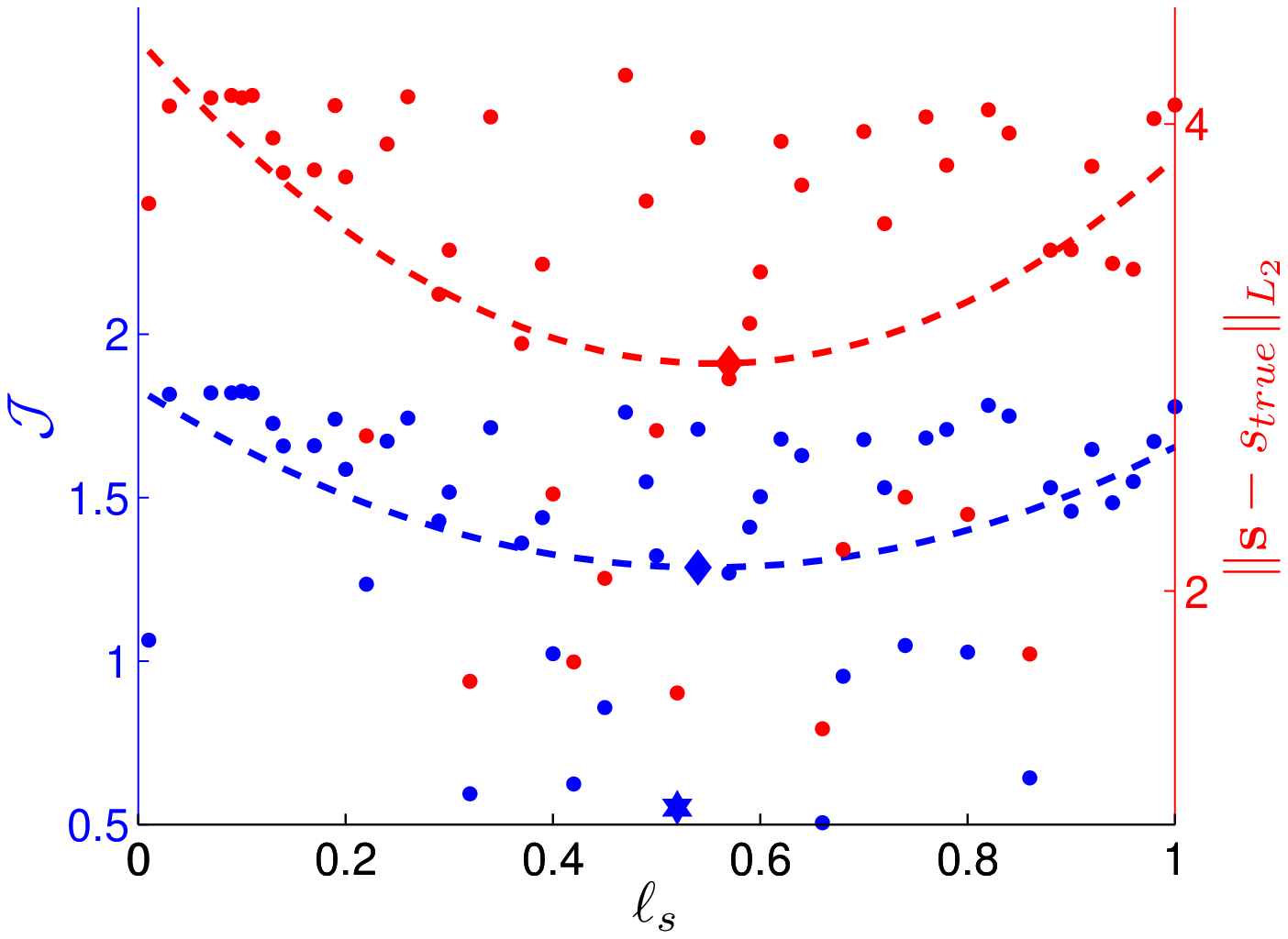}}}
  \end{center}
  \caption{Approximation of optimal parameters $\ell^*_s$ for preconditioning procedure.
    Cost functional values $\mJ$ and the solution norms $\| \bs - s_{\rm{true}} \|_{L_2}$
    are represented respectively by blue and red dots for (a) model \#1, (b) model \#2, and
    (c) model \#3. Quadratic regression models for cost functionals and solution norms
    are represented by dashed lines with minimal values shown by diamonds. The best (minimal) values
    of $\mJ$ in the proximity of approximated $\ell^*_s$ are shown by blue hexagons.}
  \label{fig:precond_s}
\end{figure}

The results of our numerical experiments confirm the fact that finding the optimal solution $\bs (t)$
is sensitive to the choice of preconditioning parameter $\ell_s$ in \eqref{eq:helm}. As noted in
Section~\ref{sec:approaches} and seen in Figure~\ref{fig:grad_precond}, small values of $\ell_s$
eliminate the difference between $\J_s'(v)$ and $\bnabla^{H^1}_s \mJ$, while large values make gradients
$\bnabla^{H^1}_s \mJ$ less informative due to ``over-smoothing''. In our strategy to find the optimal
value $\ell^*_s$, i.e.~to calibrate the preconditioning procedure, we have used two criteria. As shown
in Figure~\ref{fig:precond_s} the cost functional value $\mJ$ (blue dots) and the solution norm
$\| \bs - s_{\rm{true}} \|_{L_2}$ (red dots) are recorded after performing optimizations supplied with
different values of $\ell_s$ for each model. Both sets of points are then used to perform the least
square analysis to find the quadratic regression model (dashed lines) for each set. The quadratic
functions to model $\mJ$ are then minimized to approximate $\ell^*_s$ (blue diamonds) giving values
$\ell^*_{s,1} = 0.57$, $\ell^*_{s,2} = 0.17$ and $\ell^*_{s,3} = 0.54$ correspondingly for models
\#1, \#2 and \#3. The quadratic functions to model solution norms are also minimized to confirm the
proximity of the obtained solutions (red diamonds) to approximated $\ell^*_s$. Although the second
criterion in many cases is not available, here we use it to demonstrate the consistence of the
results obtained by both of them. Unless otherwise stated, all the computational results discussed in this
section will implement the above mentioned optimal values $\ell^*_{s,i}, \, i = 1, 2, 3$, whenever preconditioning procedure
is active for $s(t)$.

\begin{figure}[htb!]
  \begin{center}
  \mbox{
  \subfigure[model \#1]{\includegraphics[width=0.33\textwidth]{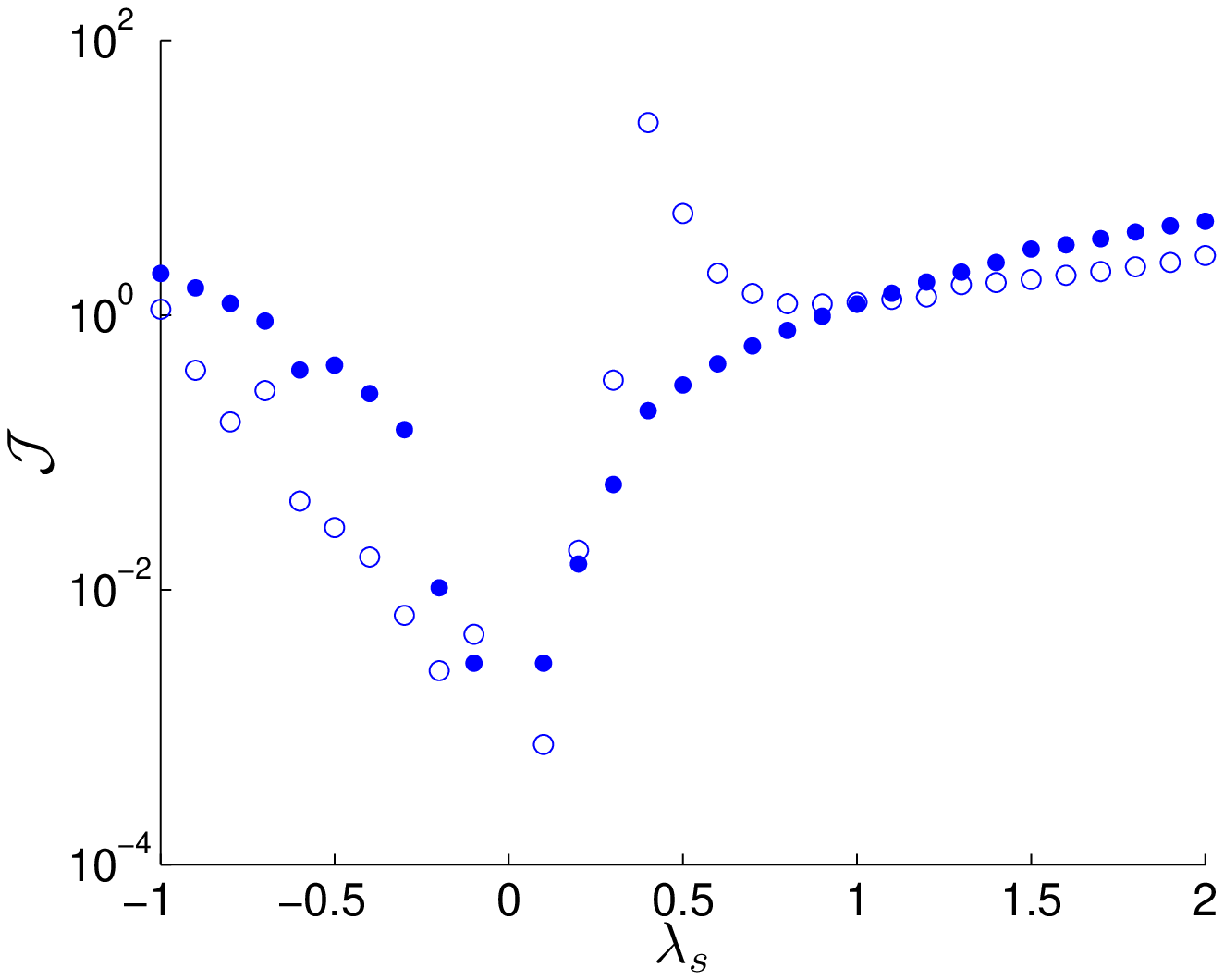}}
  \subfigure[model \#2]{\includegraphics[width=0.33\textwidth]{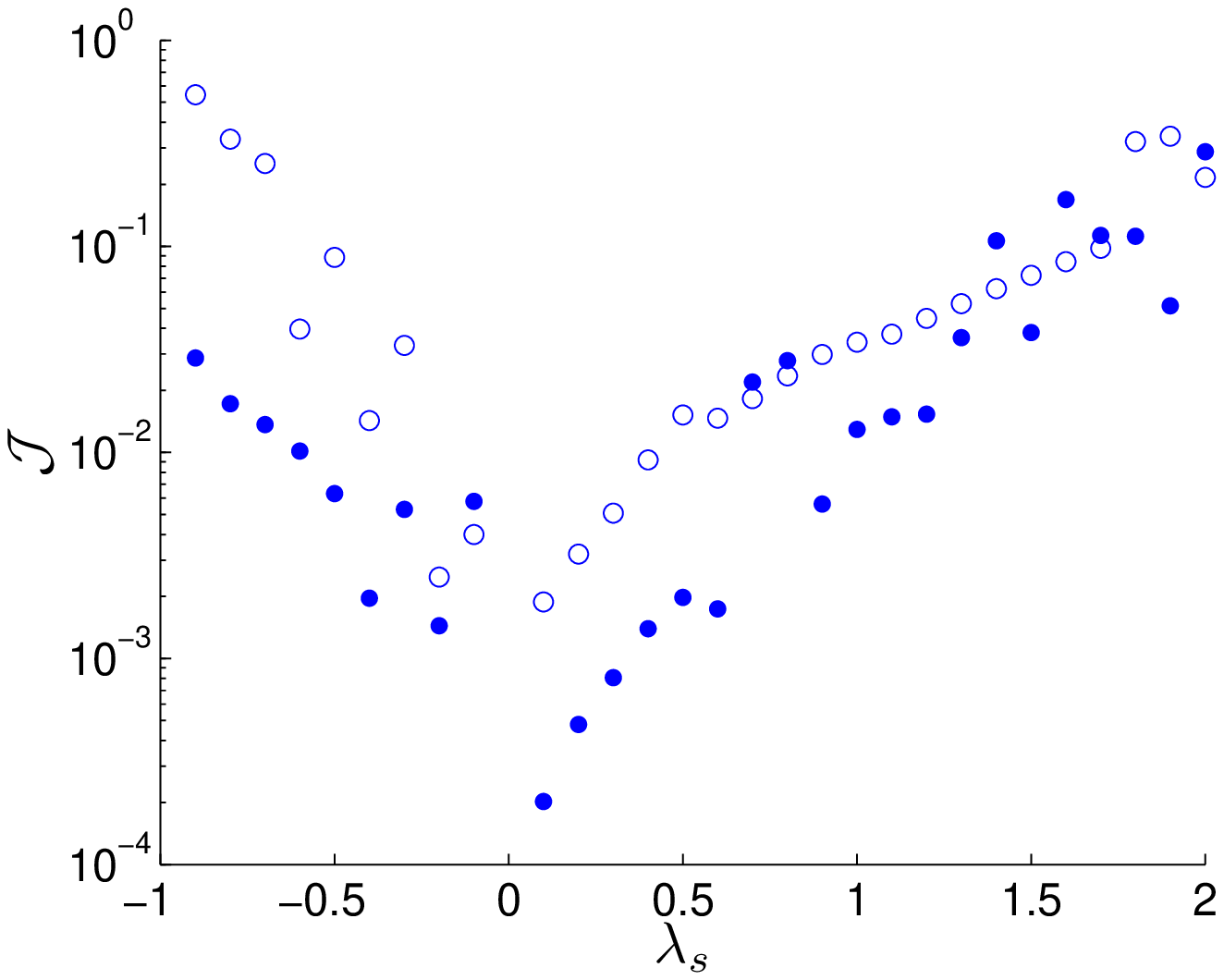}}
  \subfigure[model \#3]{\includegraphics[width=0.33\textwidth]{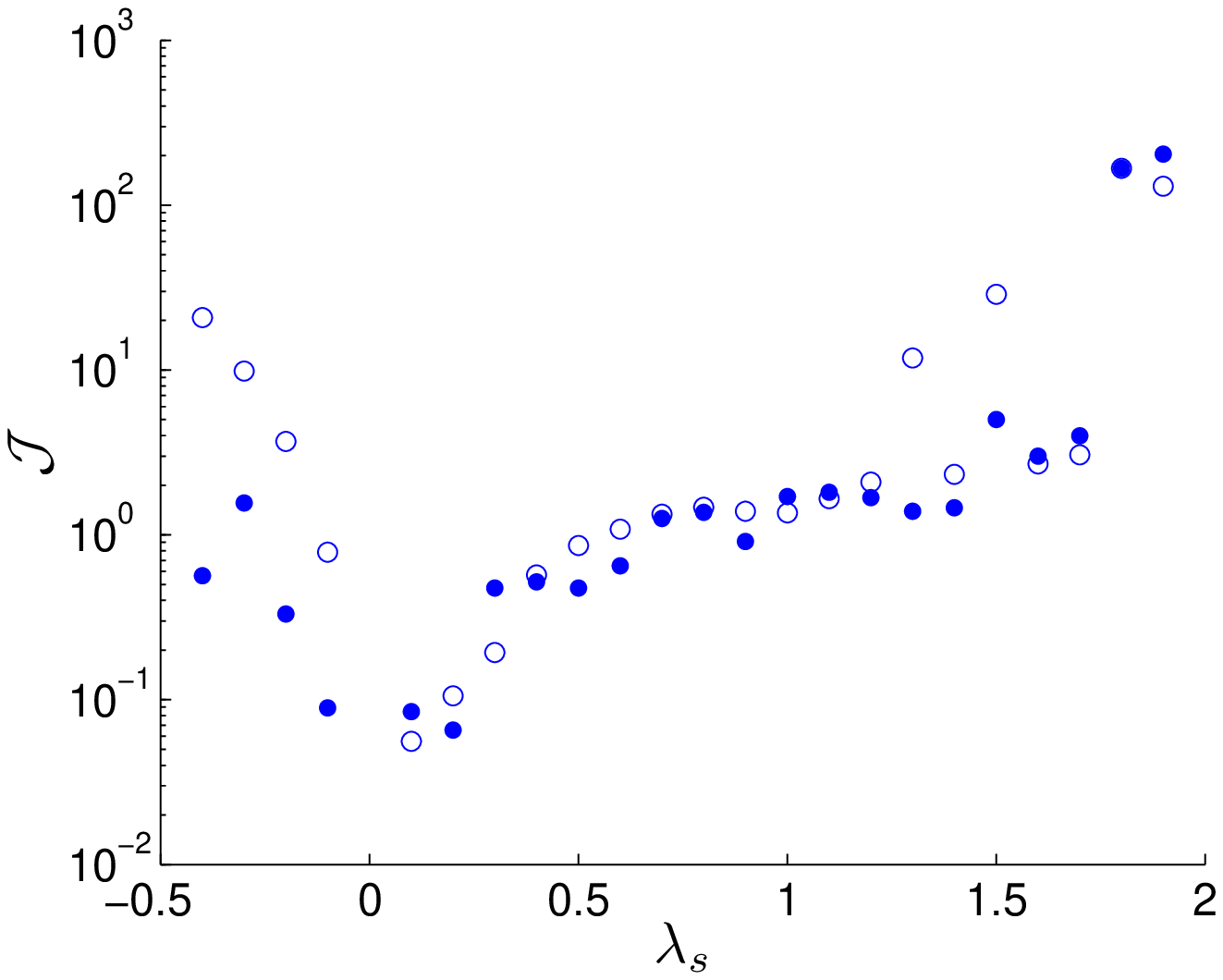}}}
  \mbox{
  \subfigure[model \#1]{\includegraphics[width=0.33\textwidth]{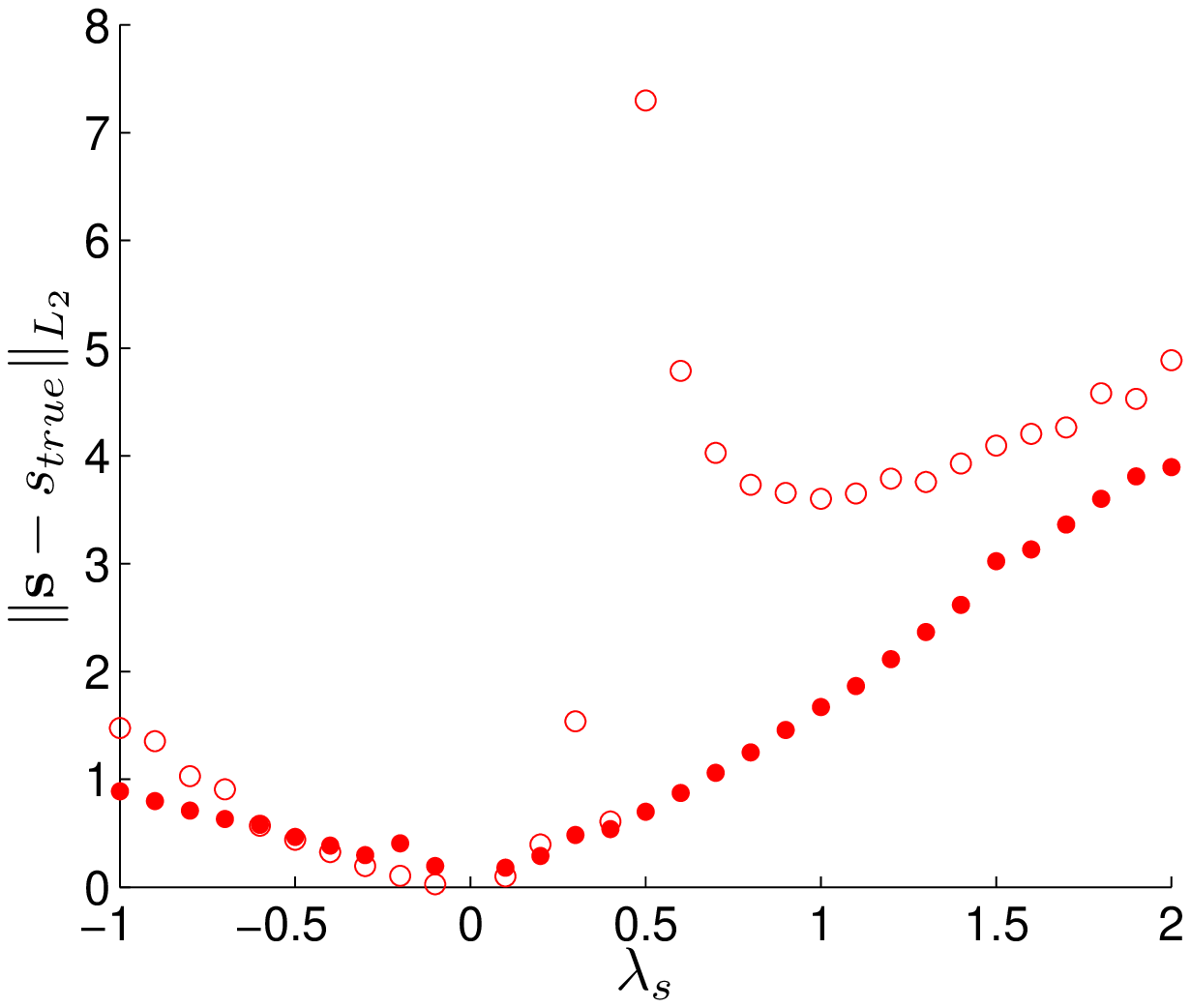}}
  \subfigure[model \#2]{\includegraphics[width=0.33\textwidth]{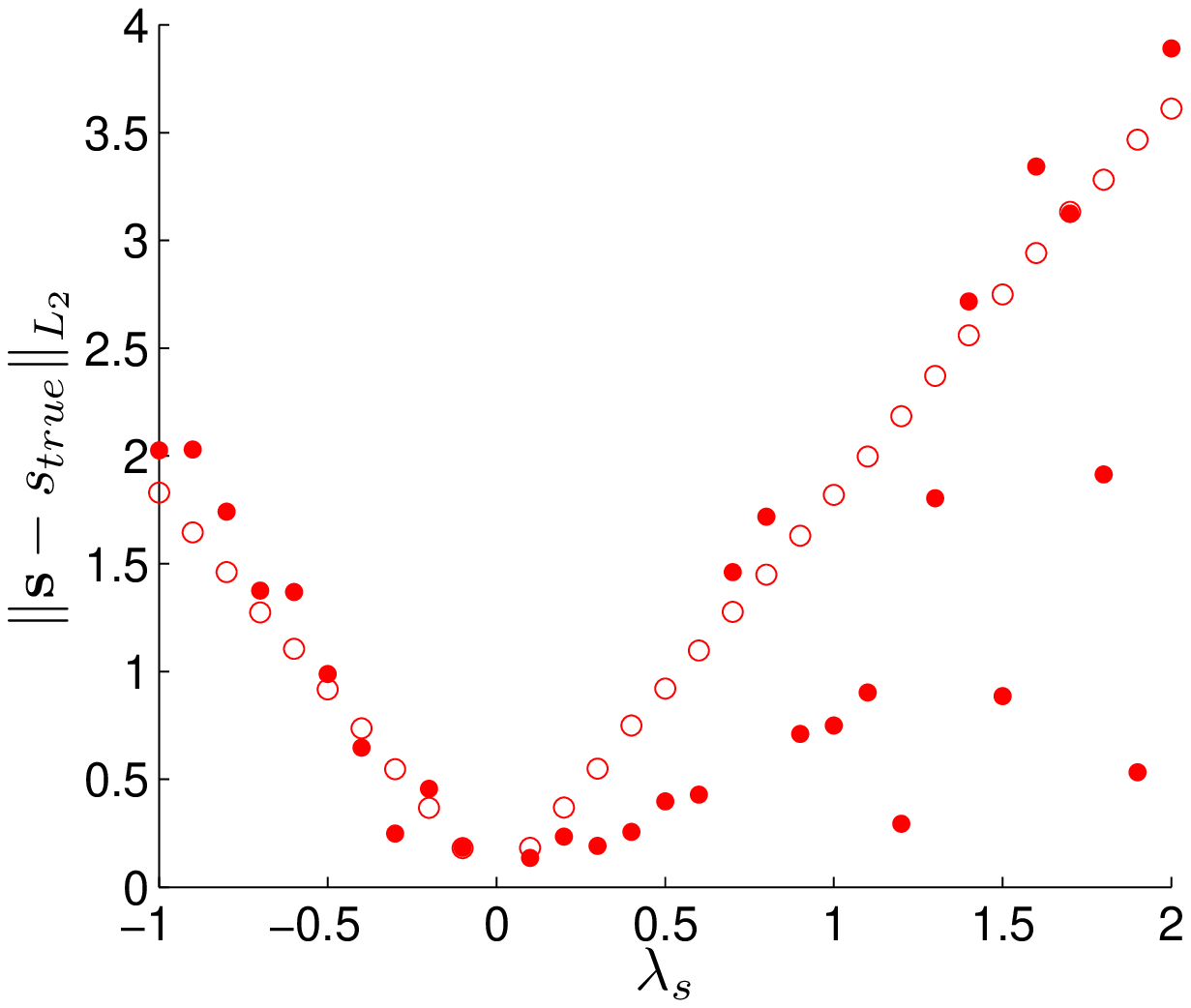}}
  \subfigure[model \#3]{\includegraphics[width=0.33\textwidth]{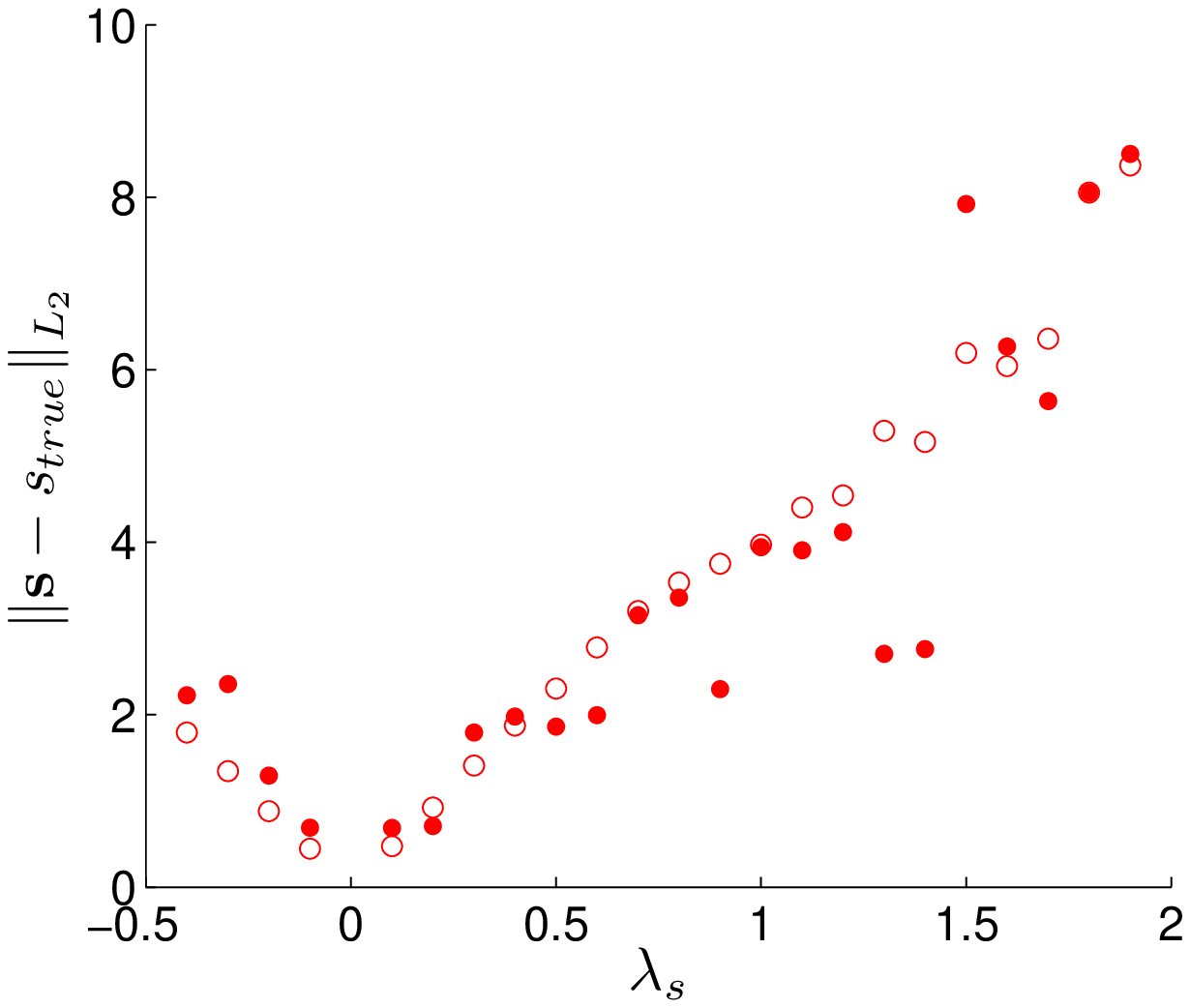}}}
  \end{center}
  \caption{Convergence analysis for $s(t)$ performed for (a,d) model \#1, (b,e) model \#2, and
    (c,f) model \#3 by evaluating (a,b,c) cost functional values $\mJ$ and (d,e,f)
    solution norms $\| \bs - s_{\rm{true}} \|_{L_2}$. Dots and circles represent results
    obtained correspondingly with and without preconditioning procedure \eqref{eq:helm}.}
  \label{fig:conv_s}
\end{figure}

We would like to reiterate that, since the inverse Stefan problem \eqref{eq:pde-1}--\eqref{eq:pde-finaltemp}
is in general nonconvex, Algorithm~\ref{alg:gen_opt} is able to find a local, rather than global,
optimal solution. To further validate our computational approach in terms of convergence to the global
optimal solution, we solve the same optimization problem for all three models starting with different
initial guesses. For consistence, these new initial guesses $s_{\rm{ini},\lambda_s}$ are parameterized
with respect to their proximity to global minimizer $s_{\rm{true}}$ in the following way
\begin{equation}
  s_{\rm{ini},\lambda_s}(t) = (1 - \lambda_s) s_{\rm{true}}(t) + \lambda_s s_{\rm{ini}}(t).
  \label{eq:s_convex_comb}
\end{equation}
We note that setting parameter $\lambda_s = 1$ recovers the regular initial guess shown in Figure~\ref{fig:models}(a,b,c),
while $\lambda_s \rightarrow 0$ moves initial guess in the close neighborhood of $s_{\rm{true}}$.

The results of the convergence test for $\lambda_s \in [-1, 2]$ (models \#1 and \#2) and $\lambda_s \in [-0.5, 2]$
(model~\#3) are shown in Figure~\ref{fig:conv_s} for cases with and without preconditioning procedure
\eqref{eq:helm} by evaluating both cost functional values $\mJ$ and solution norms $\| \bs - s_{\rm{true}} \|_{L_2}$.
As expected, our results for all three models show good convergence to global minimizer $s_{\rm{true}}$,
i.e. $\mJ \rightarrow \mJ_{\rm{min}}$ and $\| \bs - s_{\rm{true}} \|_{L_2} \rightarrow 0$ as $\lambda_s \rightarrow 0$.
We could also conclude that applying preconditioning in general benefits in improving this convergence.

\begin{figure}[htb!]
  \begin{center}
  \mbox{
  \subfigure[model \#1: $s(t)$]{\includegraphics[width=0.33\textwidth]{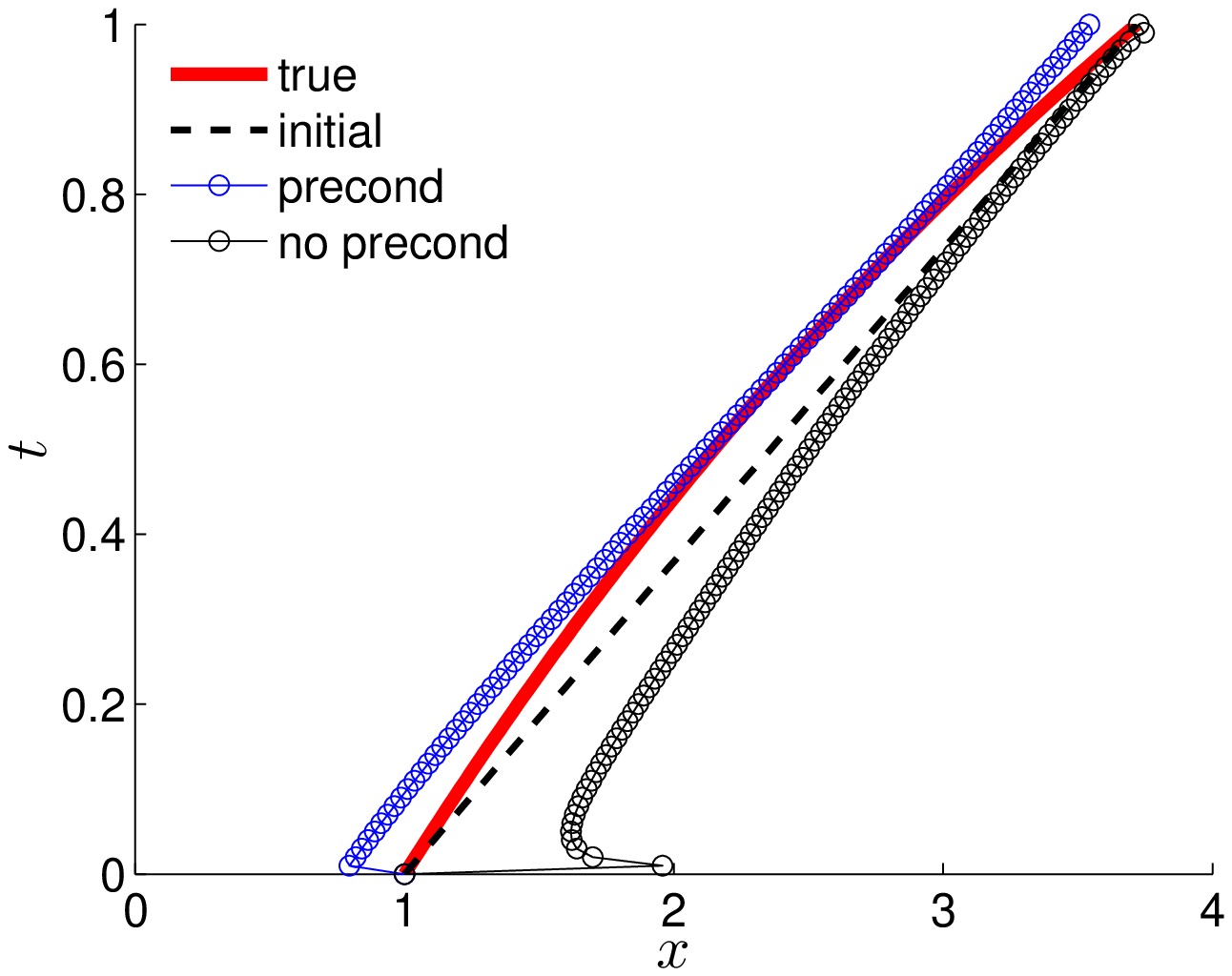}}
  \subfigure[model \#2: $s(t)$]{\includegraphics[width=0.33\textwidth]{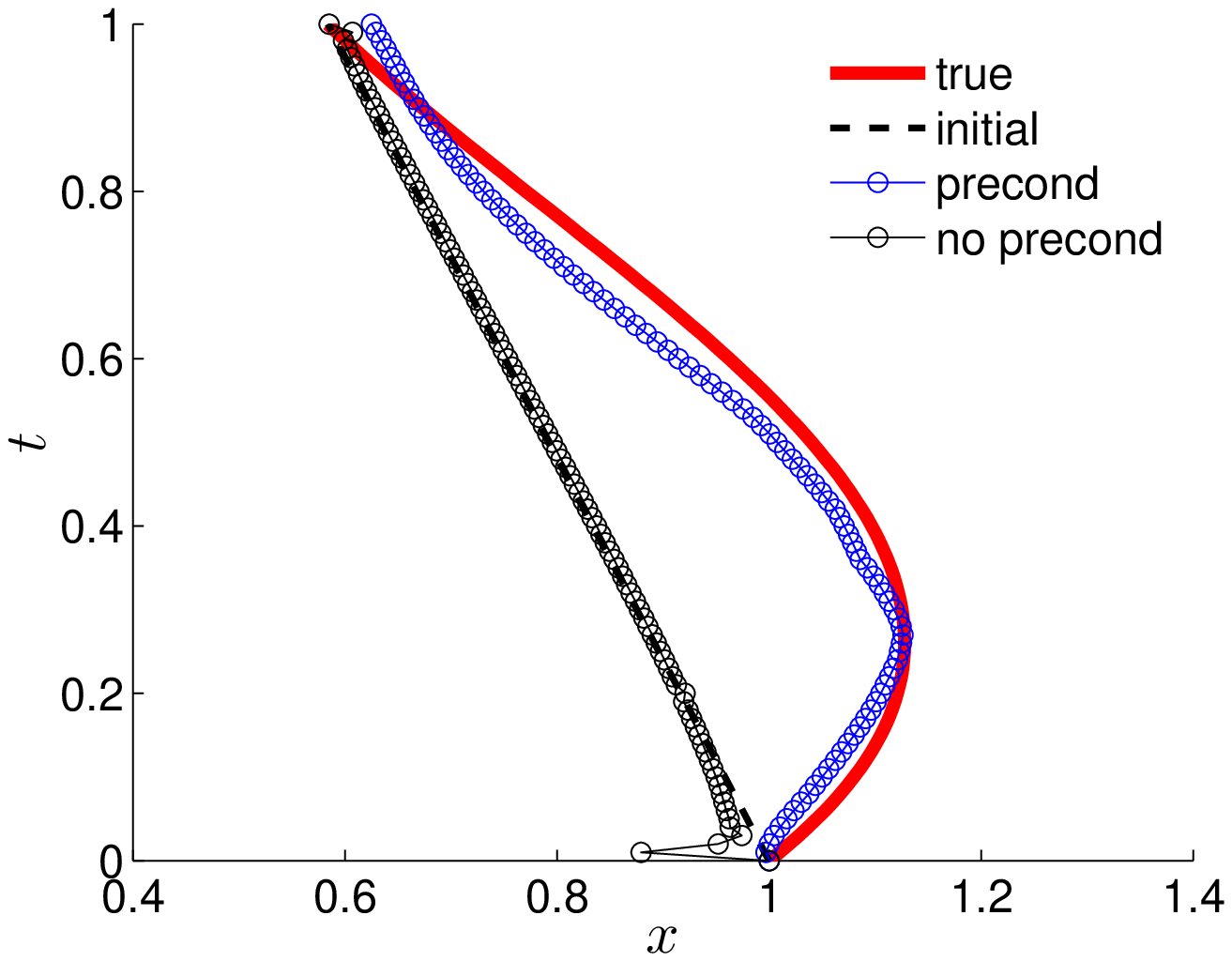}}
  \subfigure[model \#3: $s(t)$]{\includegraphics[width=0.33\textwidth]{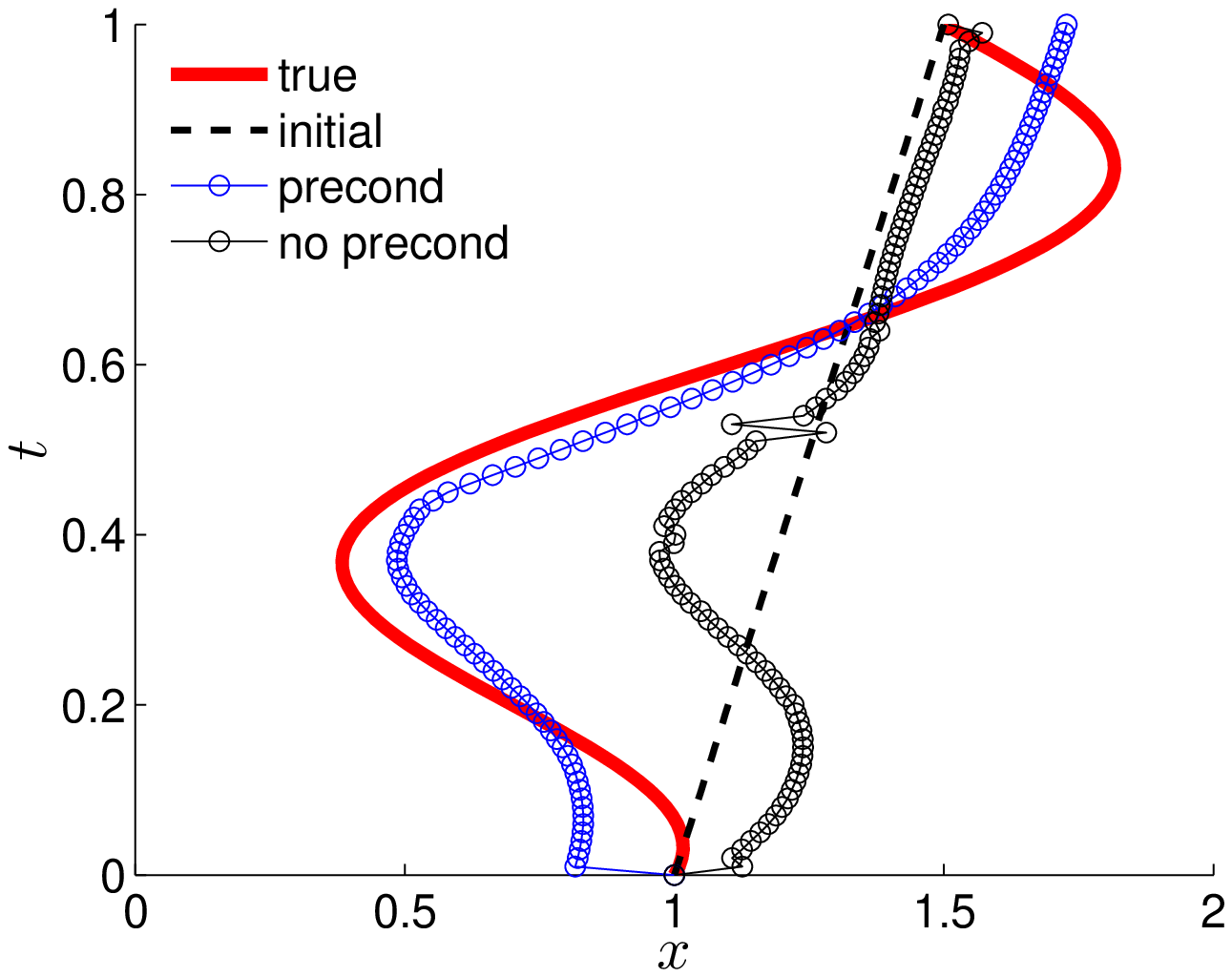}}}
  \mbox{
  \subfigure[model \#1: $w(x) \ \& \ \mu(t)$]{\includegraphics[width=0.33\textwidth]{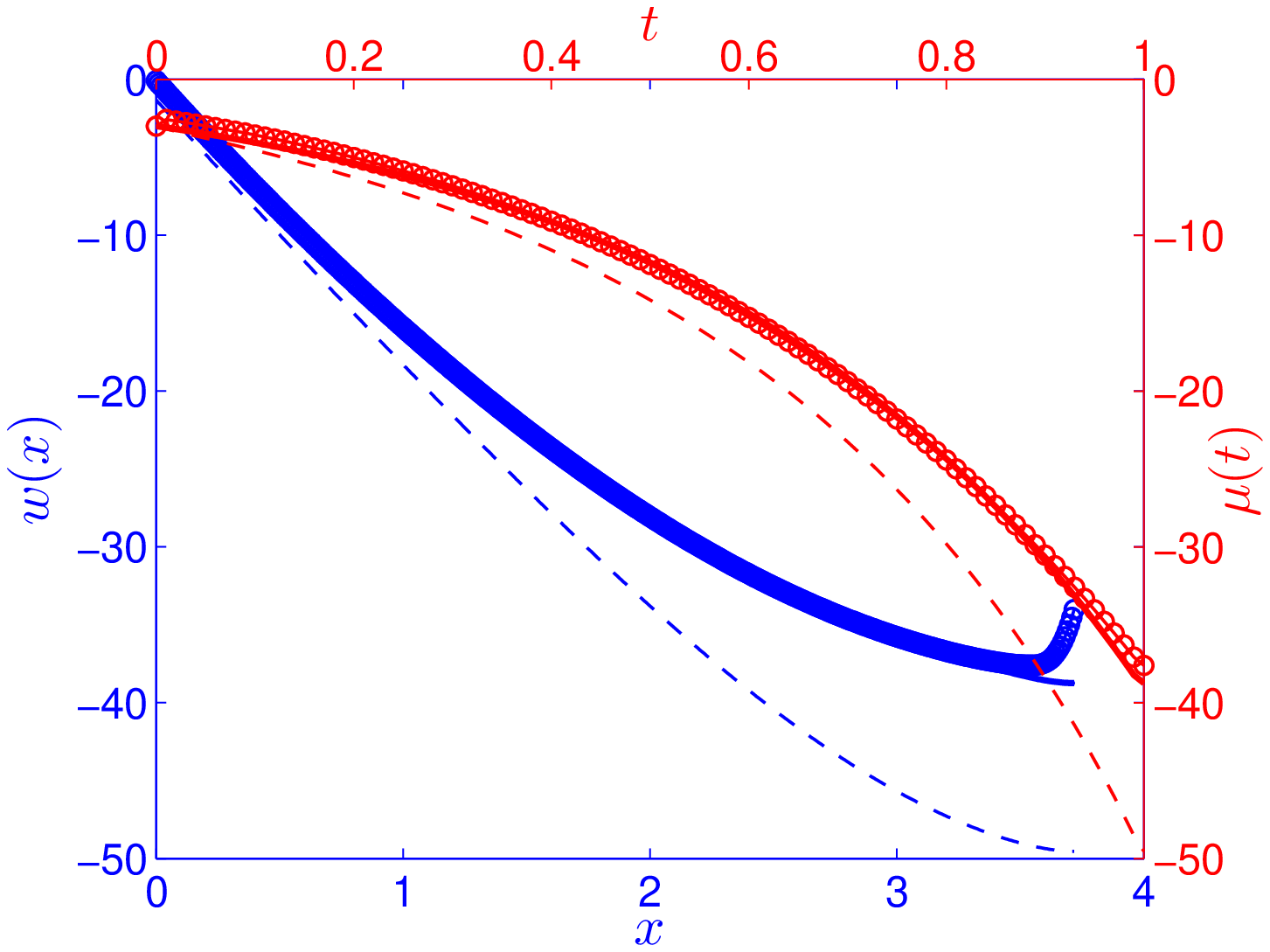}}
  \subfigure[model \#2: $w(x) \ \& \ \mu(t)$]{\includegraphics[width=0.33\textwidth]{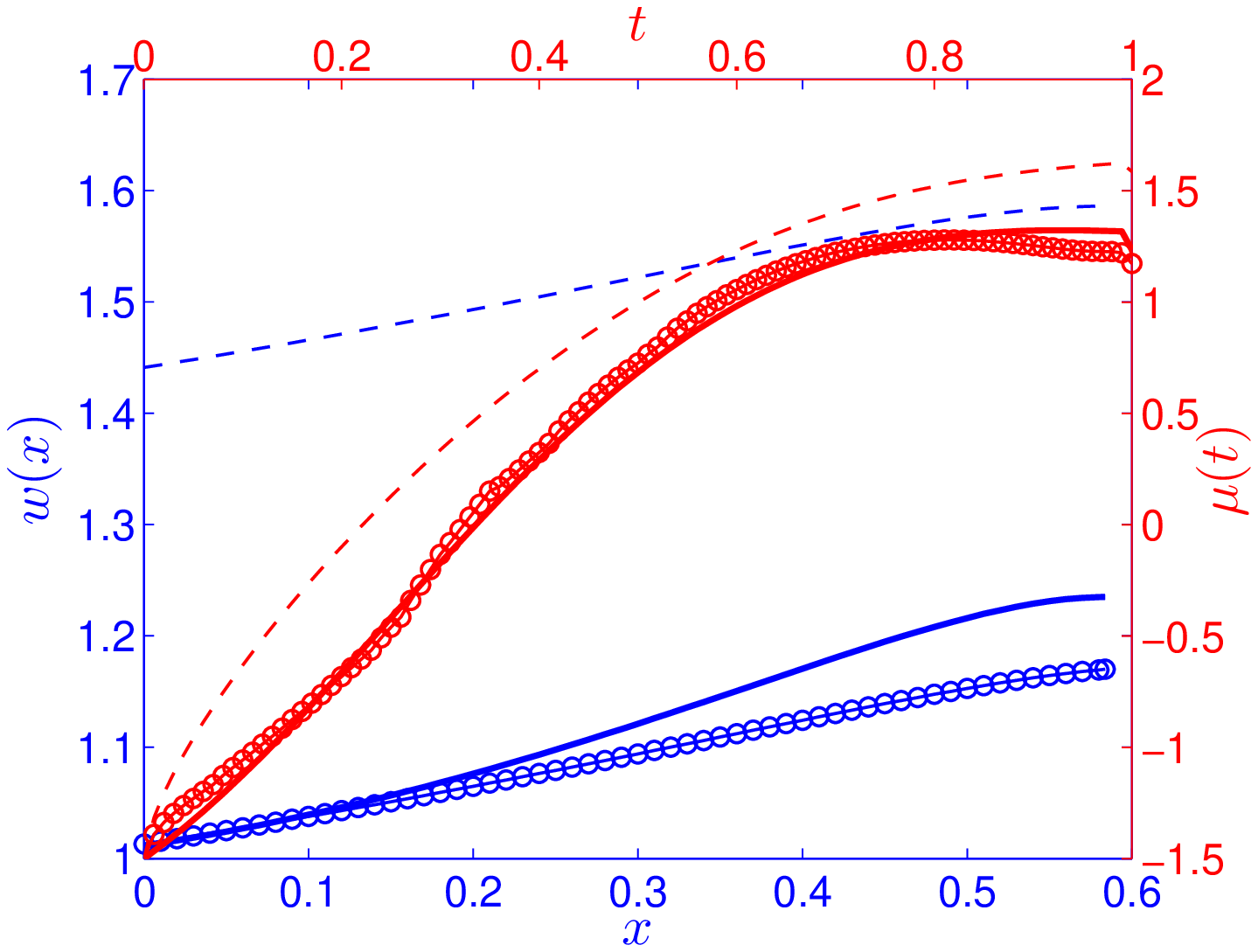}}
  \subfigure[model \#3: $w(x) \ \& \ \mu(t)$]{\includegraphics[width=0.33\textwidth]{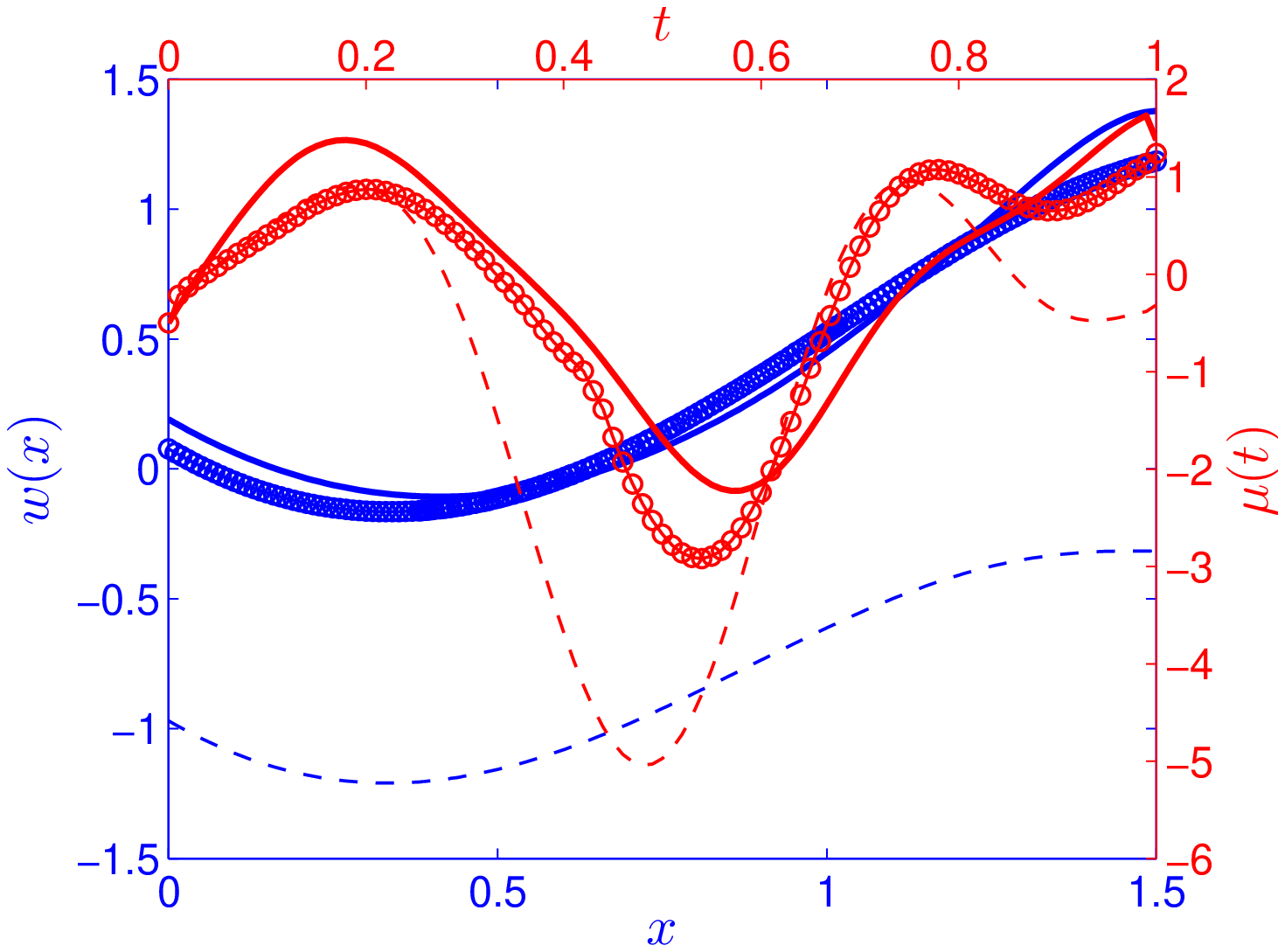}}}
  \mbox{
  \subfigure[model \#1: $\mJ_k$]{\includegraphics[width=0.33\textwidth]{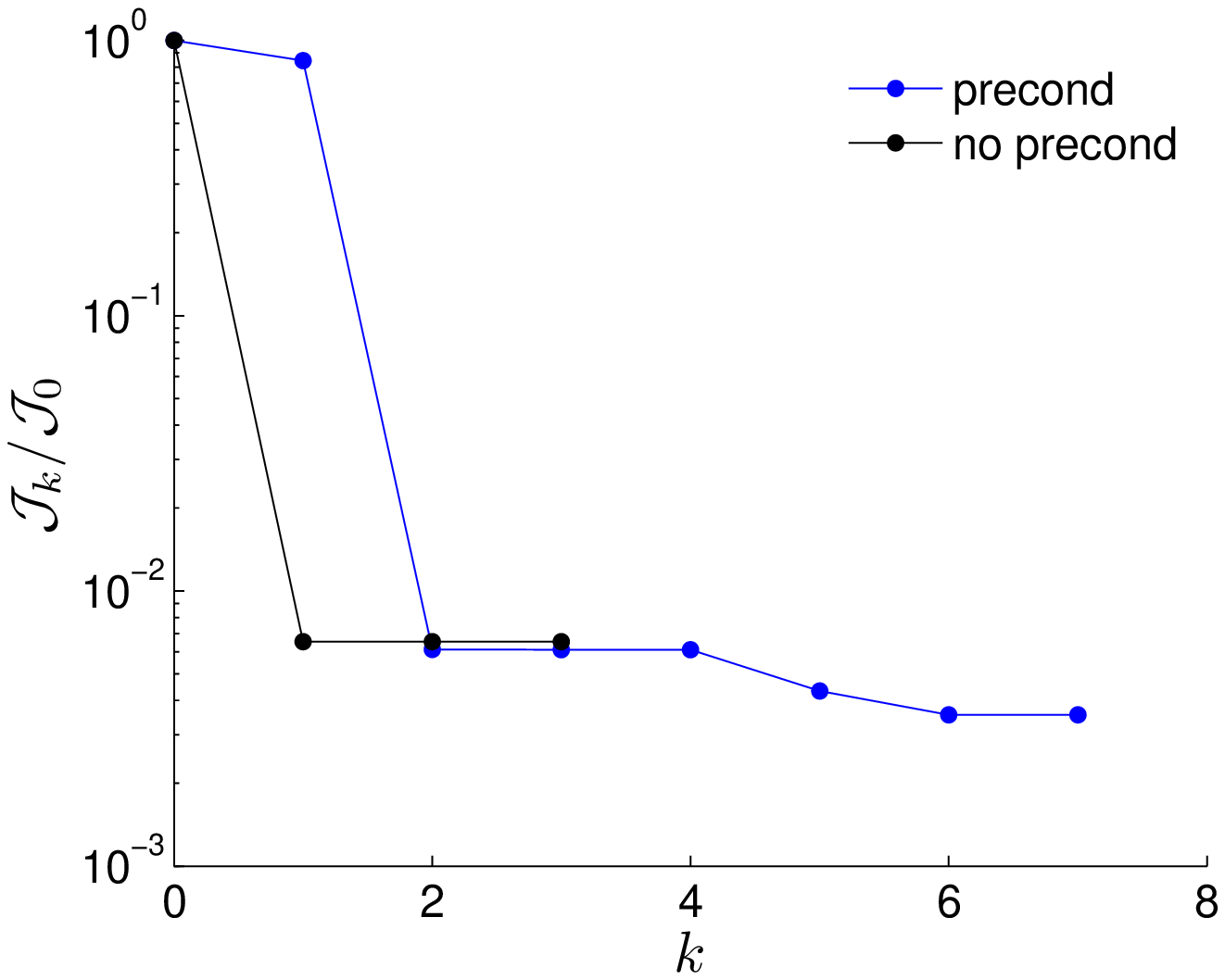}}
  \subfigure[model \#2: $\mJ_k$]{\includegraphics[width=0.33\textwidth]{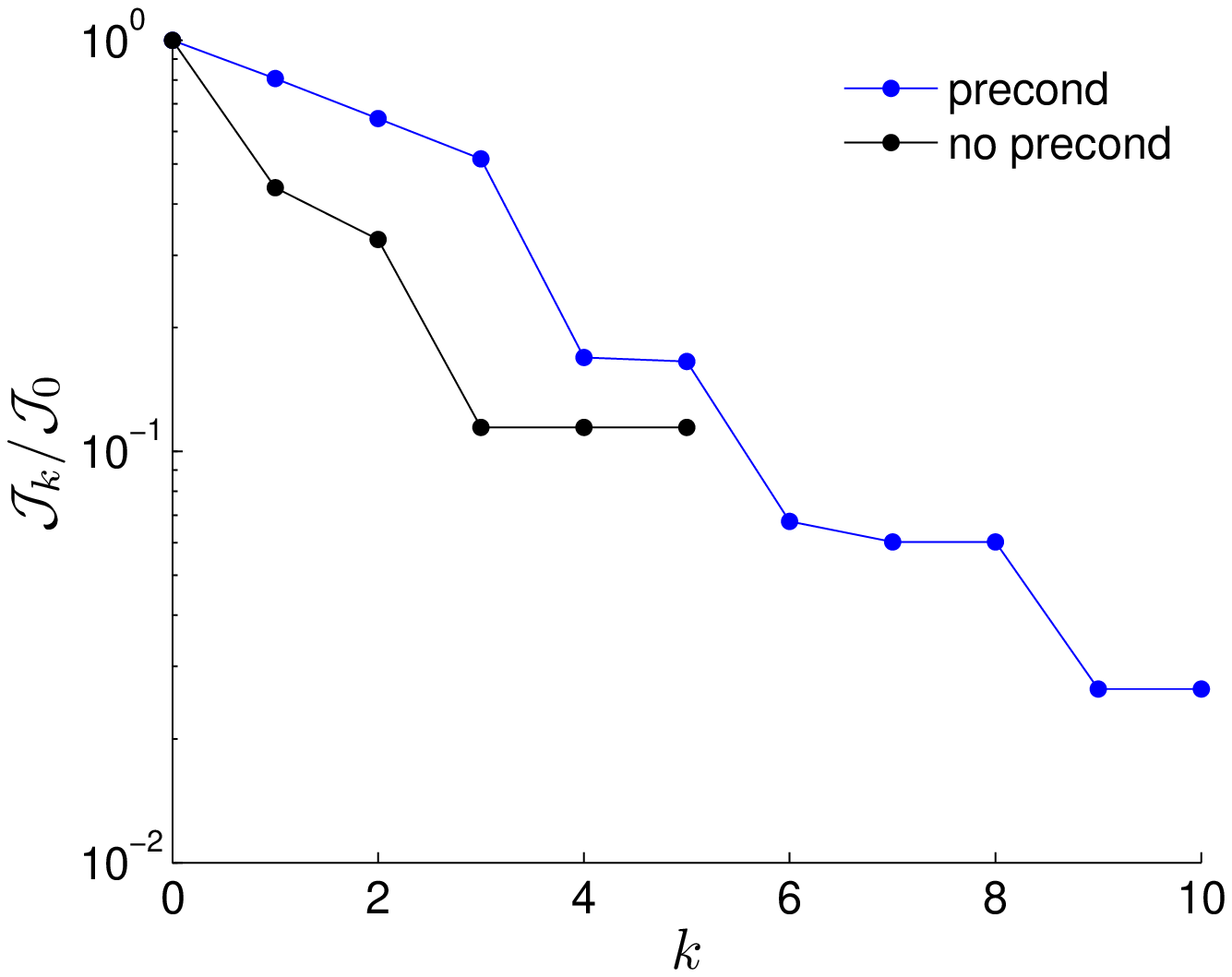}}
  \subfigure[model \#3: $\mJ_k$]{\includegraphics[width=0.33\textwidth]{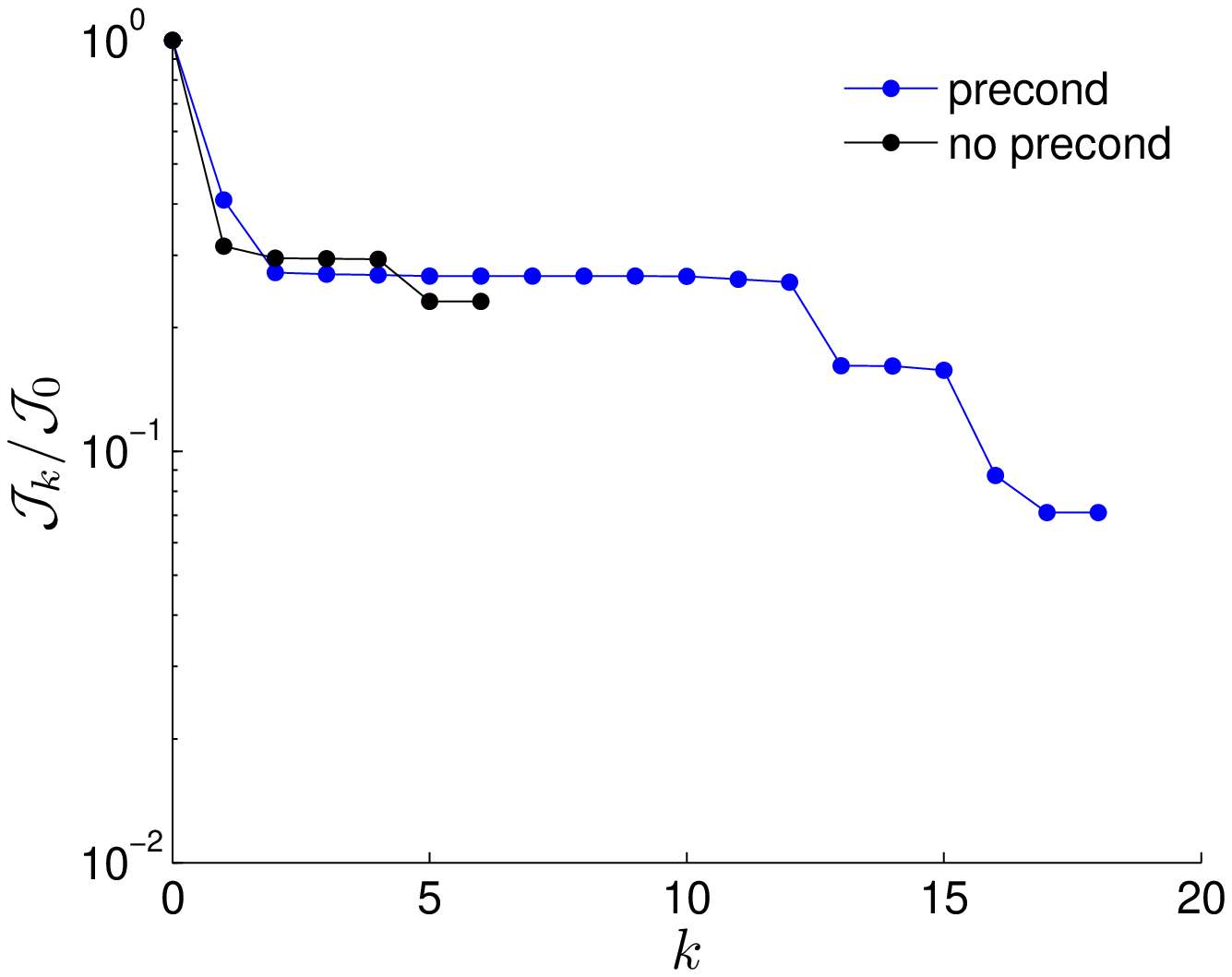}}}
  \end{center}
  \caption{Results of reconstructing free boundary $s(t)$ for (a,d,g) model \#1, (b,e,h) model \#2, and
    (c,f,i) model \#3. In (a-c) solid red and dashed black lines show respectively the shape of $s_{\rm{true}}(t)$
    and initial guess $s_{\rm{ini}}(t)$, while blue and black circles represent optimal solutions $\bs(t)$
    respectively with and without preconditioning. In (d-f) blue and red colors are used correspondingly
    for functions $w(x)$ and $\mu(t)$: dashed and solid lines represent respectively their
    values when $s = s_{\rm{ini}}(t)$ and $s = s_{\rm{true}}(t)$, while circles are used when
    $s = \bs(t)$ obtained with preconditioning. In (g-i) blue and black dots show normalized cost functionals
    $\mJ_k/\mJ_0$ as functions of iteration number $k$.}
  \label{fig:opt_ell_s}
\end{figure}

Finally, Figure~\ref{fig:opt_ell_s}(a-c) shows the outcomes of reconstructing free boundary $s(t)$ for all three
models comparing the results obtained with and without preconditioning (blue and black circles respectively).
Preconditioning procedure uses $\ell^{\star}_{s,1} = 0.47$, $\ell^{\star}_{s,2} = 0.2$ and $\ell^{\star}_{s,3} = 0.52$
obtained by finding the minimal values of $\mJ$ (shown by blue hexagons in Figure~\ref{fig:precond_s}) in the
proximity of approximated $\ell^*_s$. The superior quality of reconstruction of $s(t)$ in the preconditioned
case is obvious and it is also justified by observing how accurately the obtained solutions $u(x,T)$ and $u(s(t),t)$
match the measurements $w(x)$ and $\mu(t)$ which is seen in Figure~\ref{fig:opt_ell_s}(d-f).
In Figure~\ref{fig:opt_ell_s}(g-i) normalized cost functionals $\mJ_k/\mJ_0$
are represented as functions of iteration number $k$. As could be noted here, cases with active preconditioning
are prone to run at least two times longer with higher chances to find a ``better'' local optimizer.

\subsection{Reconstruction in the presence of noise}
\label{sec:noise}

In this section we discuss the important issue of reconstructing the free boundary $s(t)$ in the presence of noise.
This noise is incorporated into the additional measurements $\{s(t_i) \}_{i=1}^M$ made for the position of free
boundary $s(t)$ at time $t_i$ (represented by black filled circles in Figure~\ref{fig:discr_tx}(b)). In fact,
for our numerical experiments in Section~\ref{sec:single_rec} we use $M = 1, t_1 = T$ with a reference to a single
measurement $s(T) = s_{\rm true}(T)$ to create regular initial guess $s_{{\rm ini}}(t)$ as shown by black solid
line in Figure~\ref{fig:discr_tx}(b) and dashed lines in Figure~\ref{fig:models}(a,b,c). In case $M > 1$, we
assume that additional direct measurements of $s(t)$ are available and made by dividing time interval $[0, T]$
uniformly. Figure~\ref{fig:discr_tx}(b) also shows schematically the case with $M = 3$ representing three measured
values of $s(t)$ (with no noise incorporated) by three black filled circles on a red line.

To incorporate noise, say of $\eta\%$, into the measurements $\{s(t_i) \}_{i=1}^M$, we replace these measurements
at time instances $t_i, i = 1, \ldots, M$ with a new set $\{\tilde s_i^{\eta} \}_{i=1}^M$, where the independent
random variables $\tilde s_i^{\eta}$ have a normal (Gaussian) distribution with the mean $s(t_i)$ and the standard
deviation $\Delta\eta = \frac{1}{M} \sum_{i=1}^M s(t_i) \cdot \frac{\eta}{100\%}$. Unless stated otherwise, in order
to be able to directly compare reconstructions from noisy measurements with different noise levels, the same noise
realization is used after rescaling to the standard deviation $\Delta\eta$. Figure~\ref{fig:discr_tx}(b) shows
schematically the case with $M = 3$ representing three values of $\tilde s_i$ (three blue circles) with some
noise incorporated into the measurements.

These measurements could be used in different ways within the computational framework discussed previously. In the
current work we use them for two purposes. In both cases we create piecewise linear approximations of $s(t)$
as shown in Figure~\ref{fig:discr_tx}(b) by blue solid and dashed lines respectively for measurements without and
with noise. These piecewise linear approximations then could be used as
\begin{description}
  \item[case \#1:] initial guess $s_{{\rm ini}} (t) = s_{{\rm ini},M} (t)$,
  \item[case \#2:] regularization centroid $\bar s(t) = \bar s_M (t)$ in regularization term of \eqref{Eq:W:1:8reg}
    while setting initial guess in a regular way, i.e.~$s_{{\rm ini}}(t) = s_{{\rm ini},1} (t)$.
\end{description}

\begin{figure}[htb!]
  \begin{center}
  \mbox{
  \subfigure[model \#2]{\includegraphics[width=0.5\textwidth]{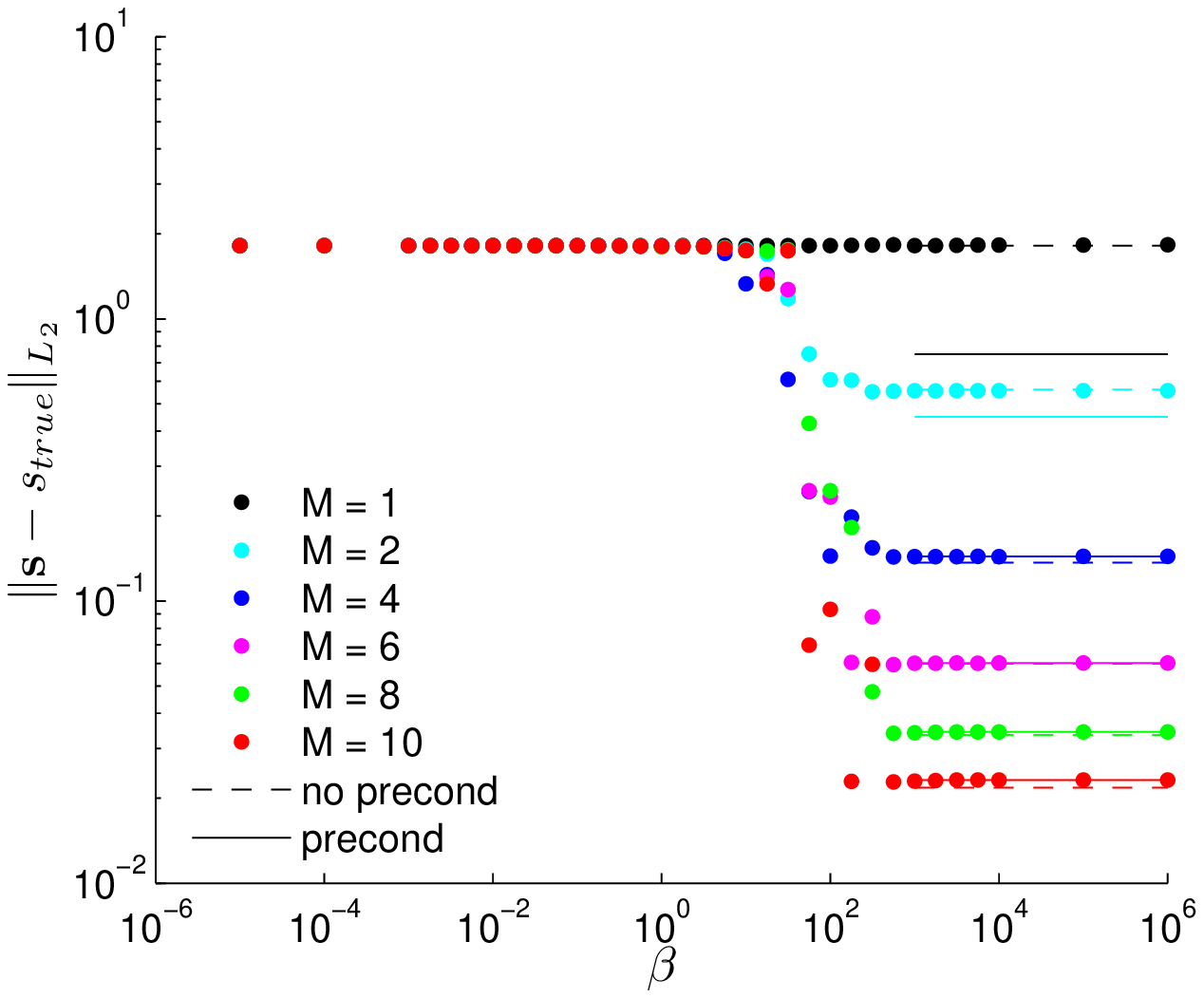}}
  \subfigure[model \#3]{\includegraphics[width=0.5\textwidth]{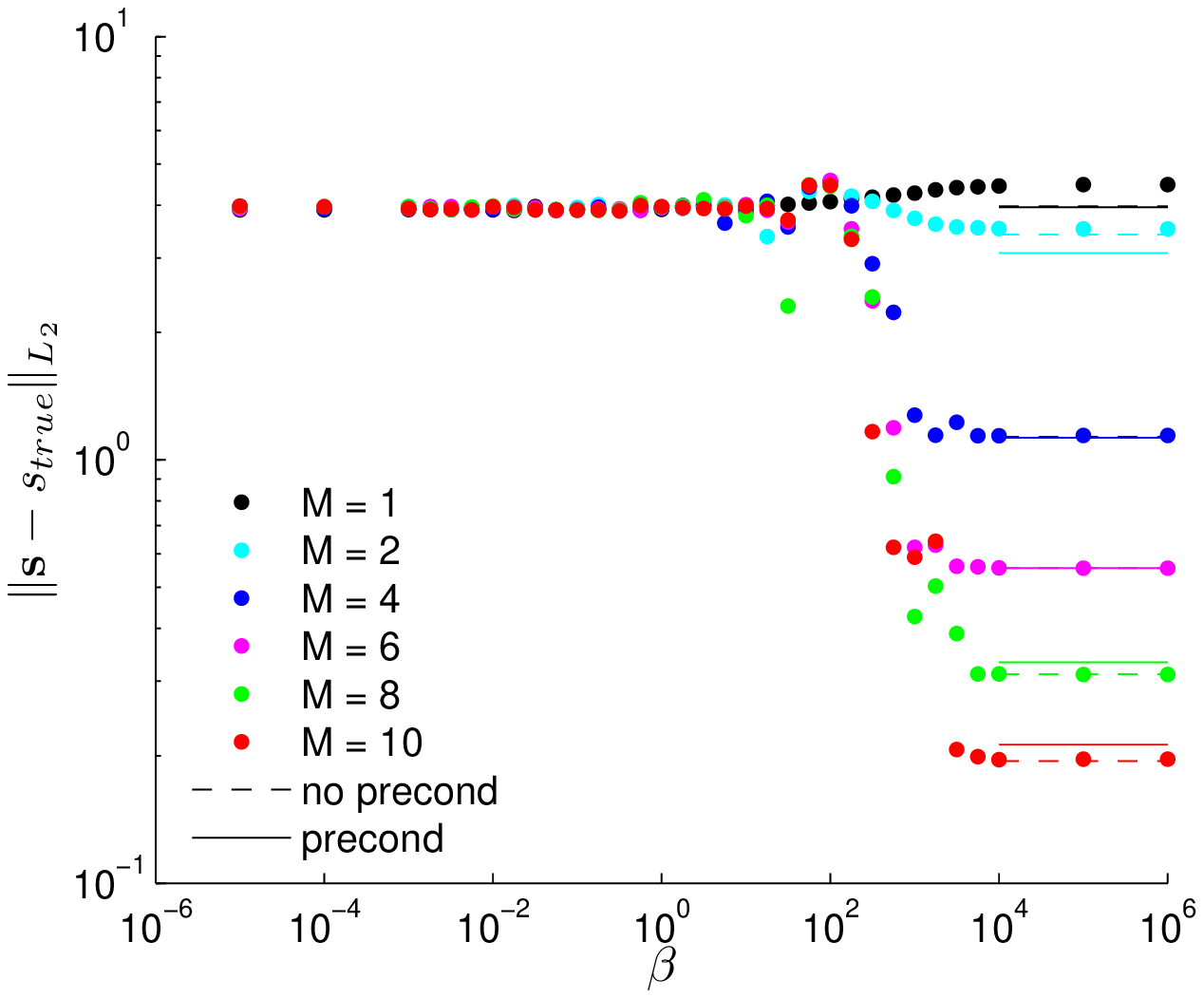}}}
  \mbox{
  \subfigure[model \#2]{\includegraphics[width=0.5\textwidth]{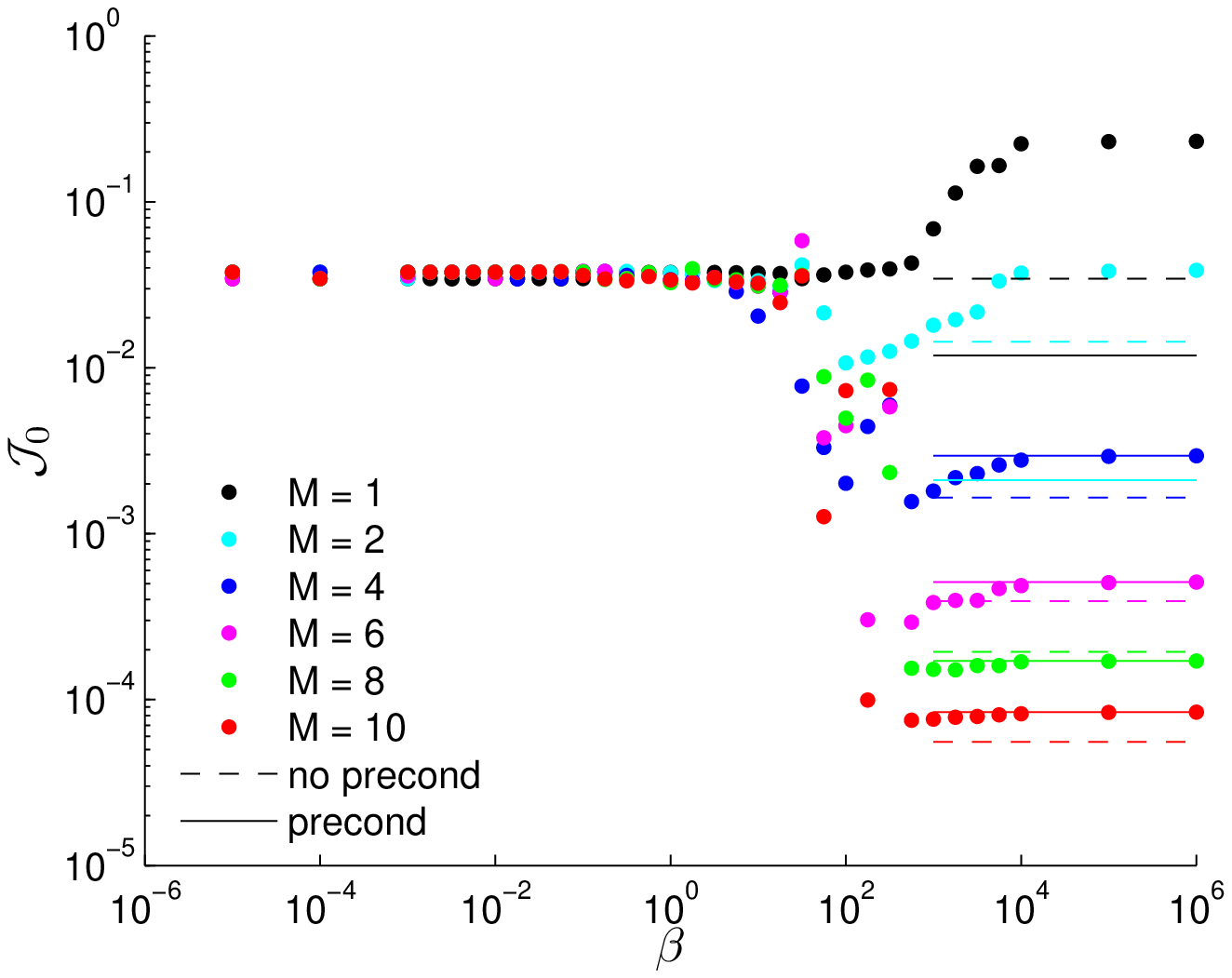}}
  \subfigure[model \#3]{\includegraphics[width=0.5\textwidth]{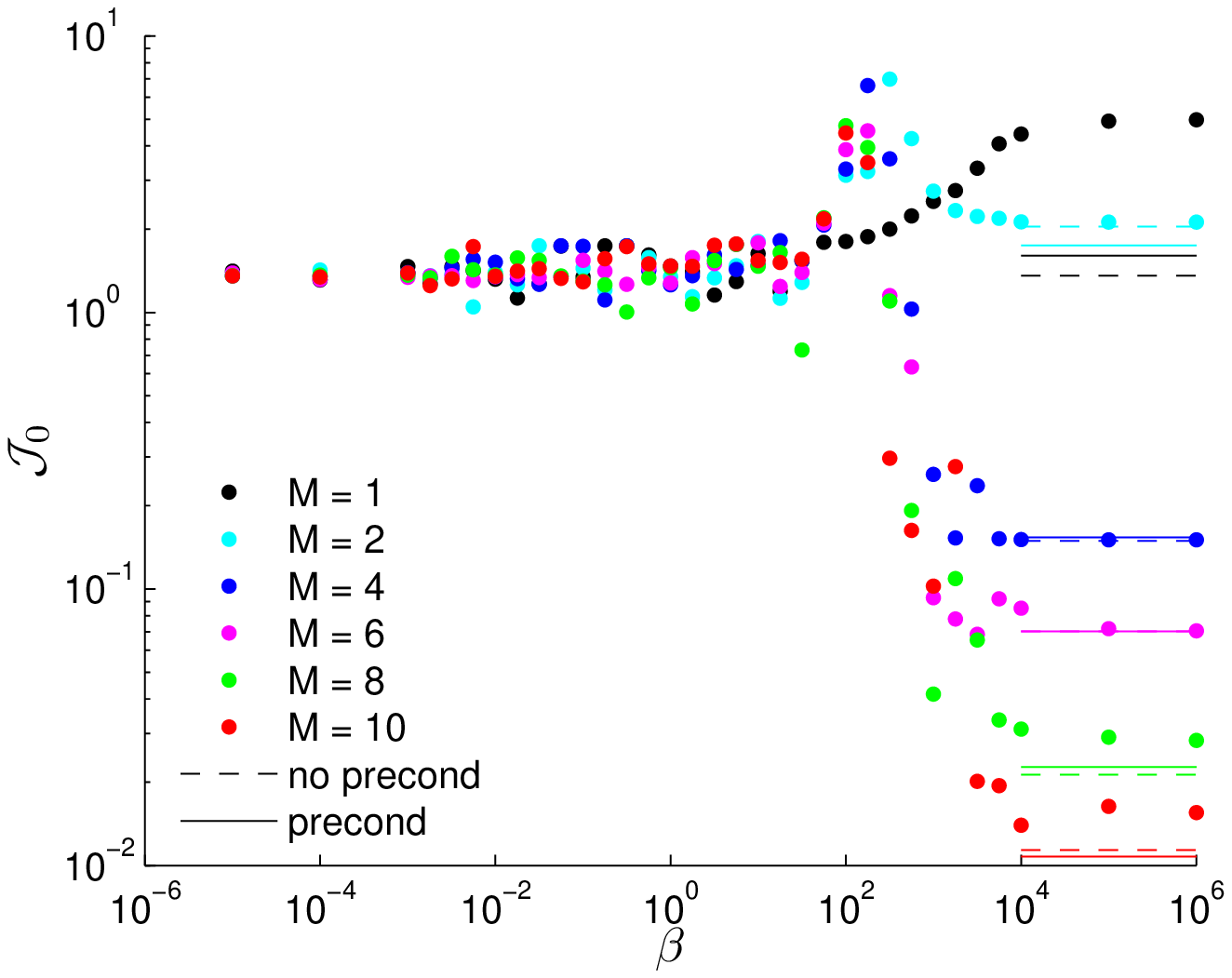}}}
  \end{center}
  \caption{Convergence analysis for $s(t)$ performed for (a,c) model~\#2 and (b,d) model~\#3 by evaluating
    (a,b) solution norms $\| \bs - s_{\rm{true}} \|_{L_2}$ and (c,d) data mismatch parts $\mJ_0$ of cost functional $\mJ$.
    Colors represent the number of available points to construct piecewise approximation to $s(t)$:
    (black)~$M=1$, (cyan)~$M=2$, (blue)~$M=4$, (pink)~$M=6$, (green)~$M=8$ and (red)~$M=10$.
    Dashed and solid lines show the values obtained for case~\#1 respectively with no preconditioning
    and using optimal preconditioning. Dots are used for case~\#2 to represent solution norms and data mismatch
    values as functions of regularization coefficient $\beta$ in \eqref{Eq:W:1:8reg}.}
  \label{fig:reg_s}
\end{figure}

Figure~\ref{fig:reg_s} shows the results of reconstructing free boundary $s(t)$ without noise in measurements for
models~\#2 and \#3 for different values of $M = \{1, 2, 4, 6, 8, 10\}$ using respective colors:
black, cyan, blue, pink, green and red. Lines represent
the values obtained for case \#1 using piecewise initial guess $s_{{\rm ini},M} (t)$ and with no preconditioning
(shown by dashed lines) and using optimal preconditioning (shown by solid lines) as discussed in Section~\ref{sec:single_rec}.
Dots represent the solution norms $\| \bs - s_{\rm{true}} \|_{L_2}$ and data mismatch parts $\mJ_0$ of cost functional $\mJ$
recorded after performing optimization each time supplied with different value of regularization coefficient $\beta$ in
\eqref{Eq:W:1:8reg}. The observed results allow us to make the following comments. First, positive effect of
preconditioning is seen for small M, e.g.~$M = 1$ and $M = 2$, while for $M \geq 4$ gradients with no preconditioning
have the same, or even better, performance than preconditioned ones. The former relates to the general effect of smoothing
gradients discussed previously in Section~\ref{sec:single_rec}. The latter could be explained by better ability of non-modified
(by smoothing) gradients to find ``better'' local optimizer if the initial guess is close to the true solution. Second,
the performance for case \#1 with no preconditioning and case \#2 with added regularization is comparable when regularization
weighting coefficient $\beta$ is sufficiently large. This also helps to identify model dependent thresholds for $\beta$ to
``calibrate'' regularization procedure used in case~\#2. For the rest of numerical experiments shown in this section we use
$\beta = 10^3$ for model~\#2 and $\beta = 10^4$ for model~\#3.

\begin{figure}[htb!]
  \begin{center}
  \mbox{
  \subfigure[model \#2: case \#1, no-precond vs.~precond]{\includegraphics[width=0.5\textwidth]{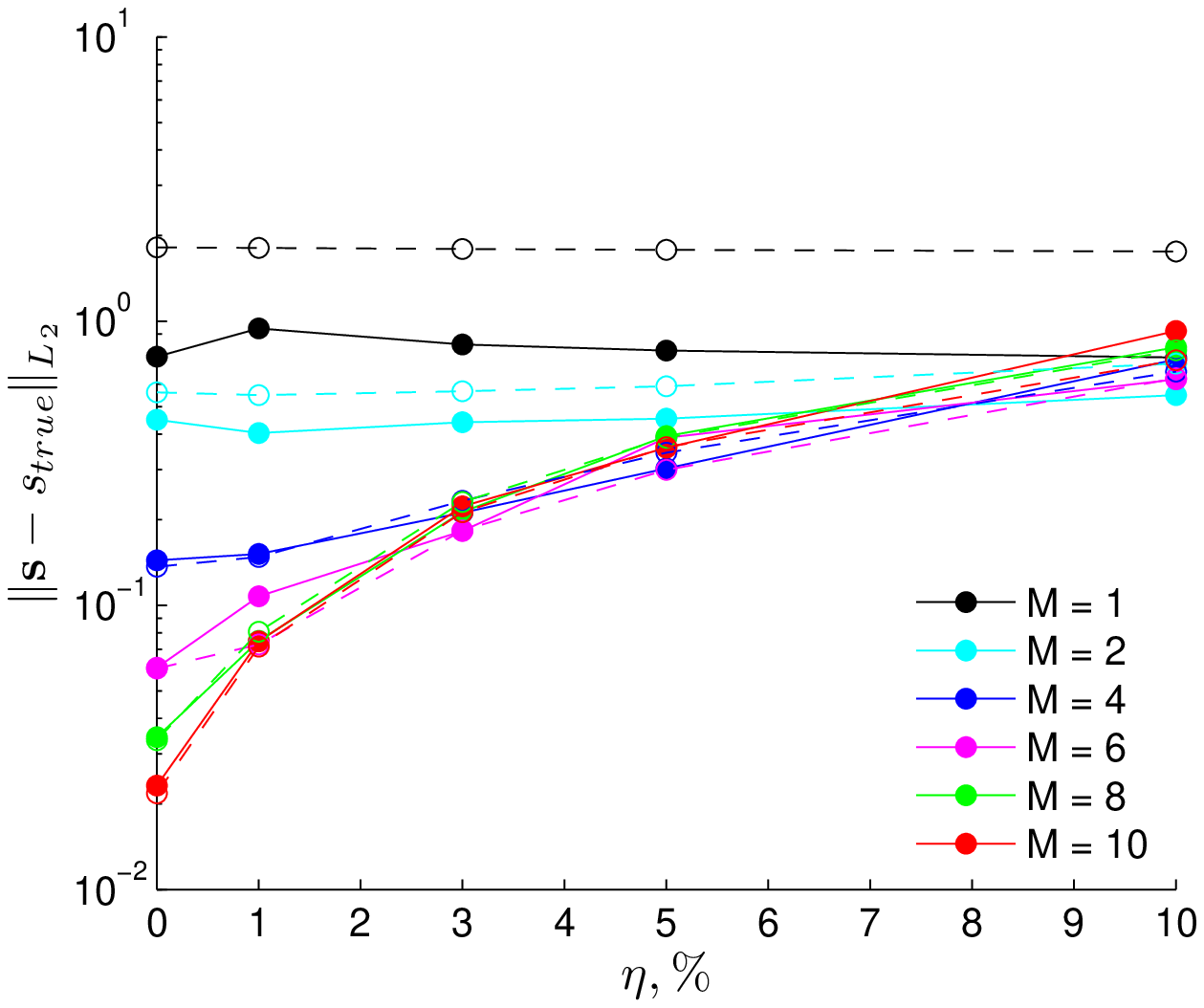}}
  \subfigure[model \#2: case \#1 (no-precond) vs.~case \#2]{\includegraphics[width=0.5\textwidth]{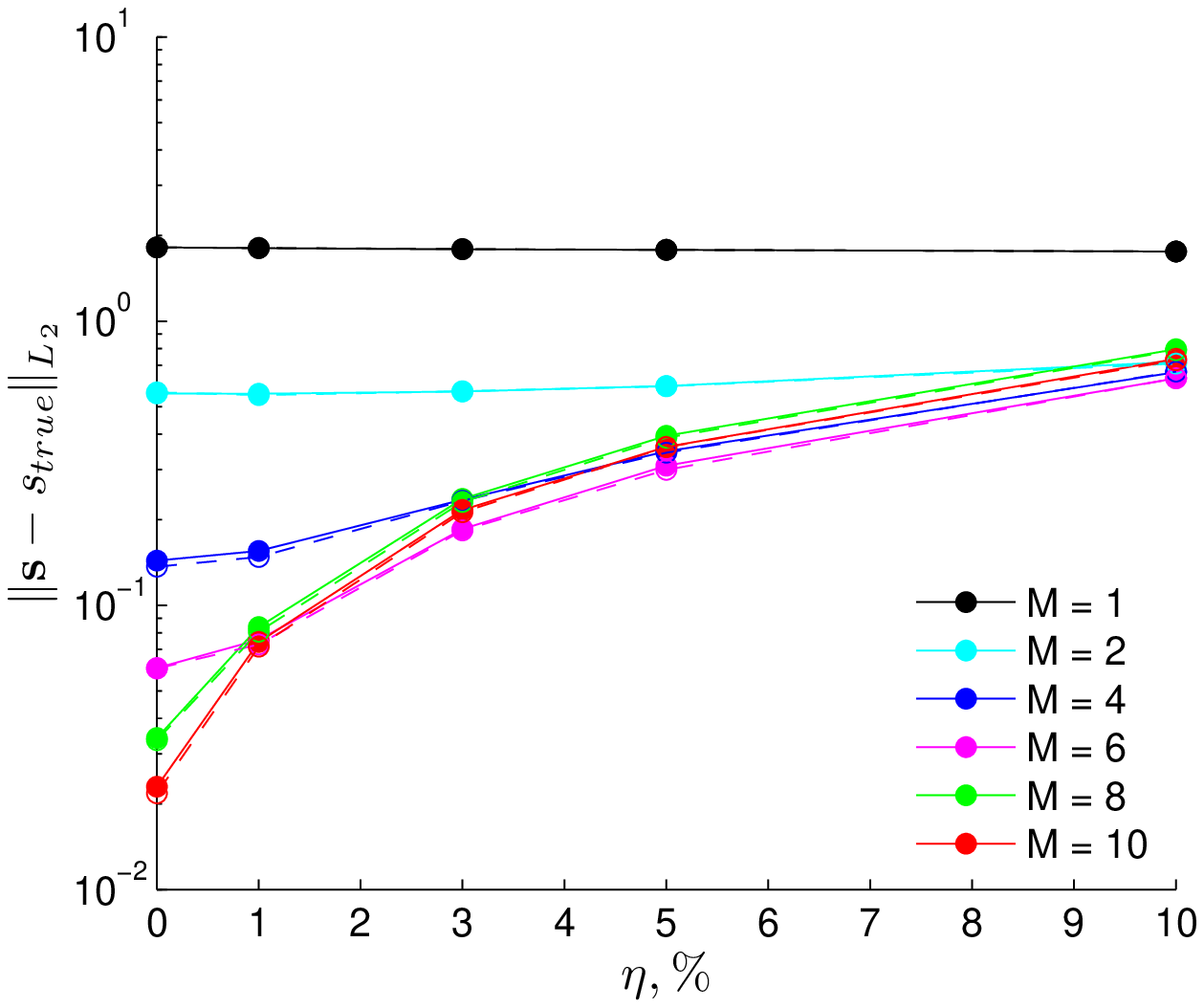}}}
  \mbox{
  \subfigure[model \#2: case \#1, no-precond vs.~precond]{\includegraphics[width=0.5\textwidth]{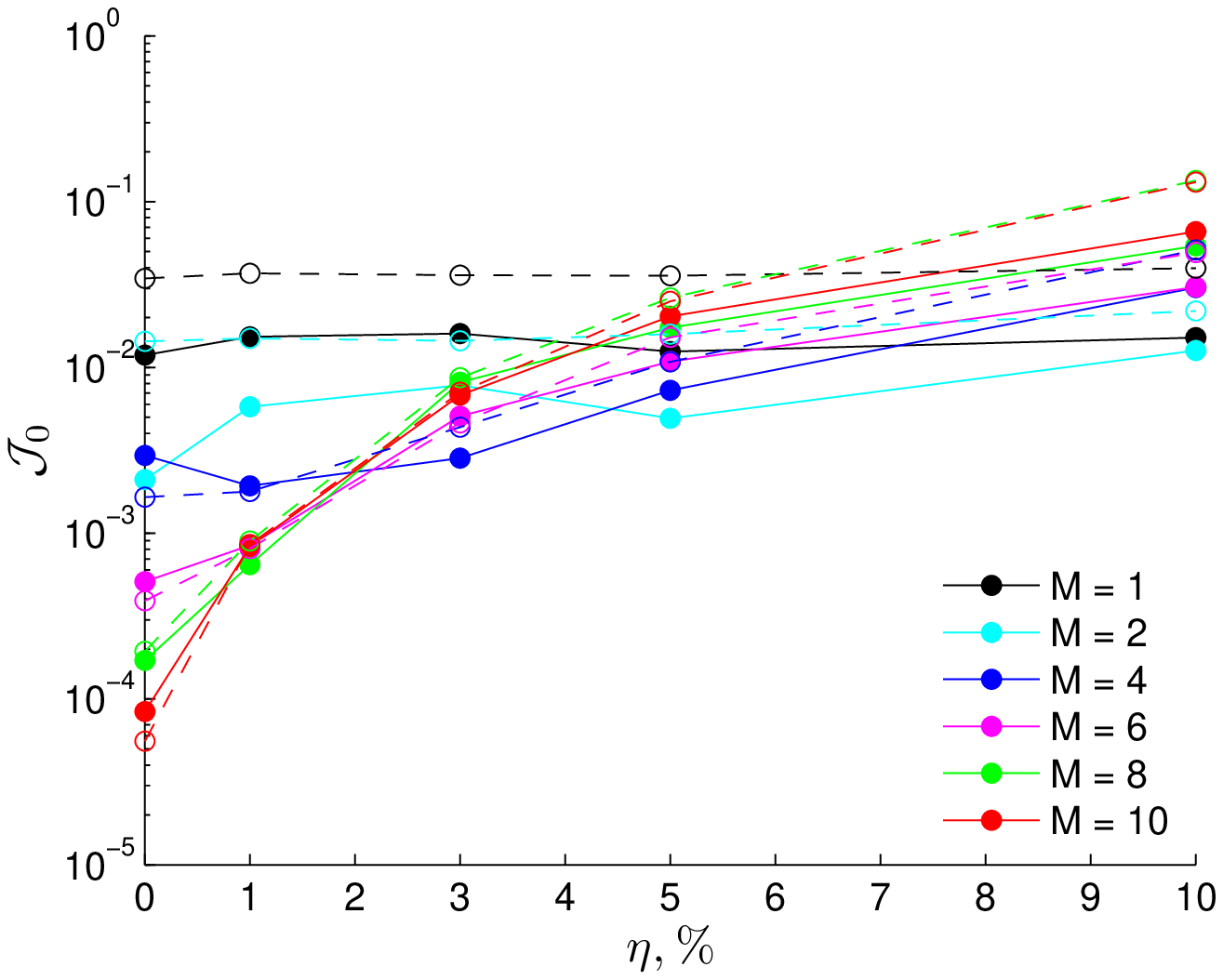}}
  \subfigure[model \#2: case \#1 (no-precond) vs.~case \#2]{\includegraphics[width=0.5\textwidth]{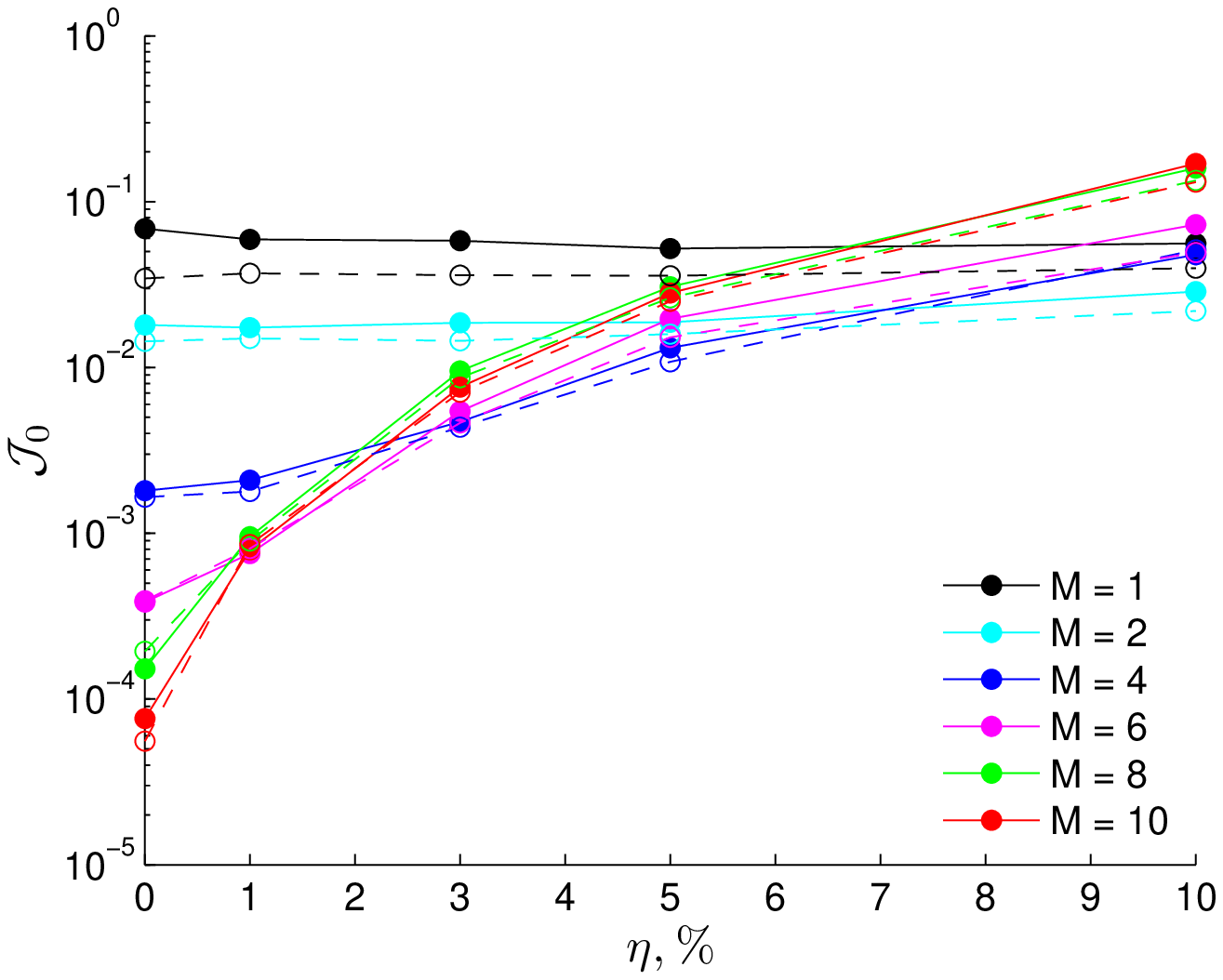}}}
  \end{center}
  \caption{Convergence analysis for $s(t)$ performed for model \#2 by evaluating (a,b) solution norms $\| \bs - s_{\rm{true}} \|_{L_2}$
    and (c,d) data mismatch parts $\mJ_0$ of cost functional $\mJ$ obtained in the presence of noise $\eta$ and averaged over 10
    noise samples. Colors represent the number of available points to construct piecewise approximation to $s(t)$:
    (black)~$M=1$, (cyan)~$M=2$, (blue)~$M=4$, (pink)~$M=6$, (green)~$M=8$ and (red)~$M=10$.
    In (a,c) the results are obtained for case~\#1 (dashed lines and empty circles) without preconditioning and (solid lines and
    filled circles) with preconditioning. In (b,d) the results are obtained for (dashed lines and empty circles) case~\#1 without preconditioning and (solid lines and filled circles) case~\#2 for regularization without preconditioning.}
  \label{fig:noise_s_model_2}
\end{figure}
\begin{figure}[htb!]
  \begin{center}
  \mbox{
  \subfigure[model \#3: case \#1, no-precond vs.~precond]{\includegraphics[width=0.5\textwidth]{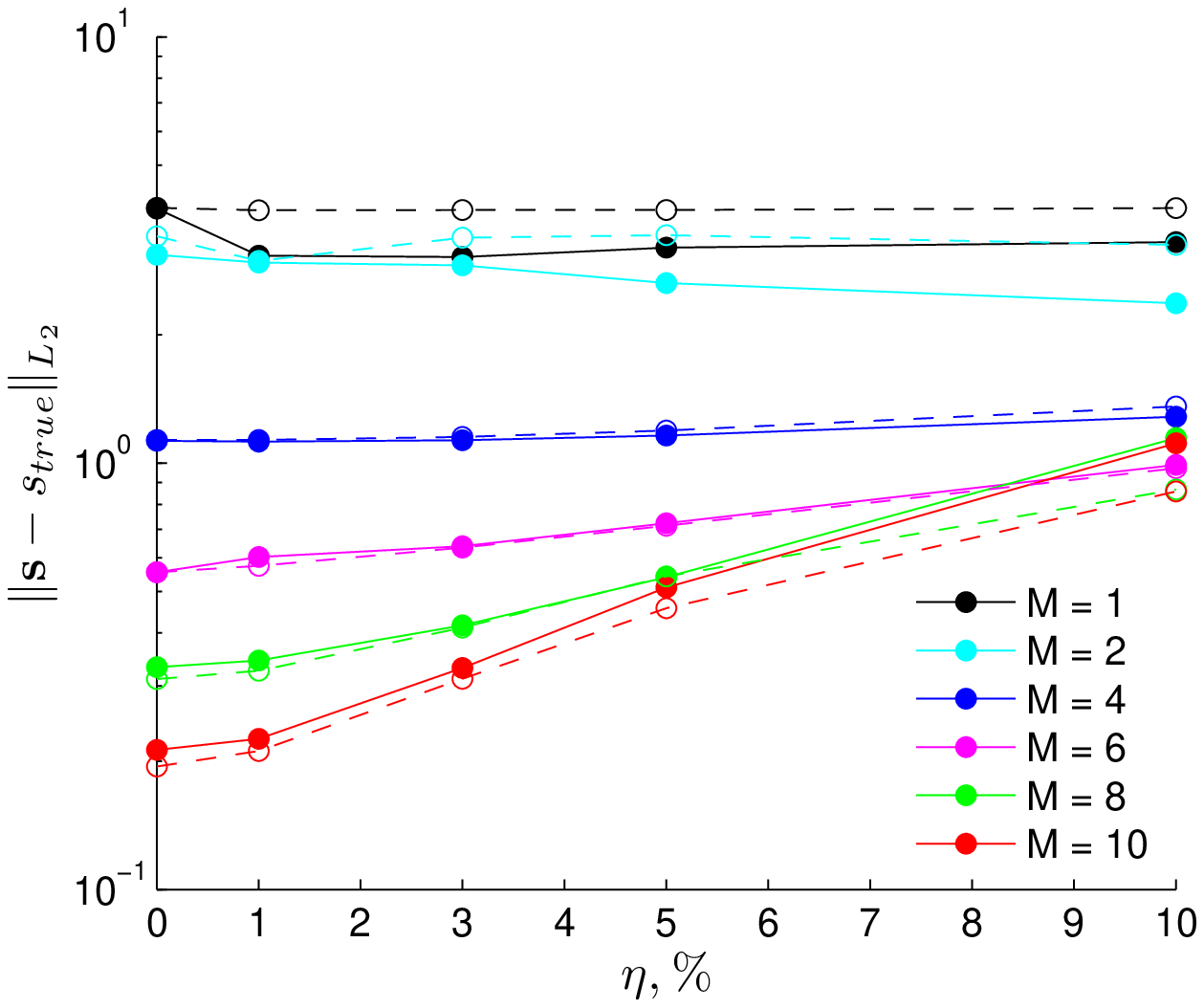}}
  \subfigure[model \#3: case \#1 (no-precond) vs.~case \#2]{\includegraphics[width=0.5\textwidth]{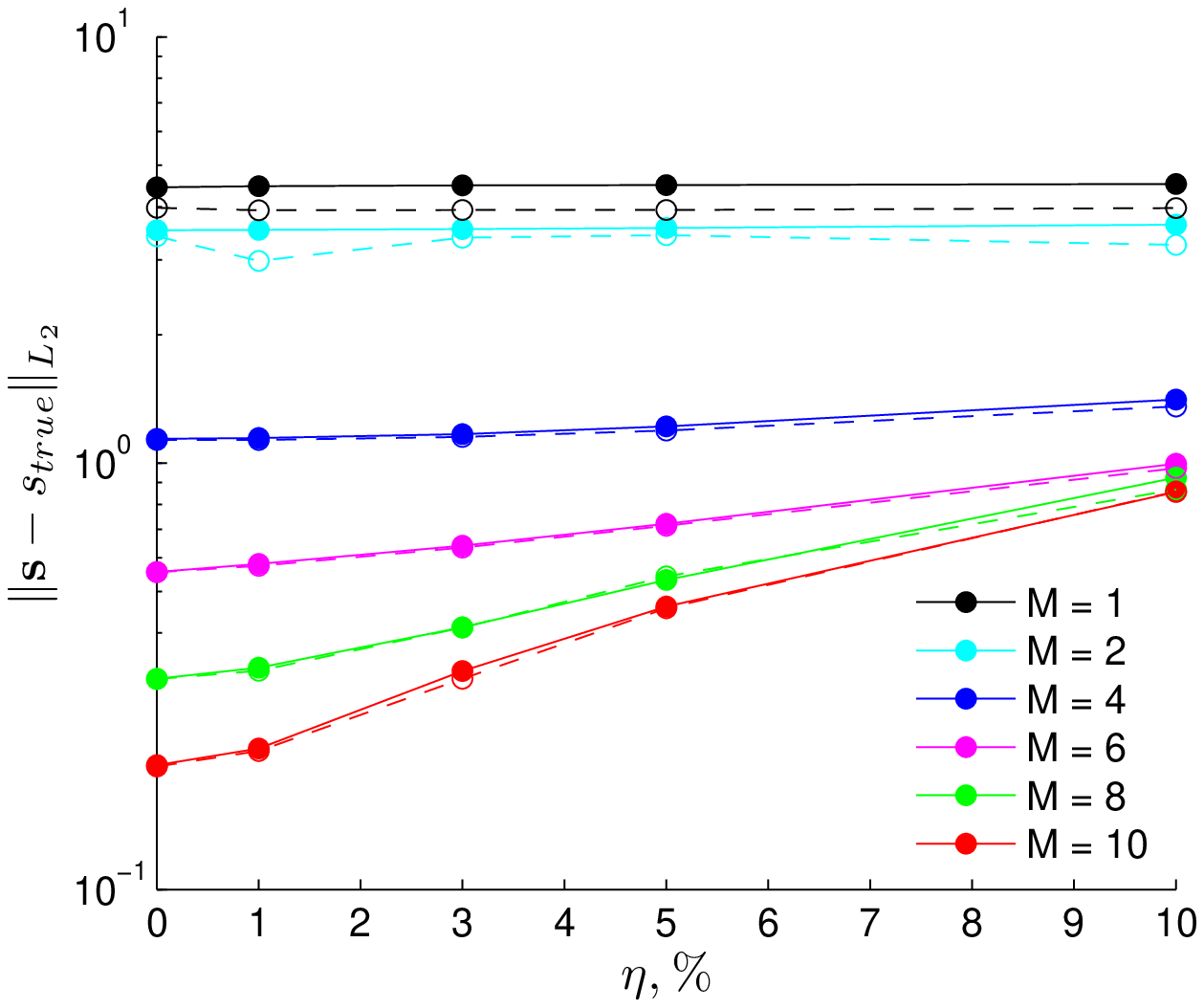}}}
  \mbox{
  \subfigure[model \#3: case \#1, no-precond vs.~precond]{\includegraphics[width=0.5\textwidth]{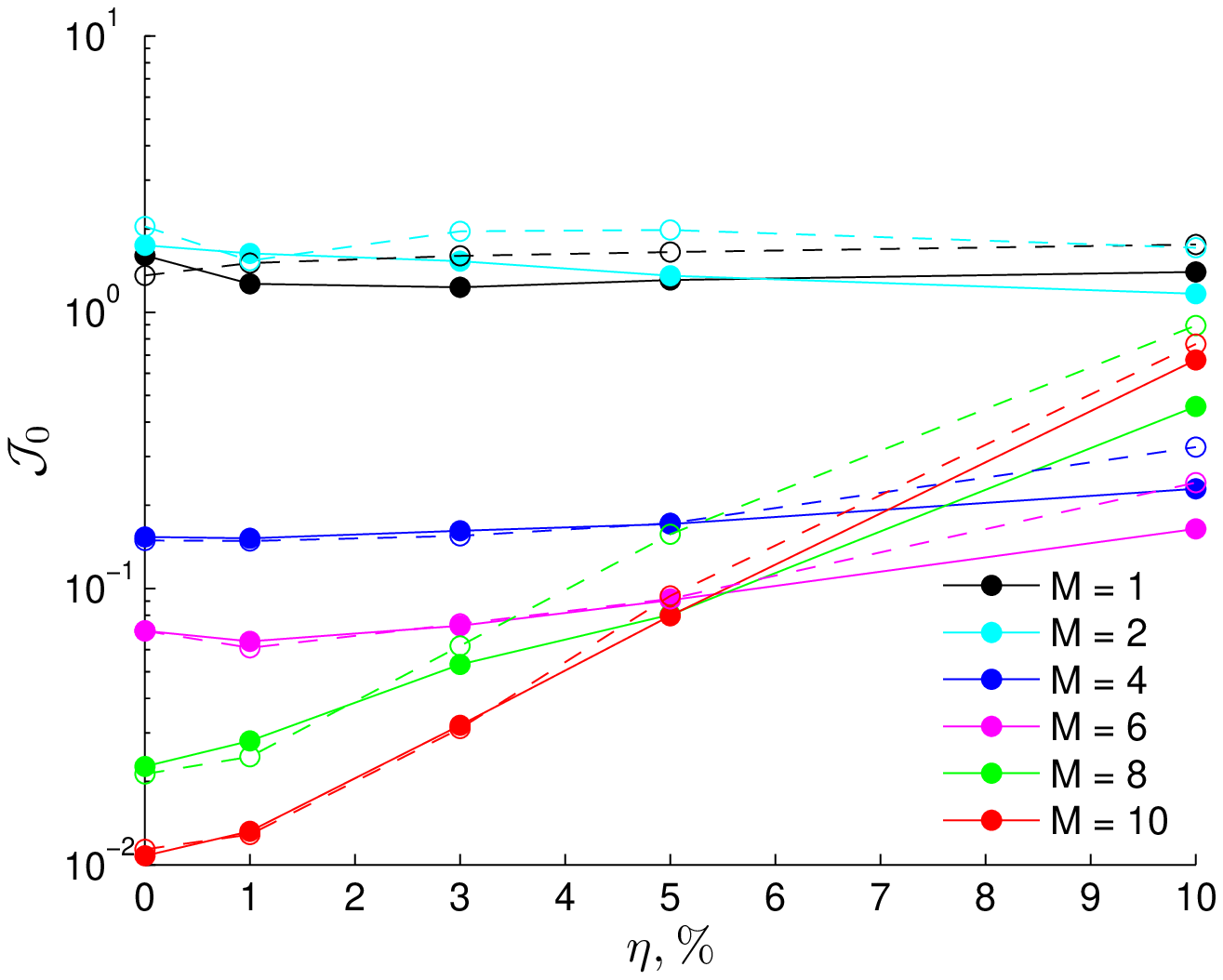}}
  \subfigure[model \#3: case \#1 (no-precond) vs.~case \#2]{\includegraphics[width=0.5\textwidth]{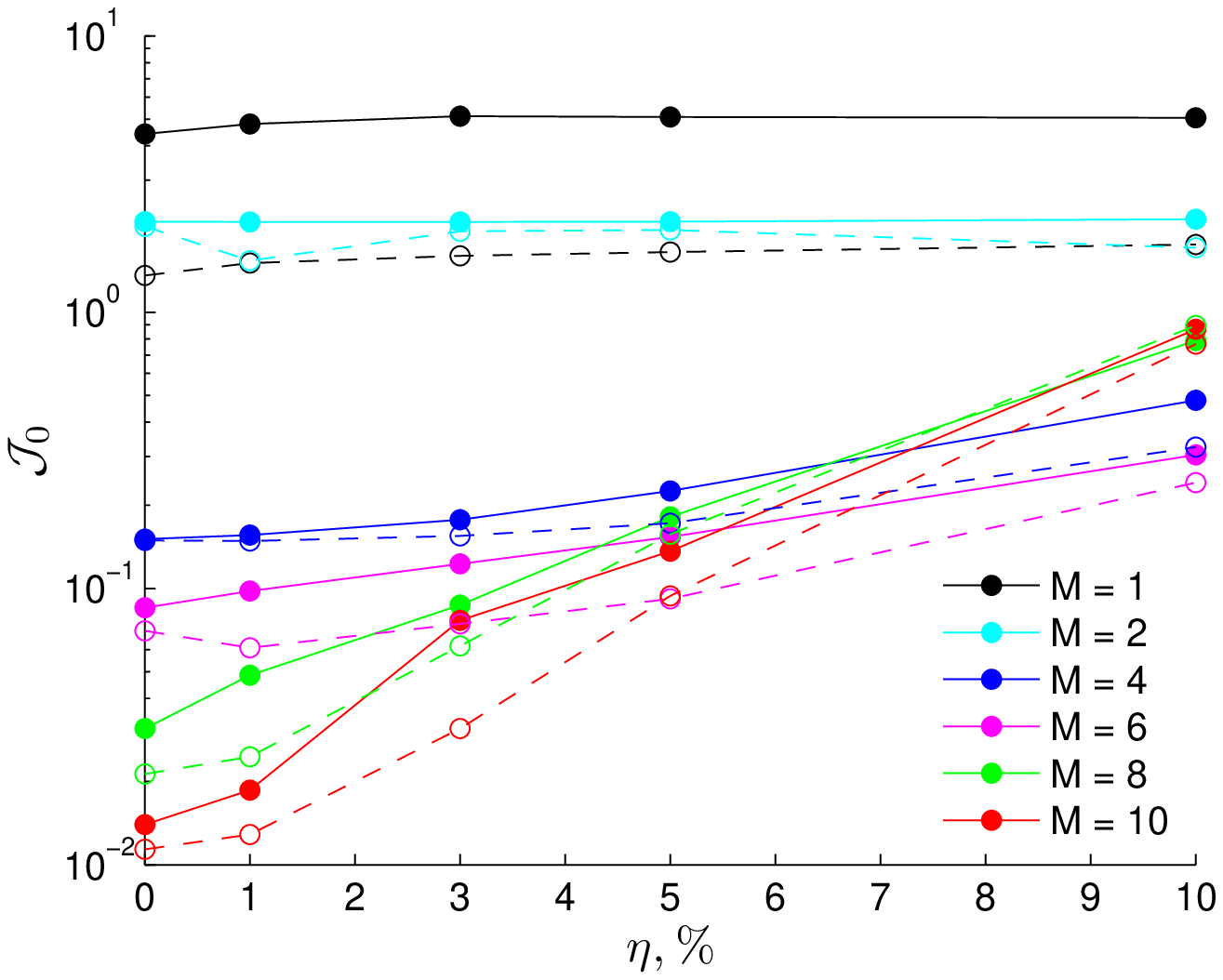}}}
  \end{center}
  \caption{Convergence analysis for $s(t)$ performed for model \#3 by evaluating (a,b) solution norms $\| \bs - s_{\rm{true}} \|_{L_2}$
    and (c,d) data mismatch parts $\mJ_0$ of cost functional $\mJ$ obtained in the presence of noise $\eta$ and averaged over 10
    noise samples. Colors represent the number of available points to construct piecewise approximation to $s(t)$:
    (black)~$M=1$, (cyan)~$M=2$, (blue)~$M=4$, (pink)~$M=6$, (green)~$M=8$ and (red)~$M=10$.
    In (a,c) the results are obtained for case~\#1 (dashed lines and empty circles) without preconditioning and (solid lines and
    filled circles) with preconditioning. In (b,d) the results are obtained for (dashed lines and empty circles) case~\#1 without preconditioning and (solid lines and filled circles) case~\#2 for regularization without preconditioning.}
  \label{fig:noise_s_model_3}
\end{figure}

In Figures~\ref{fig:noise_s_model_2} and \ref{fig:noise_s_model_3} we present the results of reconstructing
free boundary $s(t)$ respectively for models~\#2 and \#3 for different values of $M$.
These results are obtained using the approaches described earlier in the current section, namely case~\#1 without
preconditioning (dashed lines and empty circles in (a-d)), case~\#1 with preconditioning (solid lines and
filled circles in (a,c)), and case~\#2 for regularization without preconditioning (solid lines and filled circles in (b,d)).
To perform optimization we use additional data $\{\tilde s_i^{\eta} \}_{i=1}^M$ contaminated with 1\%, 3\%, 5\%, 10\% normally
distributed noise and then we average the obtained results over 10 different noise samples. For both models we compare performance
of using the preconditioning technique for case~\#1 (Figures~\ref{fig:noise_s_model_2}(a,c) and \ref{fig:noise_s_model_3}(a,c)).
Positive effect of preconditioning is seen again for small M, e.g.~$M = 1$ and $M = 2$. This is consistent with our previous
statement which now could be extended also for data containing sufficiently large, up to 10\%, noise. We also
conclude that the effect of adding regularization introduced by applying case~\#2 is comparable with performance of
no-preconditioned case~\#1 for both models. We close this section by concluding that, as expected, adding additional
measurements for the position of free boundary $s(t)$ has regularizing effect on reconstructing $s(t)$ in the presence
of large amount of noise. Such systematic methodology is also seen very useful to determine the optimal number $M$
of additional measurements in case the noise level $\eta$ is a priori estimated.

\subsection{Identification of the free boundary and other control parameters}
\label{sec:multiple_rec}

As the final step in our validation procedure, here we show the performance of utilizing
Algorithm~\ref{alg:gen_opt} in full, namely for identifying simultaneously several control parameters:
free boundary $s(t)$ and left boundary heat flux $g(t)$. First, we consider reconstructing
$g(t)$ alone for the same three models described in Section~\ref{sec:models} with fixed
boundary $s(t) = s_{\rm{true}}(t)$ to calibrate preconditioning and perform convergence analysis for $g(t)$.
At this point, Algorithm~\ref{alg:gen_opt} is used to find (local) optimal solution $\bg (t)$ iteratively
starting from initial guess $g = g_{\rm ini}(t)$ and setting $s_k(t)$ for every $k = 0, 1, 2, \ldots$
to the true expressions defined analytically in \eqref{eq:model_1}--\eqref{eq:model_3}.
Thus, similarly as used for numerical experiments in Section~\ref{sec:single_rec},
Algorithm~\ref{alg:gen_opt} skips computing the corresponding part of the gradient, namely
$\bnabla_s \mJ$ in \eqref{eq:grads}, and sets $\alpha_k^s$ in \eqref{eq:opt_SD} to zero.
Second, we present the results of identification of full control vector $v = (s(t), g(t))$
with the discussion on the effect of preconditioning and approaches to solve \eqref{eq:opt_alpha}
to find stepsize parameters $\alpha_k^s$ and $\alpha_k^g$ in \eqref{eq:opt_SD}.

\begin{figure}[htb!]
  \begin{center}
  \mbox{
  \subfigure[model \#1]{\includegraphics[width=0.33\textwidth]{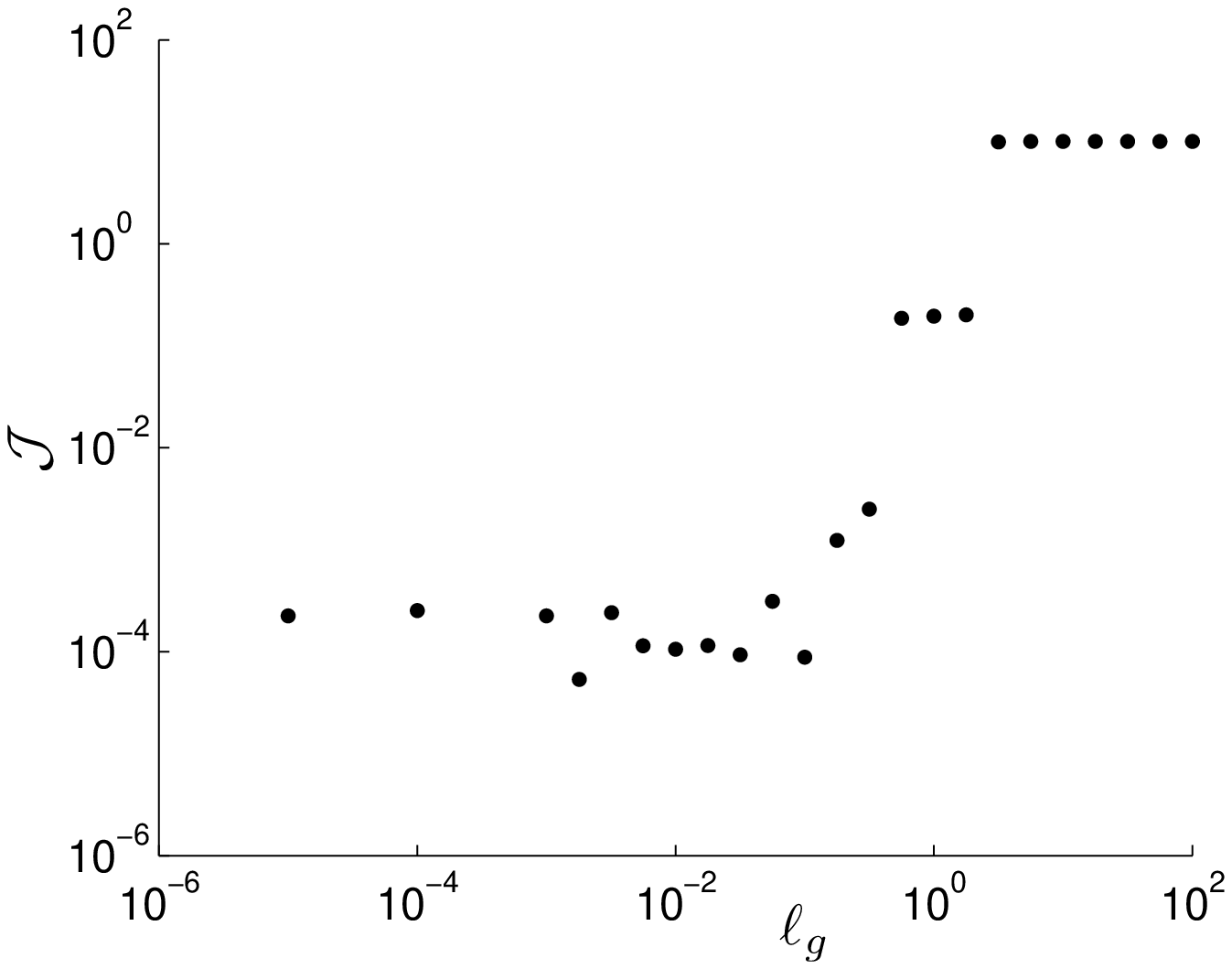}}
  \subfigure[model \#2]{\includegraphics[width=0.33\textwidth]{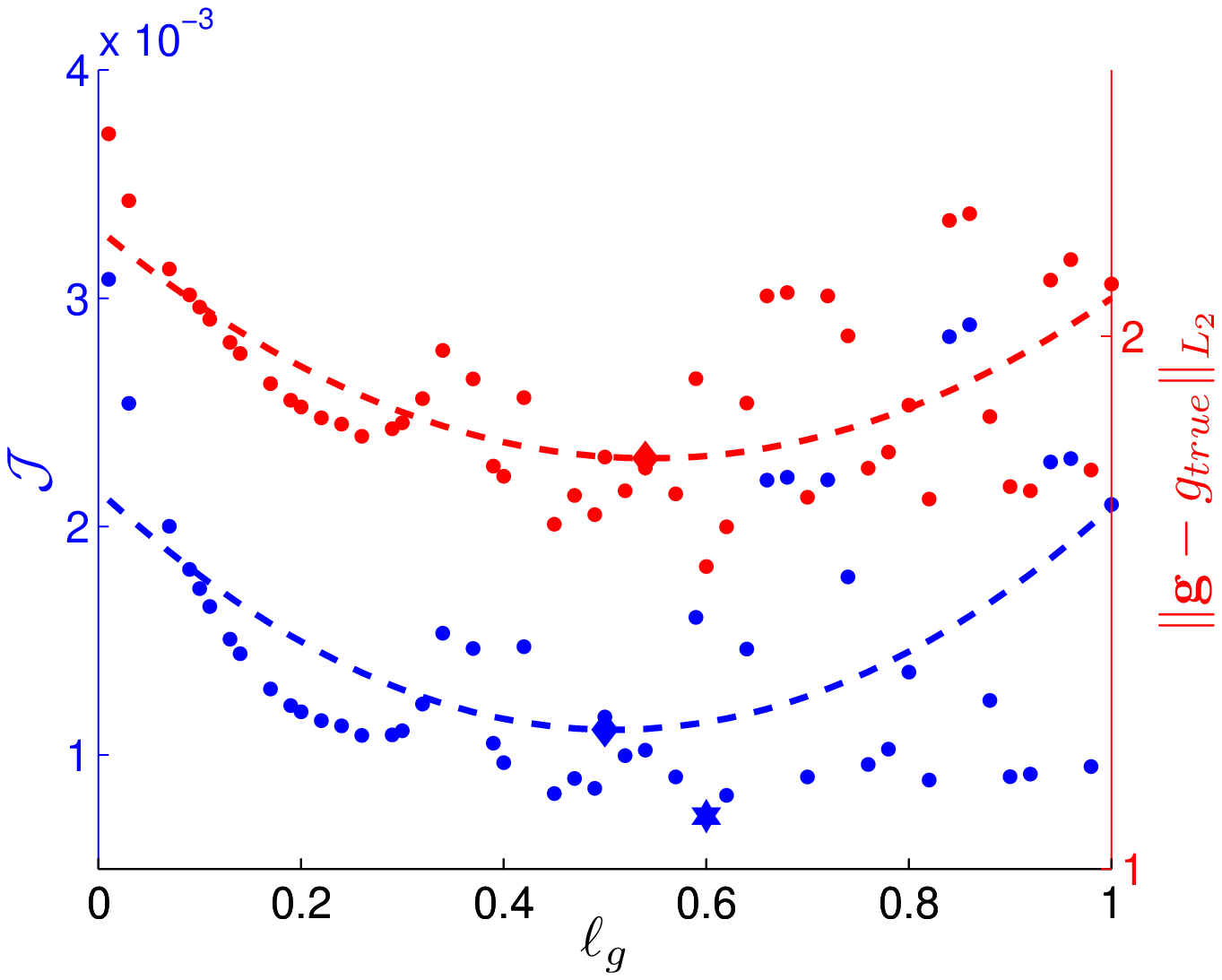}}
  \subfigure[model \#3]{\includegraphics[width=0.33\textwidth]{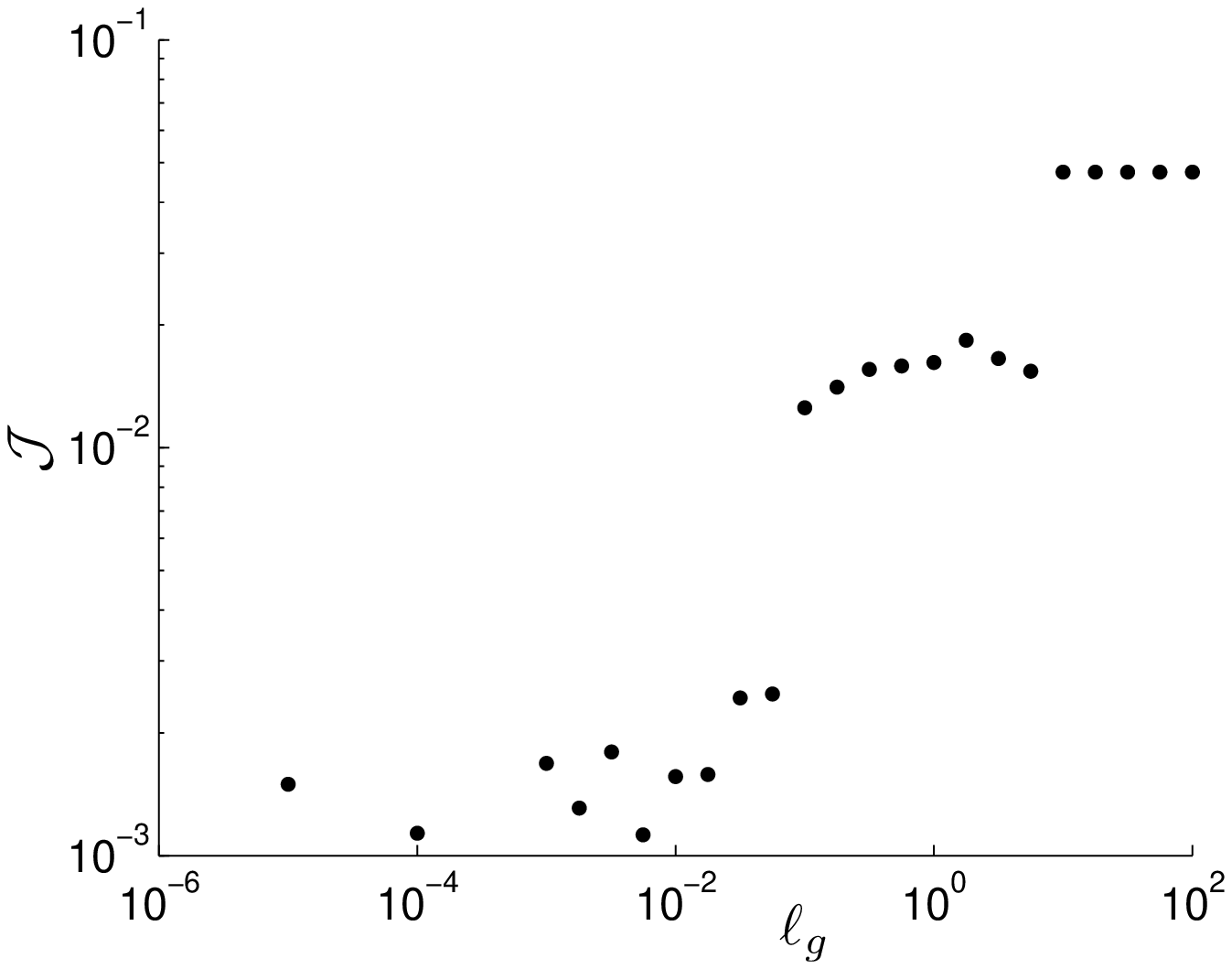}}}
  \end{center}
  \caption{(a,c)~Cost functional values $\mJ$ (black dots) computed for different orders of $\ell_g$
    for models \#1 and \#3. (b)~Approximation of optimal parameter $\ell^*_g$ for preconditioning
    procedure for model \#2. Cost functional values $\mJ$ and solution norms $\| \bg - g_{\rm{true}} \|_{L_2}$
    are represented respectively by blue and red dots. Quadratic regression models for cost functionals
    and solution norms are represented by dashed lines with minimal values shown by diamonds. The best (minimal) value
    of $\mJ$ in the proximity of approximated $\ell^*_g$ is shown by blue hexagon.}
  \label{fig:precond_g}
\end{figure}

As described previously in Section~\ref{sec:single_rec}, we intend to repeat the calibration of the preconditioning
procedure to determine the sensitivity of optimal solution $\bg (t)$ to the choice of preconditioning parameter $\ell_g$
in \eqref{eq:helm}. As before, we have used the same two criteria, namely evaluation of cost functional values $\mJ$ and
solution norms $\| \bg - g_{\rm{true}} \|_{L_2}$. As shown in Figure~\ref{fig:precond_g}(b) cost functional values
(blue dots) and solution norms (red dots) are recorded after performing optimizations supplied with different values
of $\ell_g$ for model~\#2. Both sets of points are then used to perform the least square analysis to find the
quadratic regression model (dashed lines) for each set. The quadratic function to model $\mJ$ is
then minimized to approximate $\ell^*_g$ (blue diamond) giving value $\ell^*_{g,2} = 0.5$.
The quadratic function to model the solution norm is also minimized to confirm the proximity
of the obtained solution (red diamond) to approximated $\ell^*_g$ and demonstrate consistence of the
obtained results.

In fact, such approach cannot work for models \#1 and \#3, as we believe, due to respectively their simplicity and
complexity. As shown in Figure~\ref{fig:precond_g}(a,c), based on cost functional values (black dots) computed
for different orders of $\ell_g$ we do not see noticeable improvement in applying preconditioning, and hence we are
not able to identify the interval where approximation via quadratic regression model could suggest any solution for $\ell^*_g$.
Anyway, when applying preconditioning we use $\ell^*_{g,1} = \ell^*_{g,3} = 10^{-2}$ for both models.
Unless stated otherwise, computational results discussed further in this section use mentioned above values
$\ell^*_{g,i}, \, i = 1, 2, 3$, whenever preconditioning procedure is active for $g(t)$.

Next, similarly to $s(t)$, we validate our computational approach in terms of convergence of $g(t)$ to the global
optimal solution by solving the same optimization problem for all three models starting with different initial
guesses. Again, these new initial guesses $g_{\rm{ini},\lambda_g}$ are parameterized with respect to their proximity
to global minimizer $g_{\rm{true}}$ in the following way
\begin{equation}
  g_{\rm{ini},\lambda_g}(t) = (1 - \lambda_g) g_{\rm{true}}(t) + \lambda_g g_{\rm{ini}}(t).
  \label{eq:g_convex_comb}
\end{equation}
We note that setting parameter $\lambda_g = 1$ recovers the regular initial guess shown in Figure~\ref{fig:models}(d,e,f),
while $\lambda_g \rightarrow 0$ moves initial guess in the close neighborhood of $g_{\rm{true}}$.

\begin{figure}[htb!]
  \begin{center}
  \mbox{
  \subfigure[model \#1]{\includegraphics[width=0.33\textwidth]{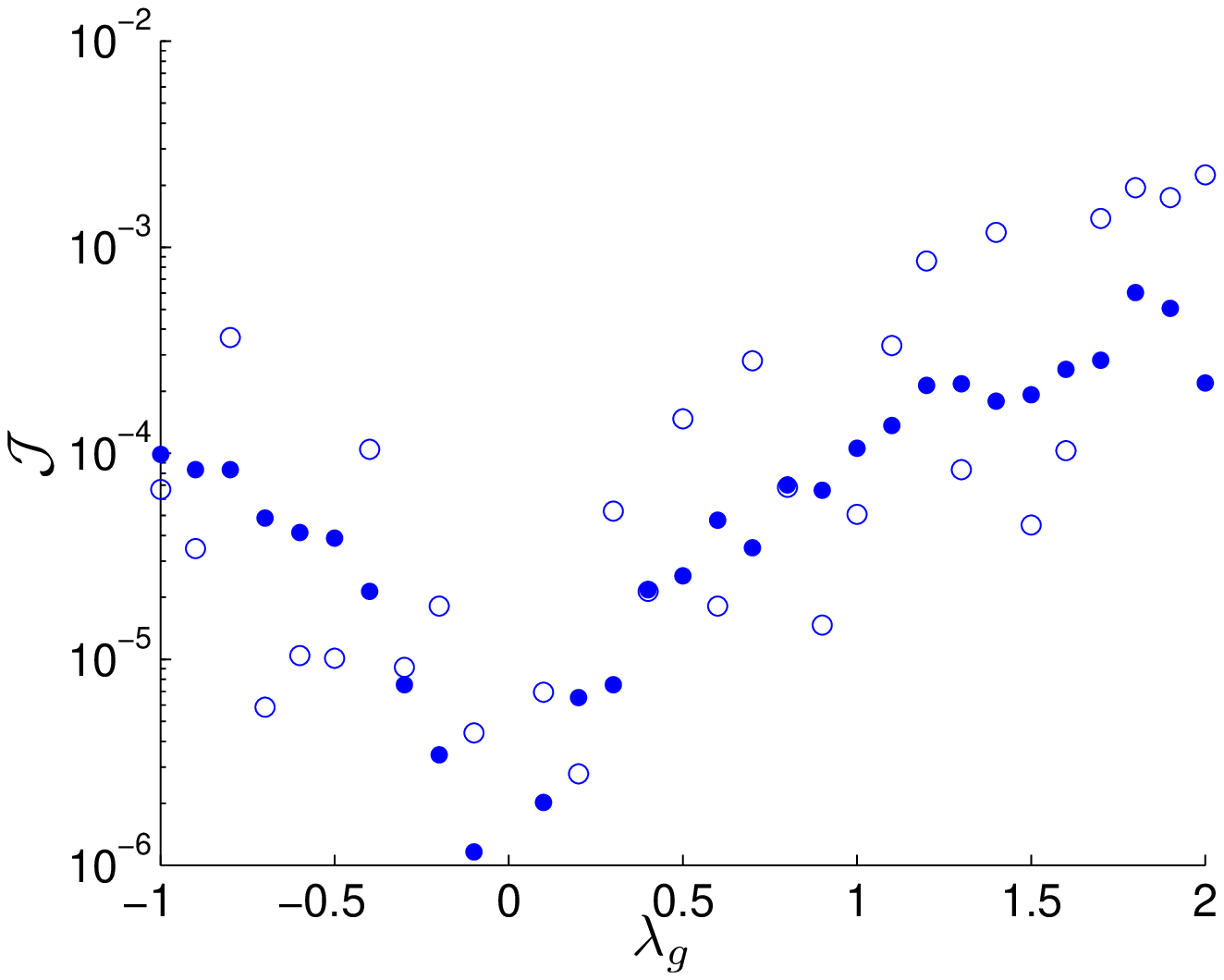}}
  \subfigure[model \#2]{\includegraphics[width=0.33\textwidth]{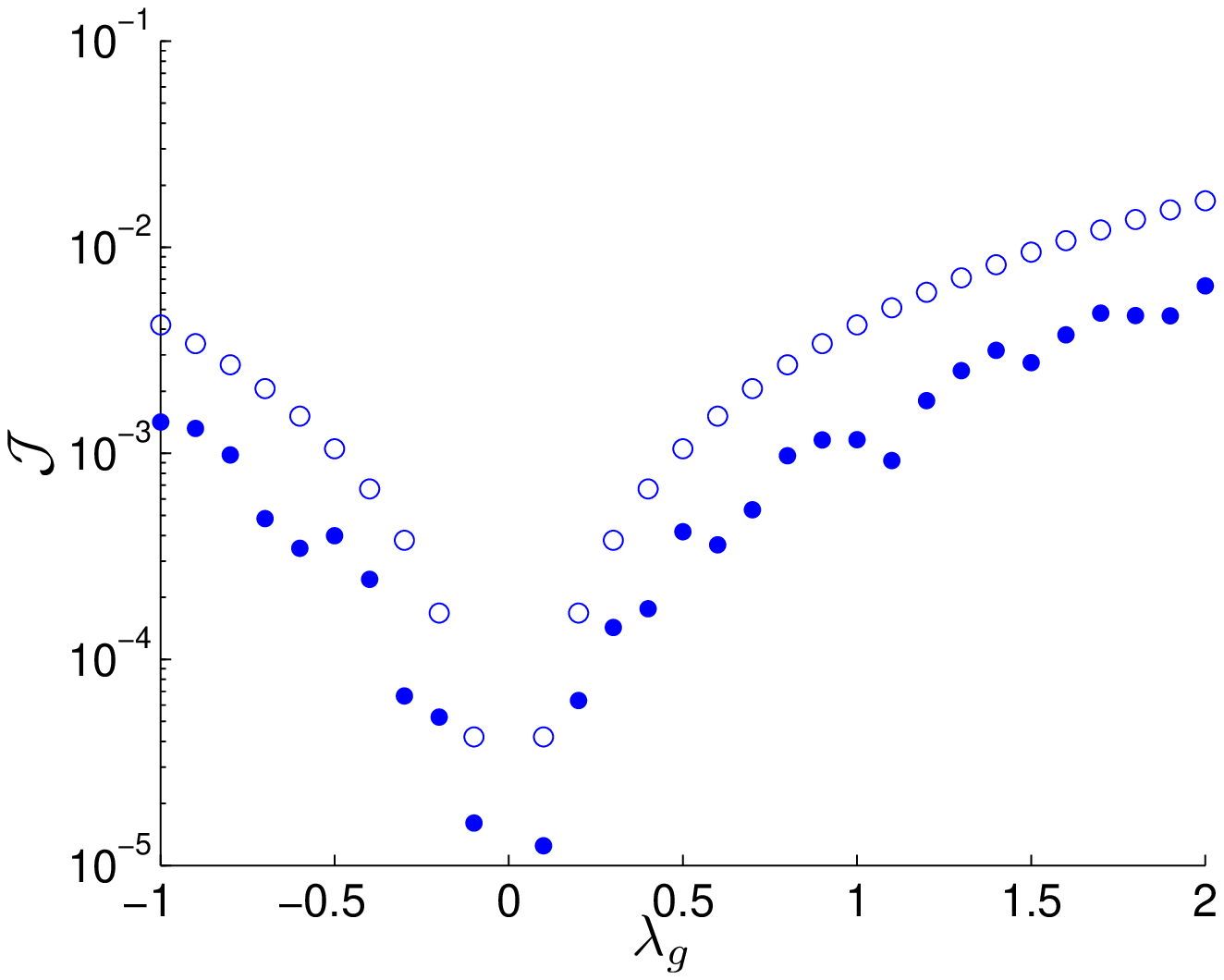}}
  \subfigure[model \#3]{\includegraphics[width=0.33\textwidth]{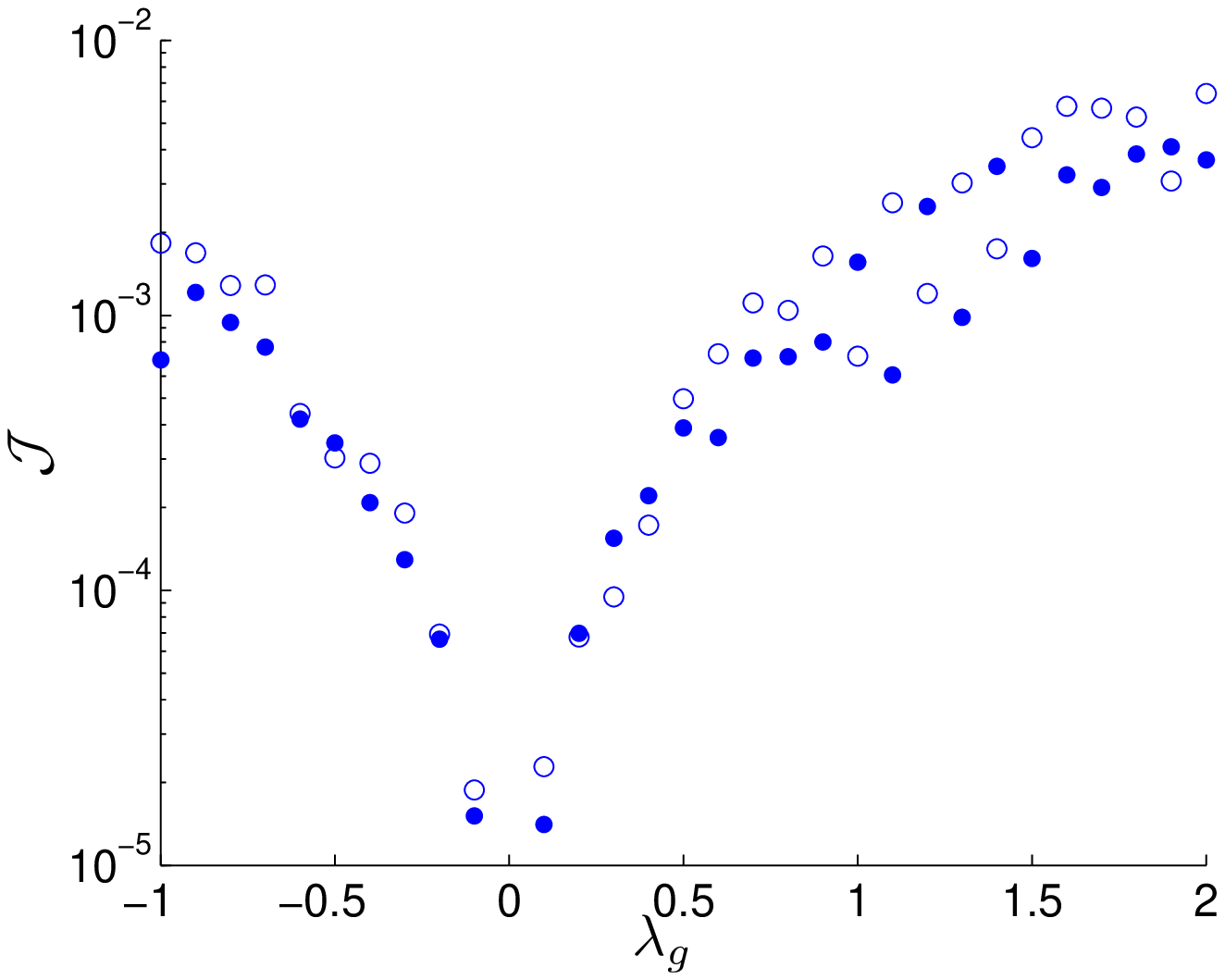}}}
  \mbox{
  \subfigure[model \#1]{\includegraphics[width=0.33\textwidth]{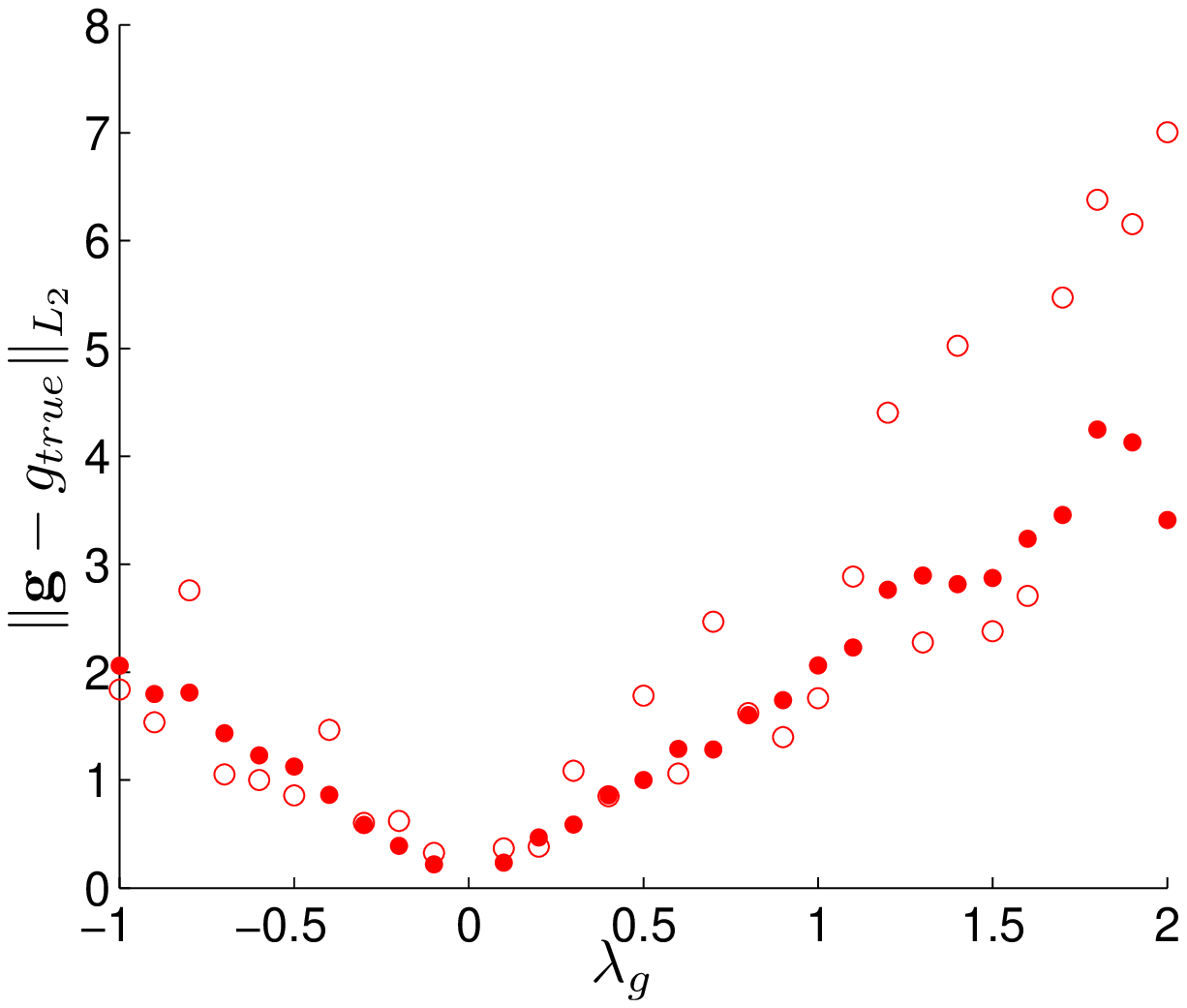}}
  \subfigure[model \#2]{\includegraphics[width=0.33\textwidth]{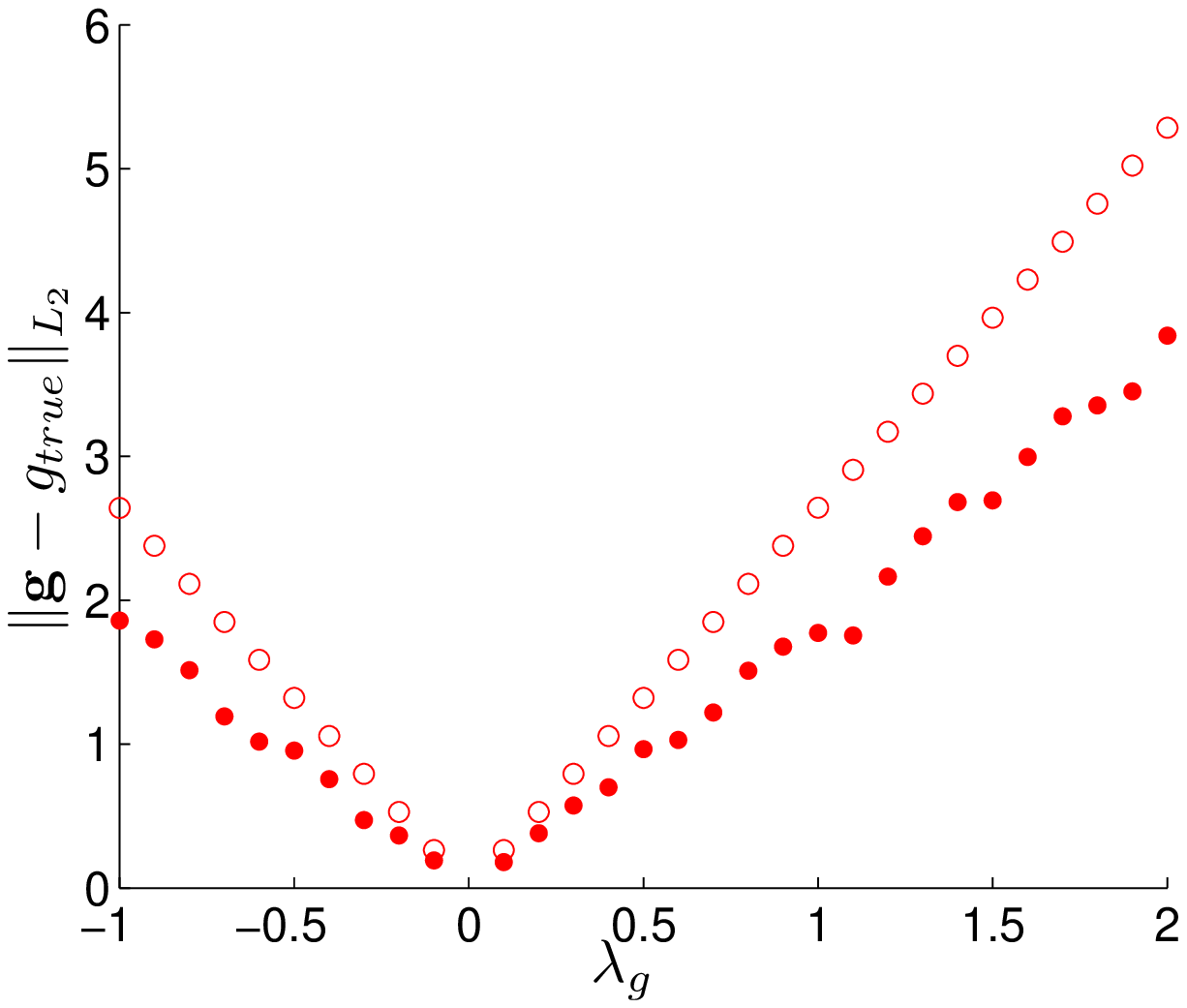}}
  \subfigure[model \#3]{\includegraphics[width=0.33\textwidth]{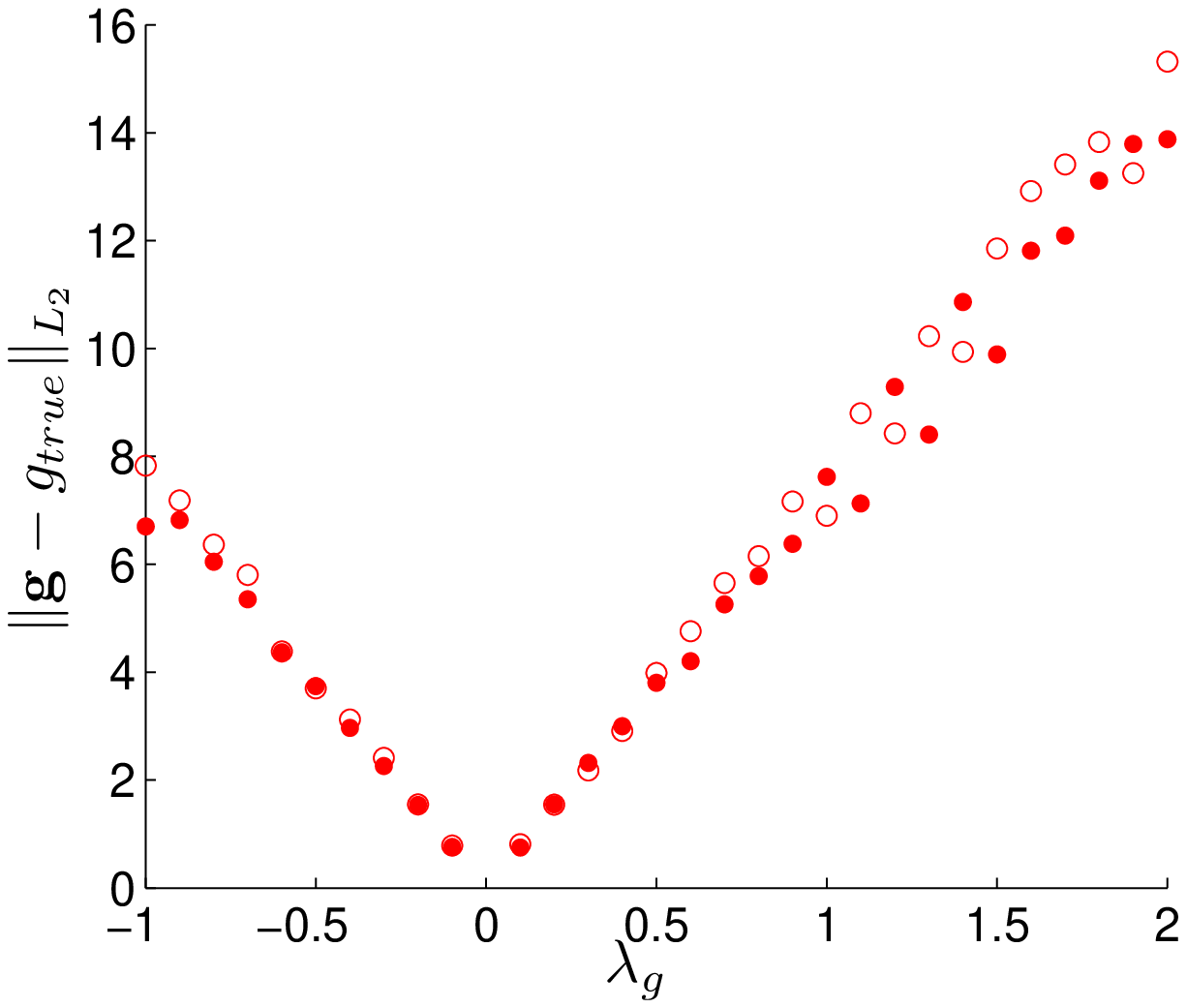}}}
  \end{center}
  \caption{Convergence analysis for $g(t)$ performed for (a,d)~model~\#1, (b,e)~model~\#2, and
    (c,f)~model~\#3 by evaluating (a,b,c) cost functional values $\mJ$ and (d,e,f)
    solution norms $\| \bg - g_{\rm{true}} \|_{L_2}$. Dots and circles represent results
    obtained correspondingly with and without preconditioning procedure \eqref{eq:helm}.}
  \label{fig:conv_g}
\end{figure}

The results of the convergence test for $\lambda_g \in [-1, 2]$ for all three models are shown in Figure~\ref{fig:conv_g}
for cases with and without preconditioning procedure \eqref{eq:helm} by evaluating both cost functional
values $\mJ$ and solution norms $\| \bg - g_{\rm{true}} \|_{L_2}$. The results for all three models show good convergence
to global minimizer $g_{\rm{true}}$, i.e. $\mJ \rightarrow \mJ_{\rm{min}}$ and $\| \bg - g_{\rm{true}} \|_{L_2} \rightarrow 0$
as $\lambda_g \rightarrow 0$. We could also conclude that applying preconditioning benefits in improving convergence significantly
for model~\#2 for which approximated optimal parameter $\ell^*_g$ is found. We should also mention that setting $\ell^*_g$ to a
rather small value makes a stabilizing impact on this convergence as clearly seen, e.g., in model~\#1.

\begin{figure}[htb!]
  \begin{center}
  \mbox{
  \subfigure[model \#1: $g(t)$]{\includegraphics[width=0.33\textwidth]{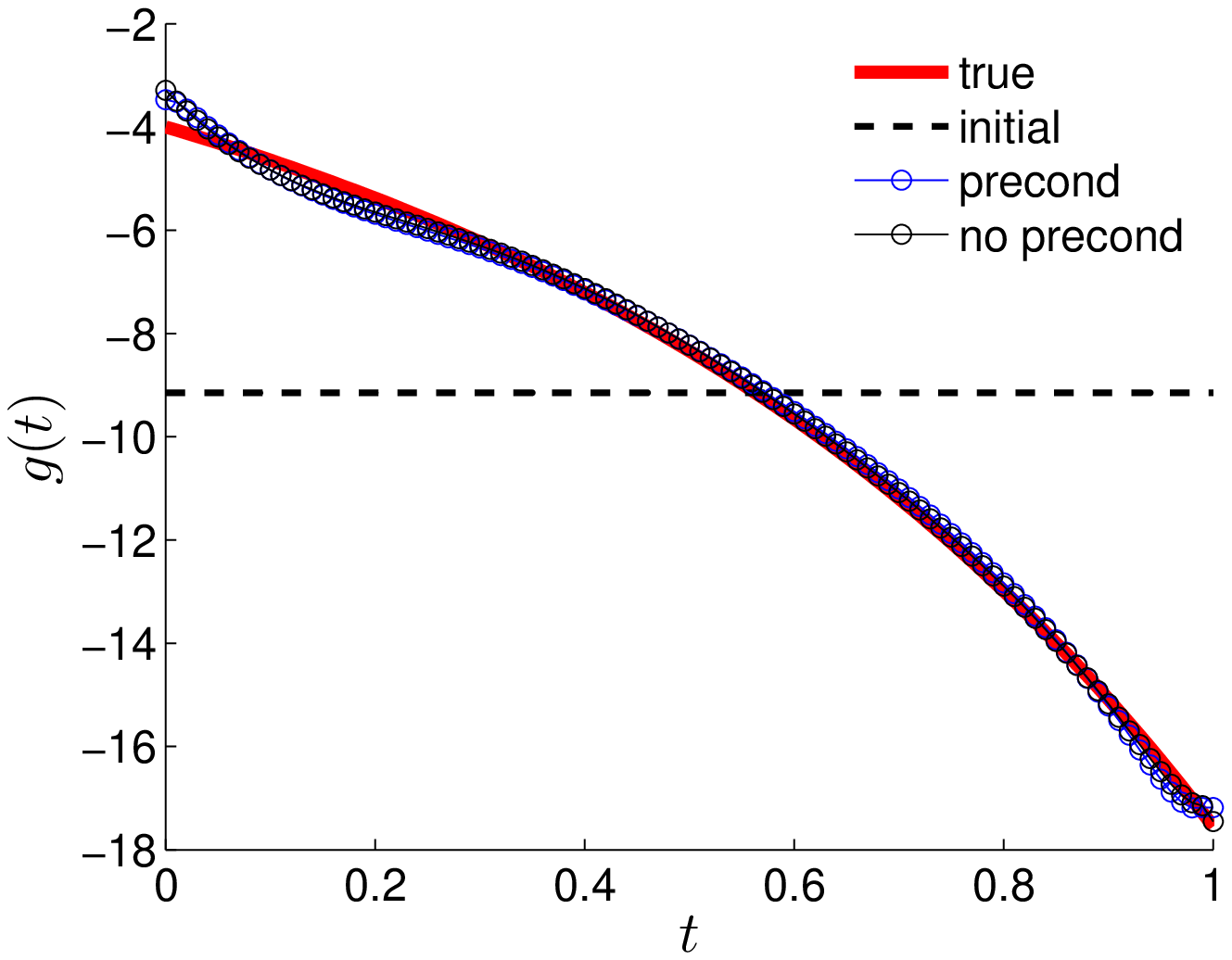}}
  \subfigure[model \#2: $g(t)$]{\includegraphics[width=0.33\textwidth]{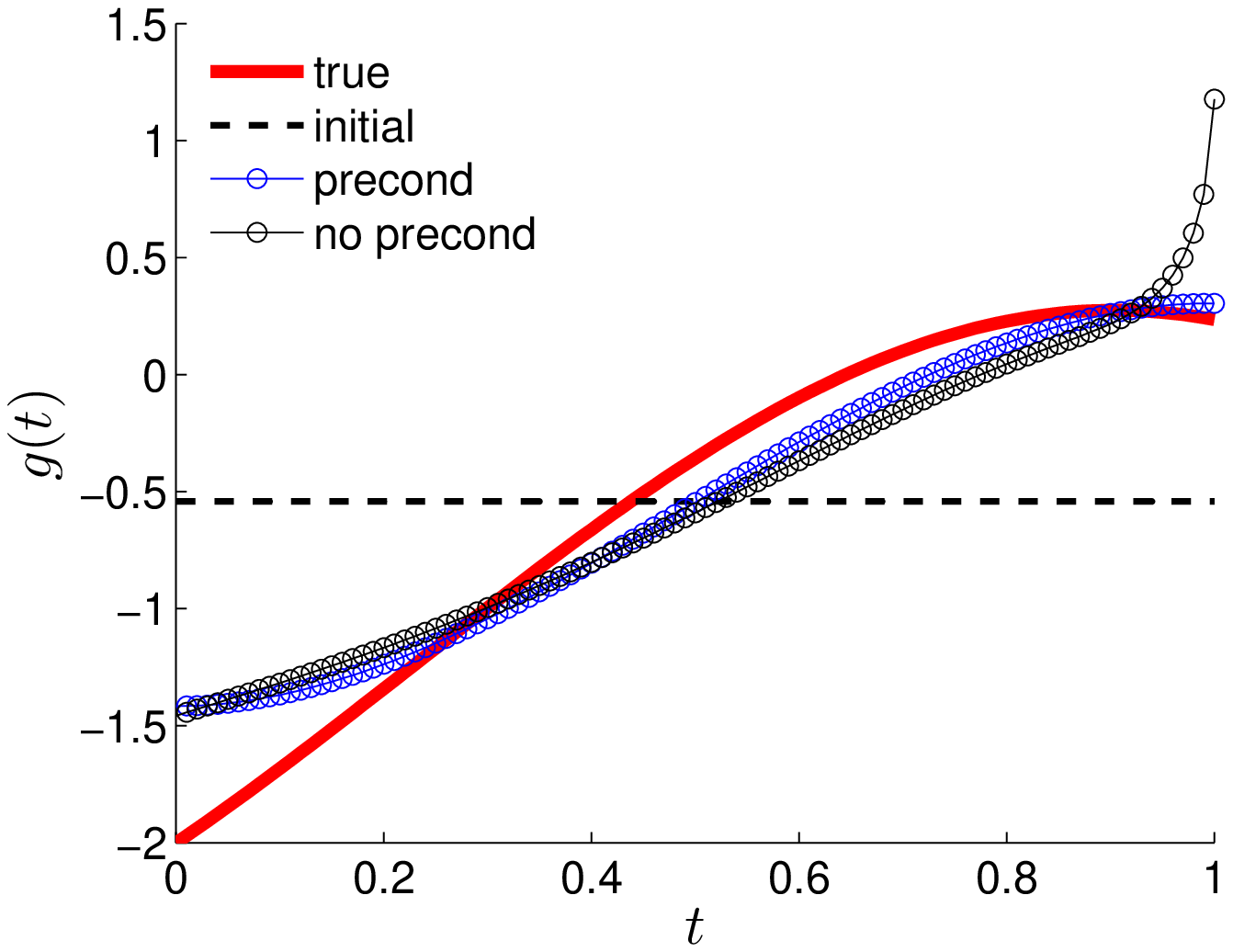}}
  \subfigure[model \#3: $g(t)$]{\includegraphics[width=0.33\textwidth]{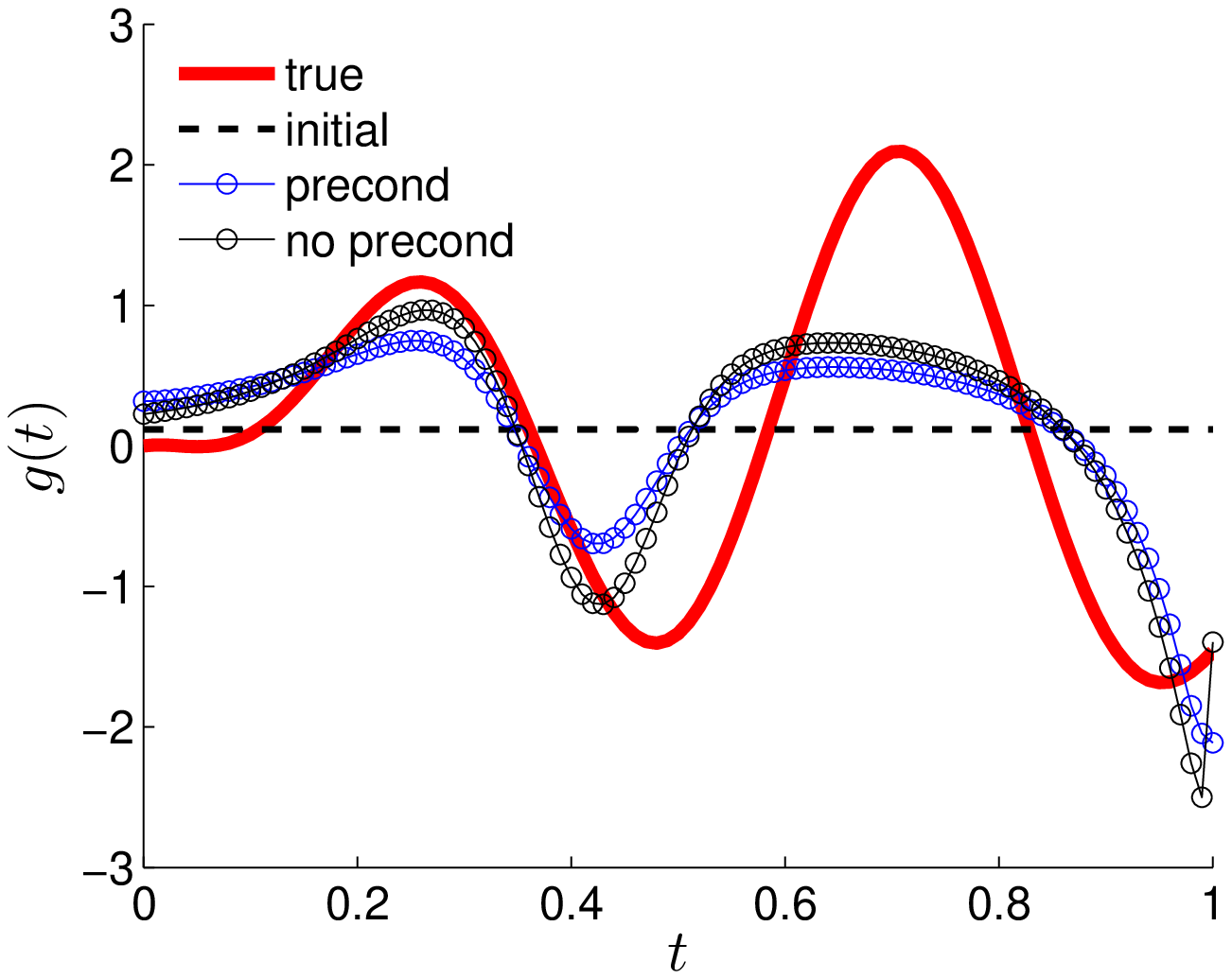}}}
  \mbox{
  \subfigure[model \#1: $w(x) \ \& \ \mu(t)$]{\includegraphics[width=0.33\textwidth]{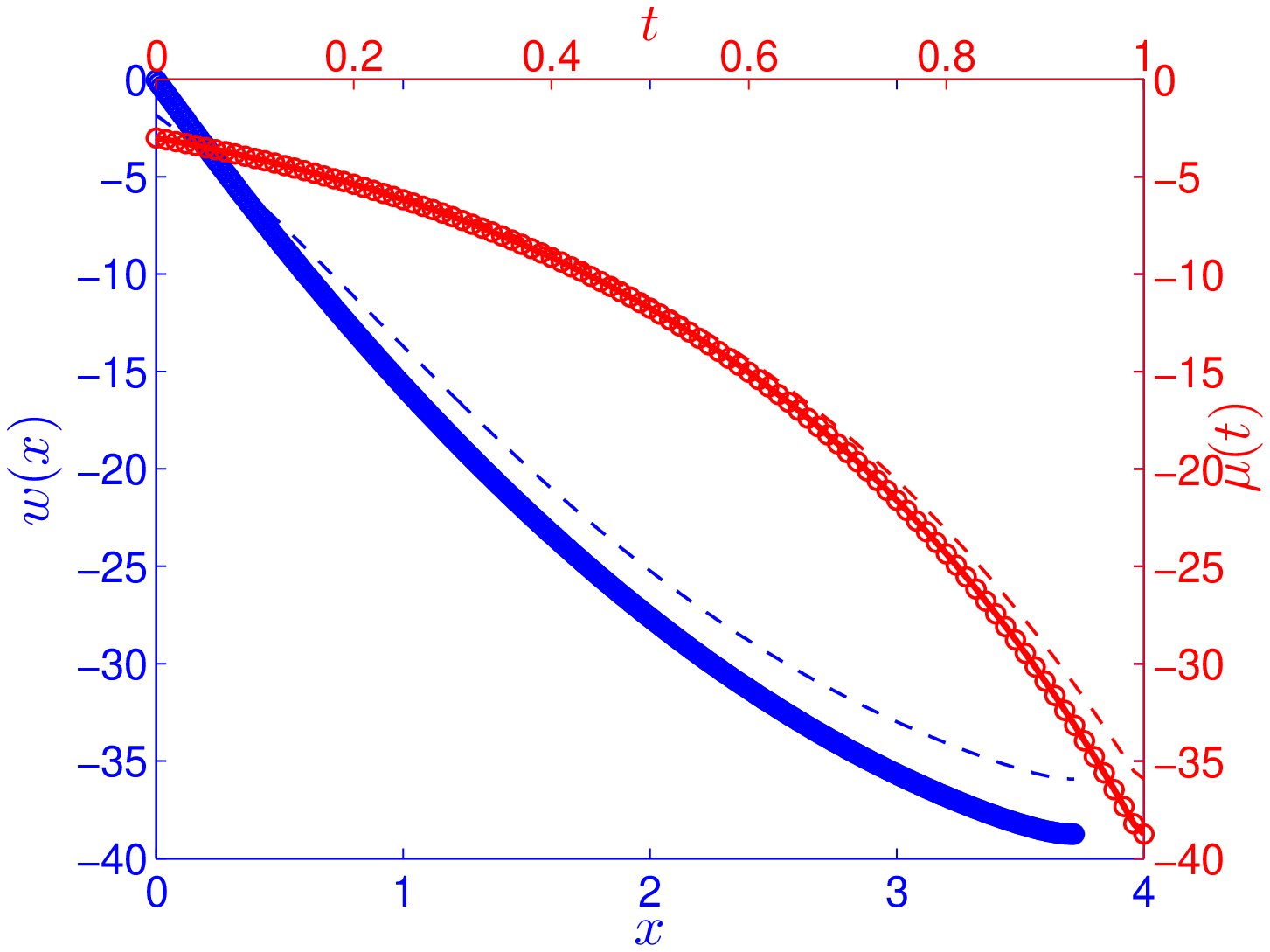}}
  \subfigure[model \#2: $w(x) \ \& \ \mu(t)$]{\includegraphics[width=0.33\textwidth]{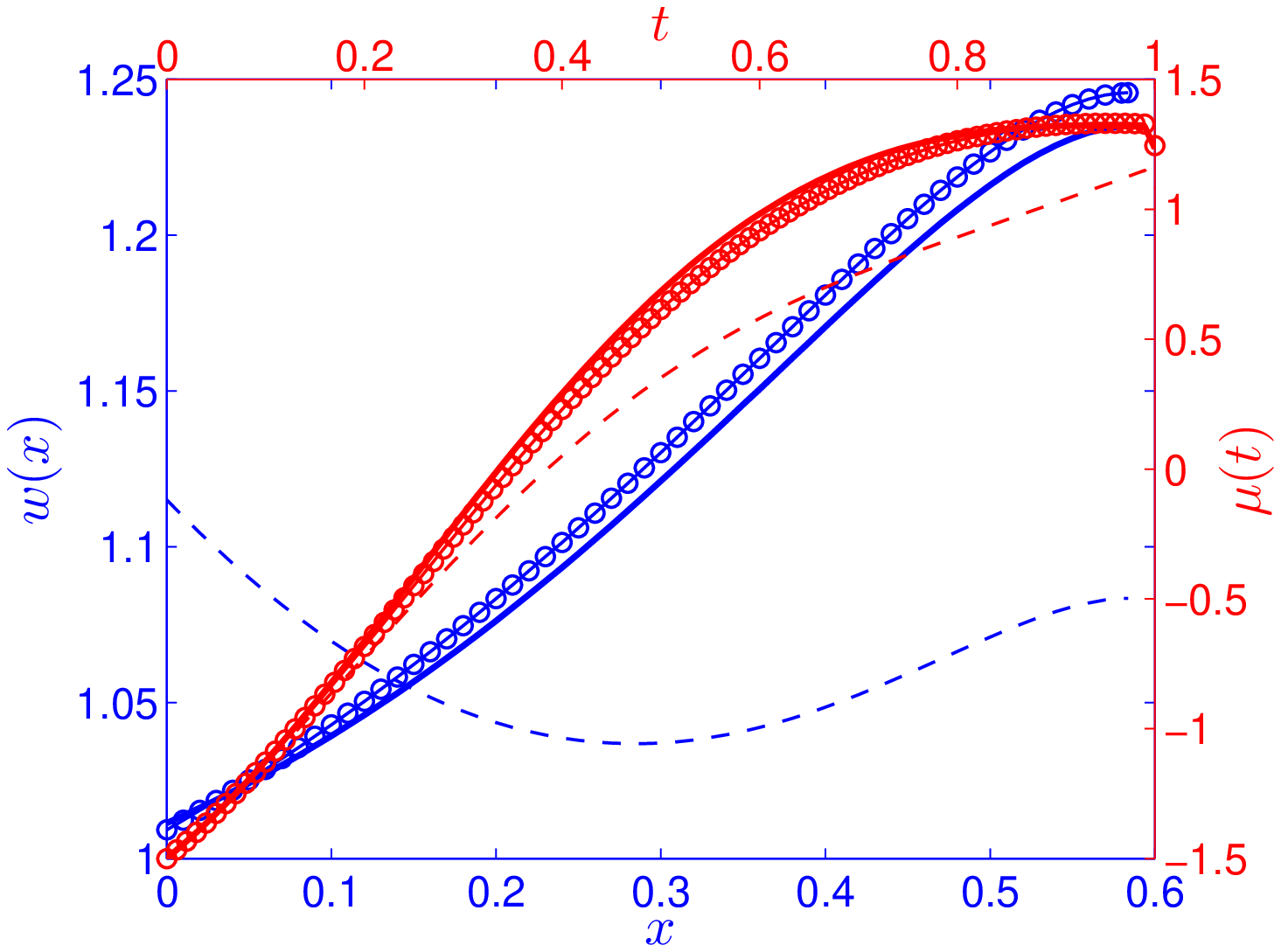}}
  \subfigure[model \#3: $w(x) \ \& \ \mu(t)$]{\includegraphics[width=0.33\textwidth]{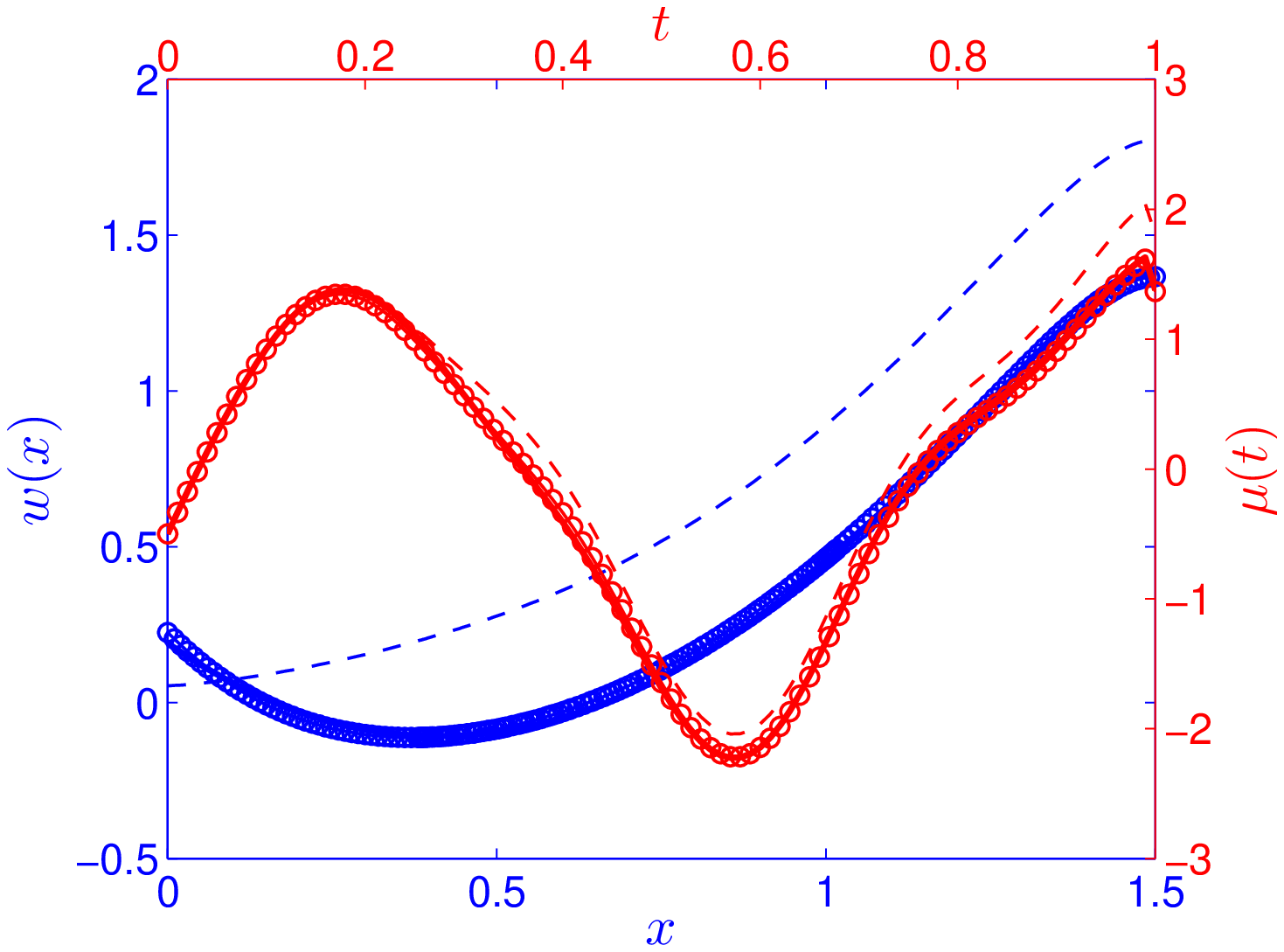}}}
  \mbox{
  \subfigure[model \#1: $\mJ_k$]{\includegraphics[width=0.33\textwidth]{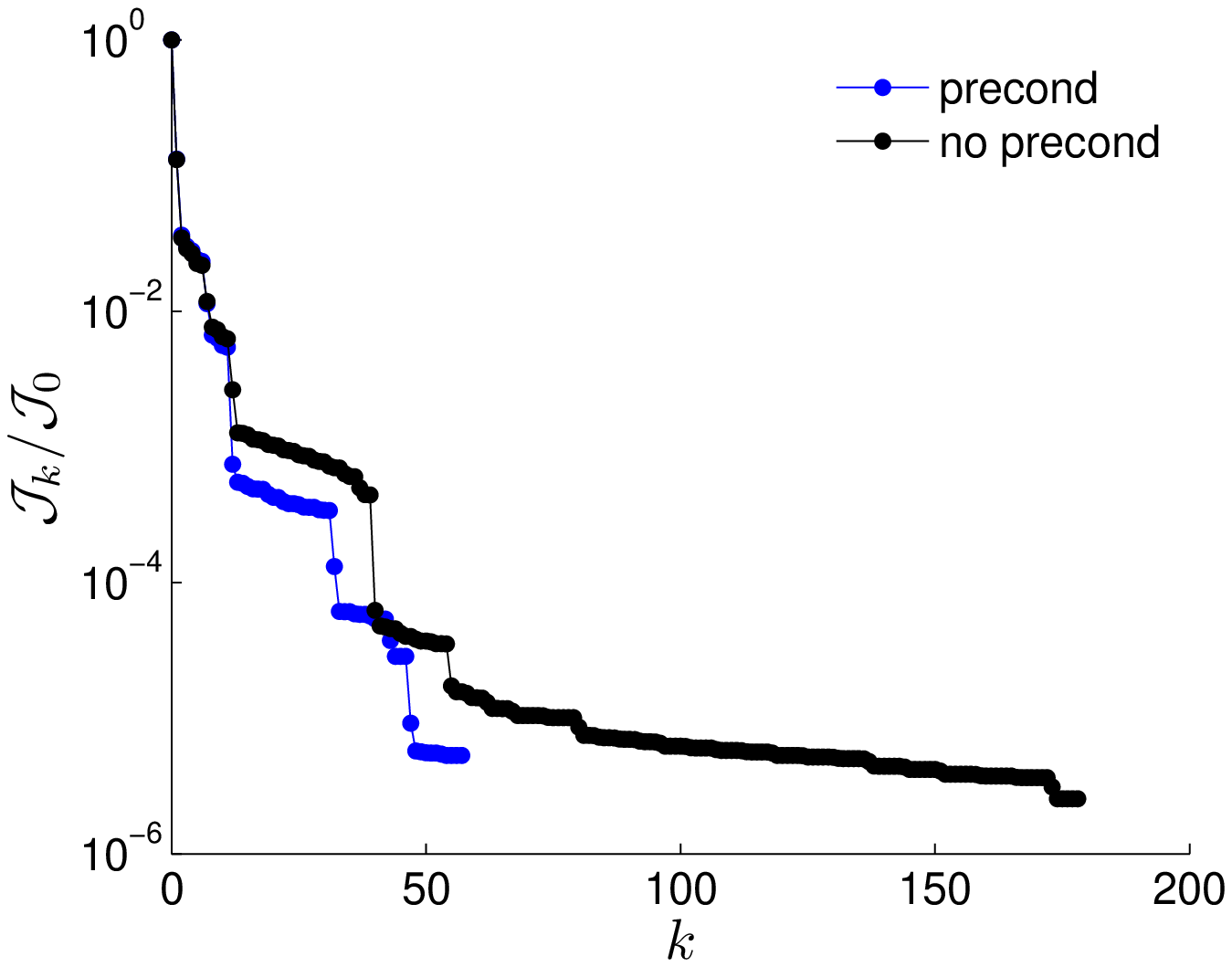}}
  \subfigure[model \#2: $\mJ_k$]{\includegraphics[width=0.33\textwidth]{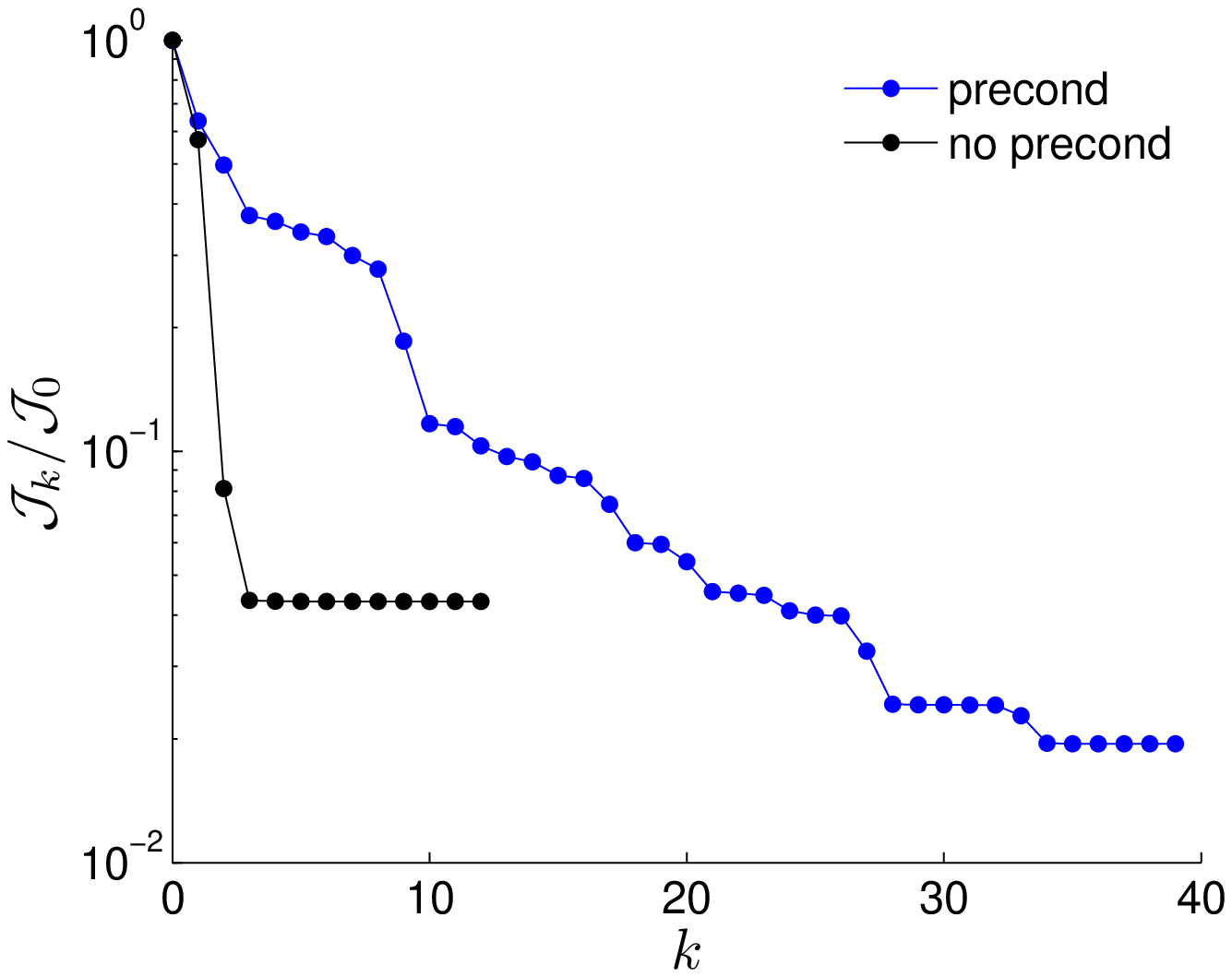}}
  \subfigure[model \#3: $\mJ_k$]{\includegraphics[width=0.33\textwidth]{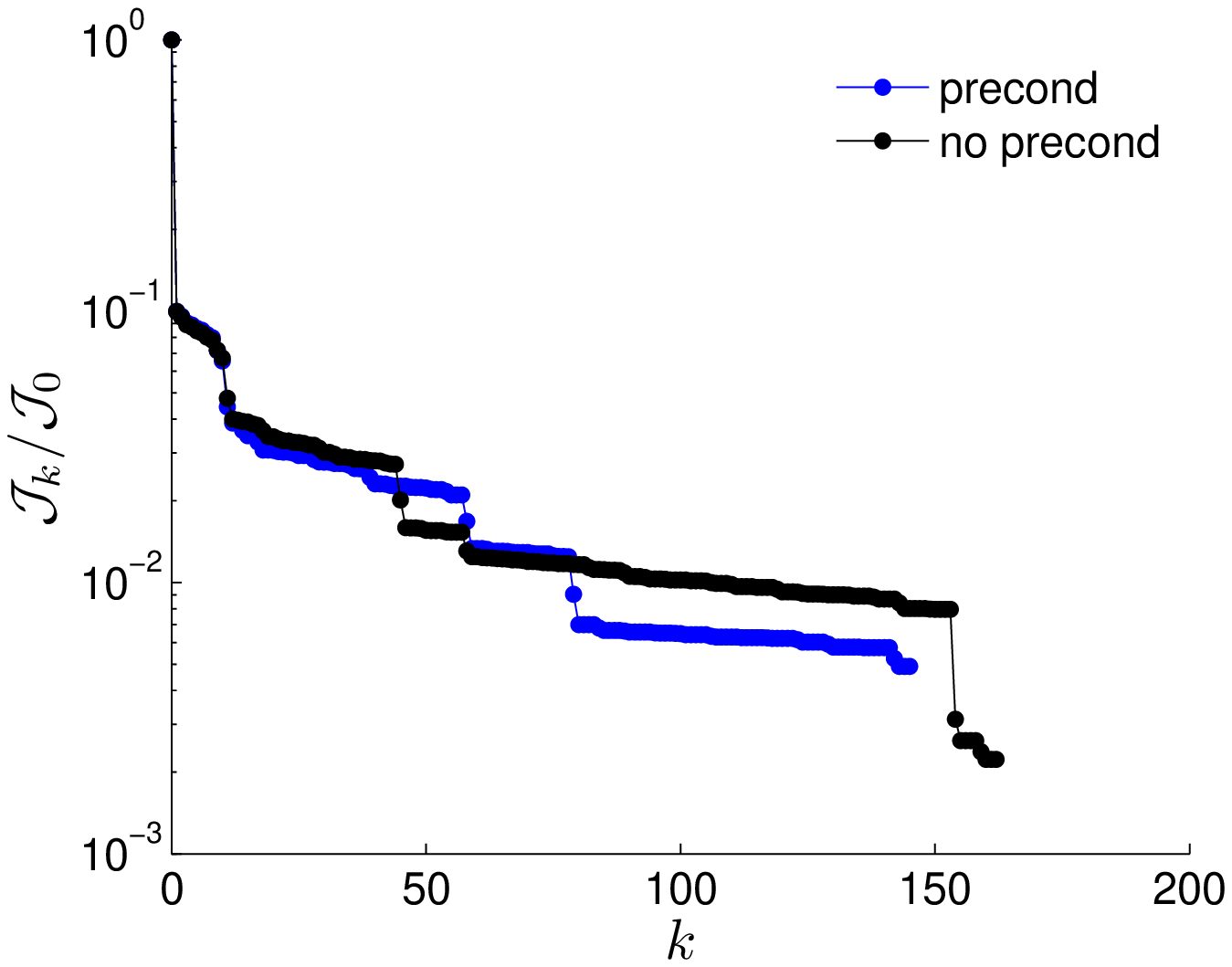}}}
  \end{center}
  \caption{Results of reconstructing left boundary heat flux $g(t)$ for (a,d,g) model \#1, (b,e,h) model \#2, and
    (c,f,i) model \#3. In (a-c) solid red and dashed black lines show respectively the shape of $g_{\rm{true}}(t)$
    and initial guess $g_{\rm{ini}}(t)$, while blue and black circles represent optimal solutions $\bg(t)$
    respectively with and without preconditioning. In (d-f) blue and red colors are used correspondingly
    for functions $w(x)$ and $\mu(t)$: dashed and solid lines represent correspondingly their
    values when $g = g_{\rm{ini}}(t)$ and $g = g_{\rm{true}}(t)$, while circles are used when $s = \bg(t)$
    obtained with preconditioning. In (g-i) blue and black dots show normalized cost functionals
    $\mJ_k/\mJ_0$ as functions of iteration number $k$.}
  \label{fig:opt_ell_g}
\end{figure}

Figure~\ref{fig:opt_ell_g}(a-c) shows the outcomes of reconstructing left boundary heat flux $g(t)$ for
all three models comparing the results obtained with and without preconditioning (blue and black circles
respectively). Preconditioning procedure for model~\#2 uses $\ell^{\star}_{g,2} = 0.6$ obtained by finding
the minimal value of $\mJ$ in the proximity of approximated $\ell^*_g$. This value is shown by blue hexagon
in Figure~\ref{fig:precond_g}(b). Preconditioning procedures for models \#1 and \#3 use
$\ell^{\star}_{g,1} = \ell^{\star}_{g,3} = 10^{-2}$ as discussed before. The quality of reconstruction
of $g(t)$ depends obviously on the complexity of the model. Simple model~\#1 shows very accurate results,
while models \#2 and \#3 are stuck on the local solutions $\bg(t)$. But even in the absence of perfect match,
these solutions are close to true functions $g_{\rm{true}}$. Such results could be explained by non-uniqueness
of the solved problem of finding heat flux at a left boundary based on measurements obtained at final time
$T$ and at free (right-side) boundary $s(t)$. This non-uniqueness is also justified by observing how accurately
the obtained solutions $u(x,T)$ and $u(s(t),t)$ match the measurements $w(x)$ and $\mu(t)$ which is seen in
Figure~\ref{fig:opt_ell_g}(d-f). In Figure~\ref{fig:opt_ell_g}(g-i) normalized cost functionals $\mJ_k/\mJ_0$
are represented as functions of iteration number $k$. As could be noted here, identification of heat flux $g(t)$
performed separately from free boundary $s(t)$ requires more optimization iterations to reach the same
termination condition $\left| \frac{\mJ_k-\mJ_{k-1}}{\mJ_{k-1}} \right| < 10^{-5}$. We find this observation
useful while discussing further results of simultaneous reconstruction $s(t)$ and $g(t)$.

The last series of our computational results shows identification of full control vector $v = (s(t), g(t))$
by using the proposed approach outlined in Algorithm~\ref{alg:gen_opt}. As it is concluded previously in the
current section and also in Section~\ref{sec:single_rec}, use of the preconditioning procedure \eqref{eq:helm}
provides much better performance in reconstructing both $s(t)$ and $g(t)$ when it is done
separately. Thus, we keep this technique on while obtaining the rest results shown in this section. As seen in
Figure~\ref{fig:opt_ell_s}(g,h,i) and Figure~\ref{fig:opt_ell_g}(g,h,i), cost functional $\mJ$ is more
sensitive to changes in the free (right) boundary $s(t)$ rather than in the left boundary heat flux $g(t)$.
Such difference in the sensitivity of $\mJ$ results in different rate of convergence for $s(t)$ and $g(t)$.
Due to this fact, we would like to compare the results obtained with different strategies for finding
optimal values for stepsize parameters $\alpha_k^s$ and $\alpha_k^g$ used in iterative descent gradient
procedure \eqref{eq:opt_SD}. In order to solve problem \eqref{eq:opt_alpha} we use the following three approaches:
\begin{description}
  \item[{\bf \#1.}] {\it Simultaneous} identification of $s(t)$ and $g(t)$ by setting $\alpha_k = \alpha_k^s = \alpha_k^g$
    while solving one-dimensional optimization problem \eqref{eq:opt_alpha} and updating both $s_k(t)$ and $g_k(t)$
    within the same $k$-th optimization iteration.
  \item[{\bf \#2.}] Identification of $s(t)$ and $g(t)$ in the {\it interchanging} order when only one control is updated
    during $k$-th optimization iteration. In other words, when $k=2n-1, \, n = 1, 2, \ldots$, we set $\alpha_k^g = 0$
    and solve \eqref{eq:opt_alpha} to find $\alpha_k = \left( \alpha_k^s, 0 \right)$ and update only $s_k(t)$
    using \eqref{eq:opt_SD_s}. Then, similarly, for $k=2n, \, n = 1, 2, \ldots$, we set $\alpha_k^s = 0$
    and solve \eqref{eq:opt_alpha} to find $\alpha_k = \left( 0, \alpha_k^g \right)$ and update only $g_k(t)$
    using \eqref{eq:opt_SD_g}.
  \item[{\bf \#3.}] Identification of $s(t)$ and $g(t)$ in the {\it $N$-interchanging} order, or using so-called
    {\it nested optimization}. This strategy utilizes the same approach \#2 to update only one control at a time,
    but changing controls every $N$ optimization iterations. In fact, approach \#2 could be seen as a method of
    the same kind when $N = 1$.
\end{description}

\begin{figure}[htb!]
  \begin{center}
  \mbox{
  \subfigure[model \#1: $s(t)$]{\includegraphics[width=0.33\textwidth]{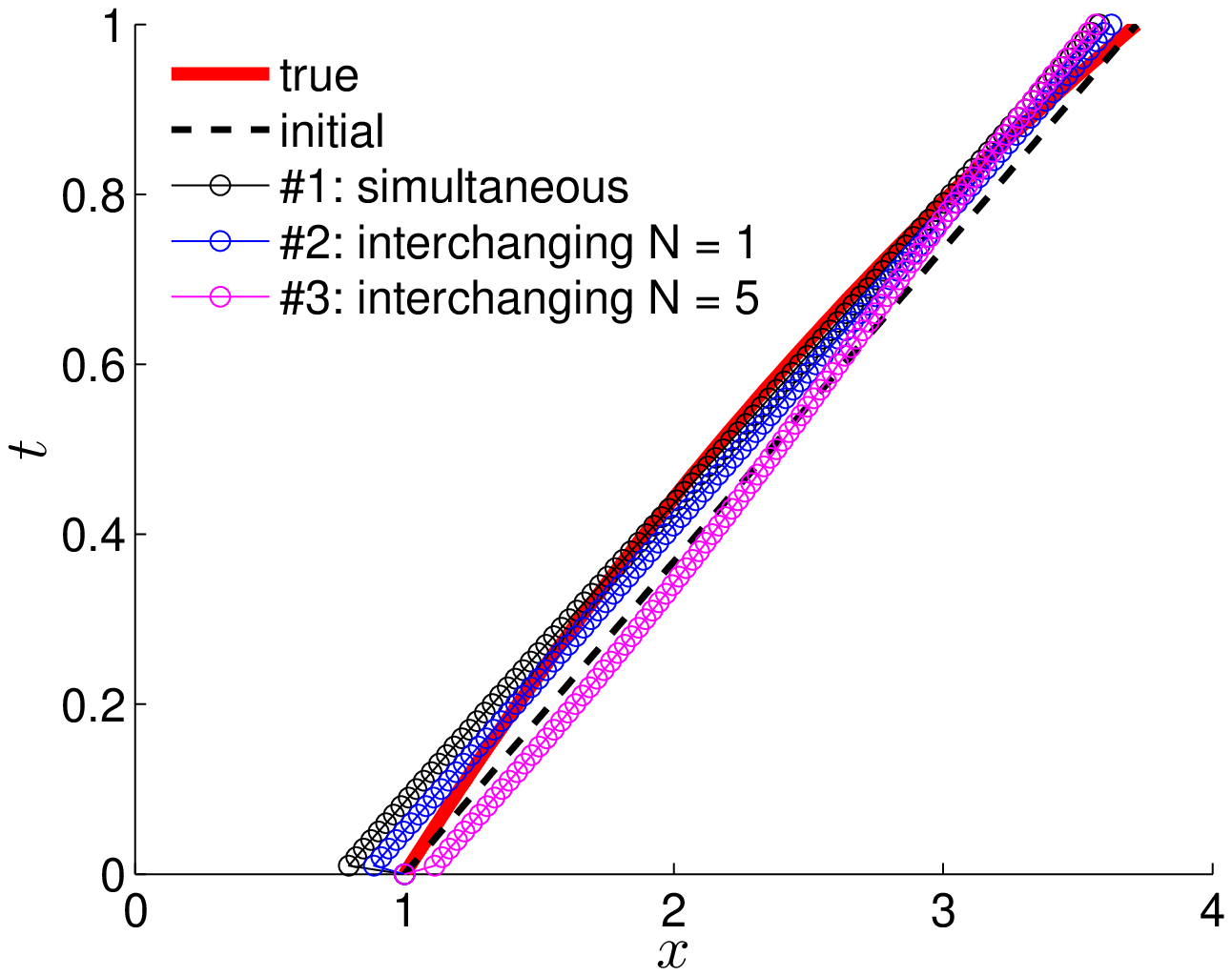}}
  \subfigure[model \#2: $s(t)$]{\includegraphics[width=0.33\textwidth]{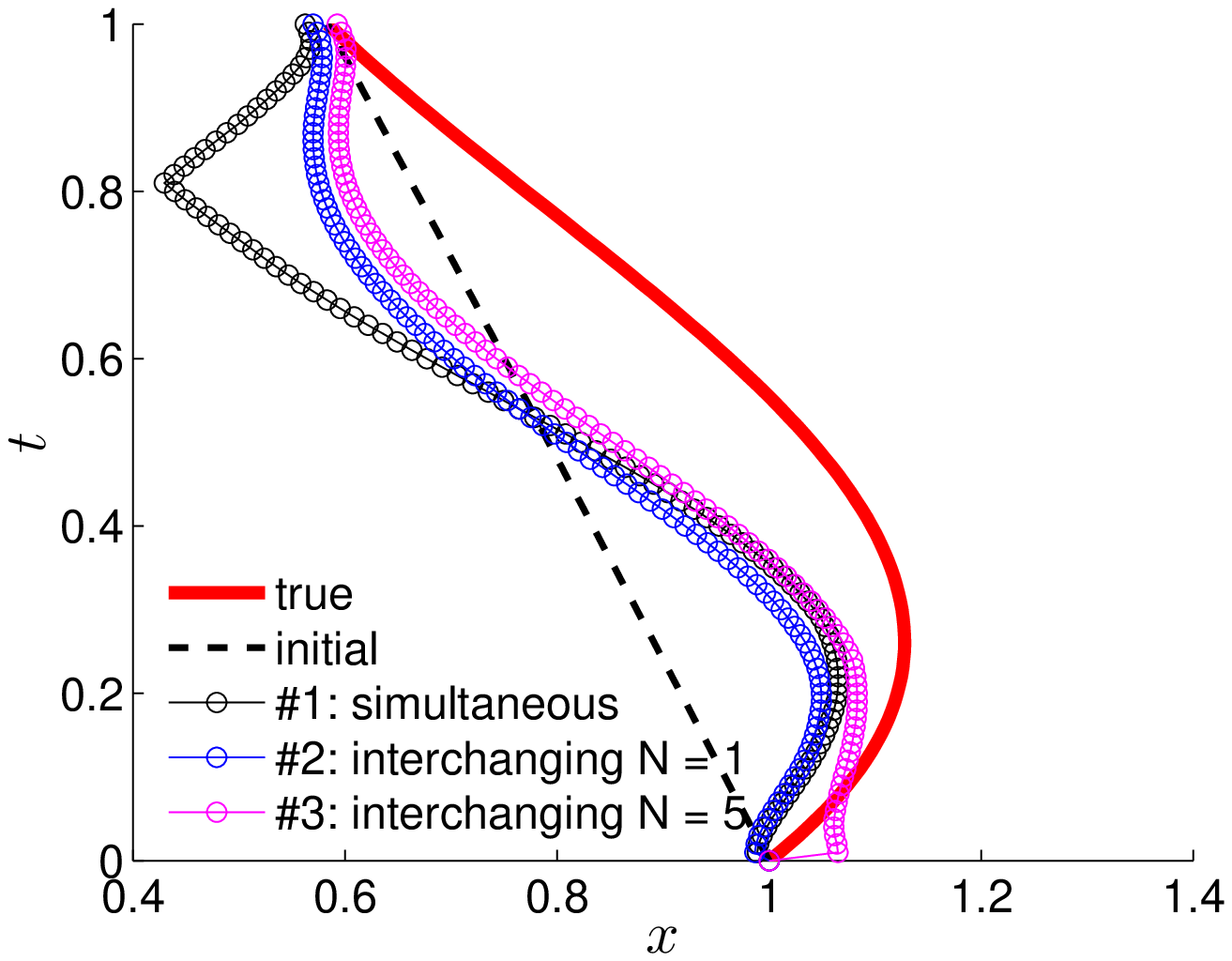}}
  \subfigure[model \#3: $s(t)$]{\includegraphics[width=0.33\textwidth]{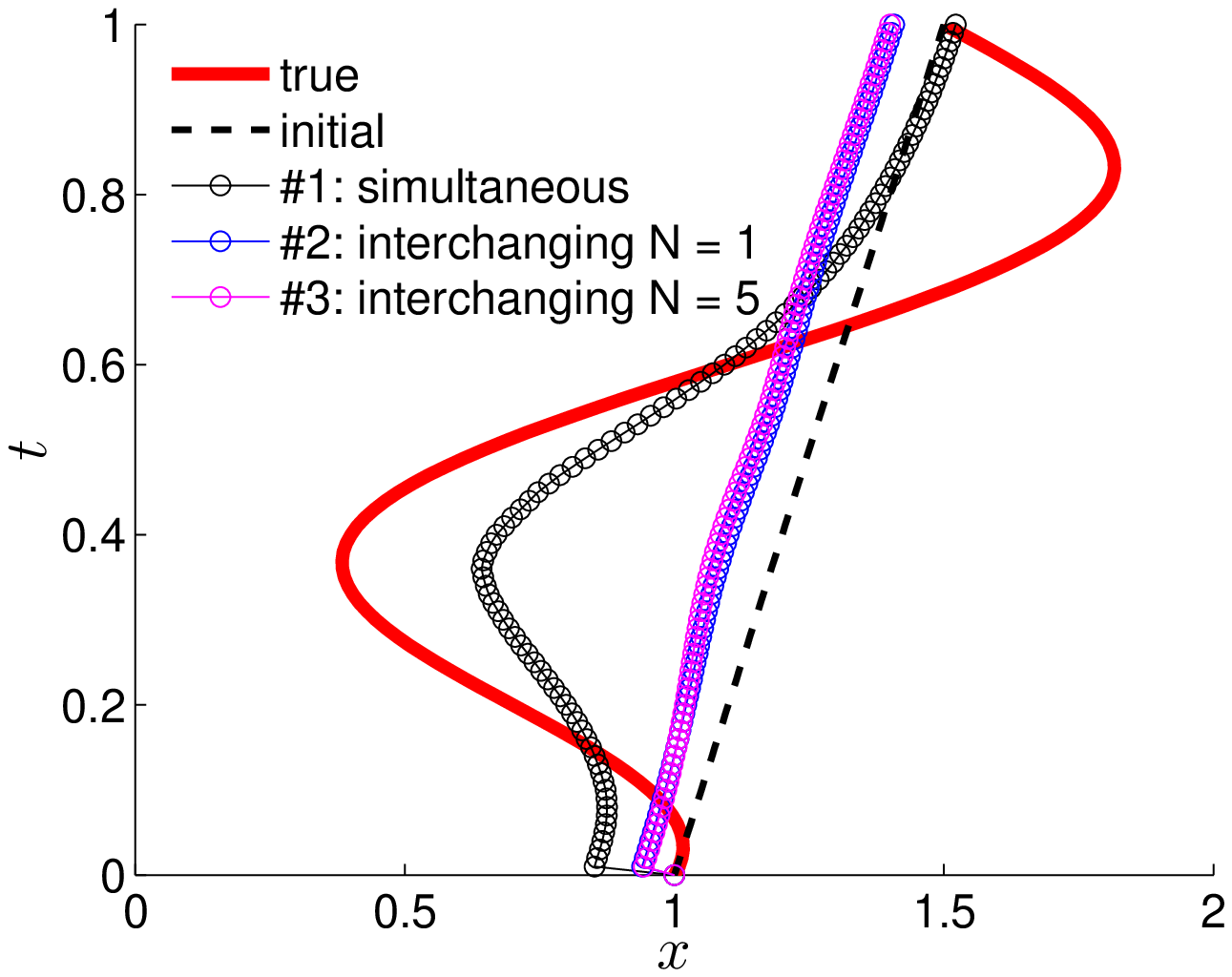}}}
  \mbox{
  \subfigure[model \#1: $g(t)$]{\includegraphics[width=0.33\textwidth]{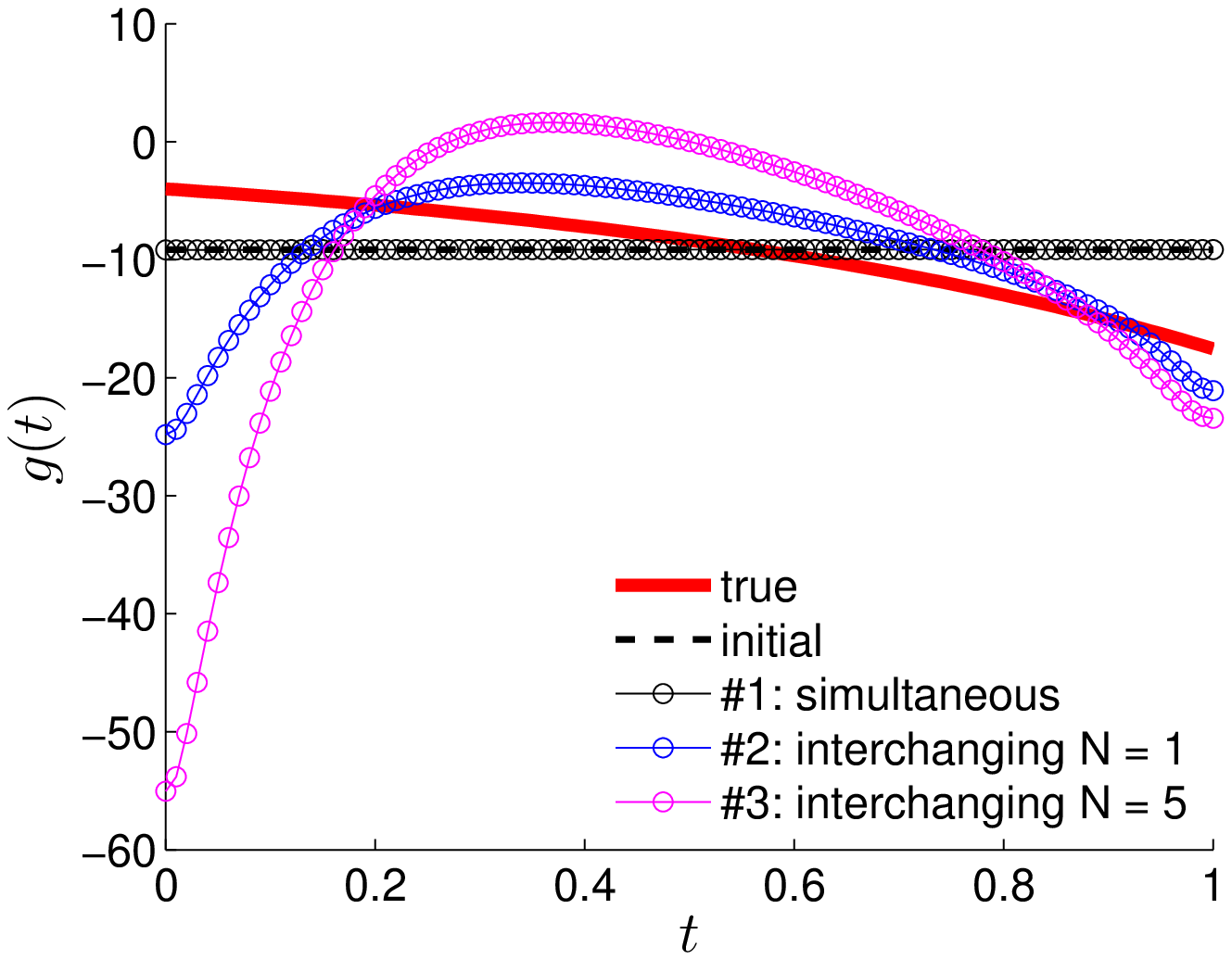}}
  \subfigure[model \#2: $g(t)$]{\includegraphics[width=0.33\textwidth]{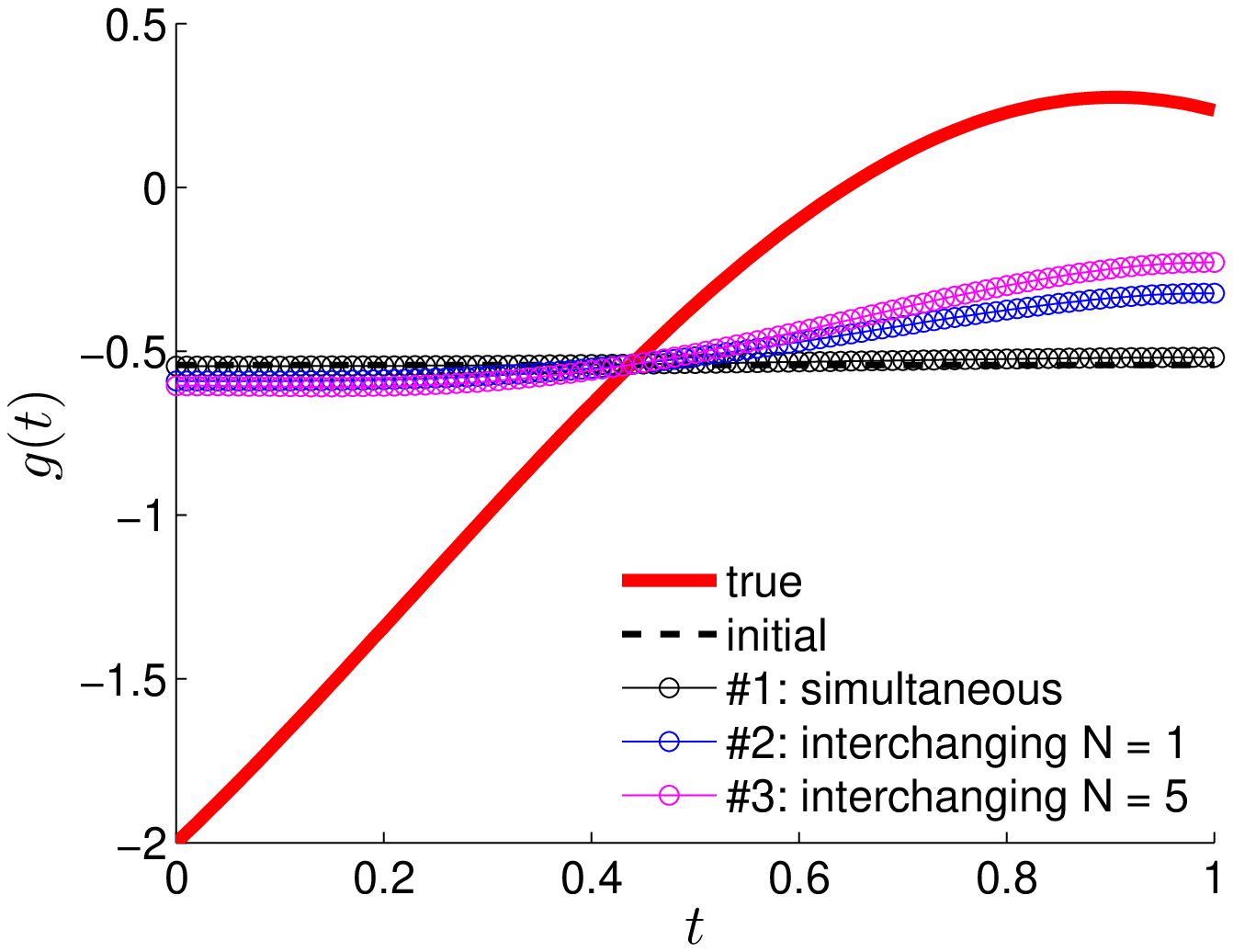}}
  \subfigure[model \#3: $g(t)$]{\includegraphics[width=0.33\textwidth]{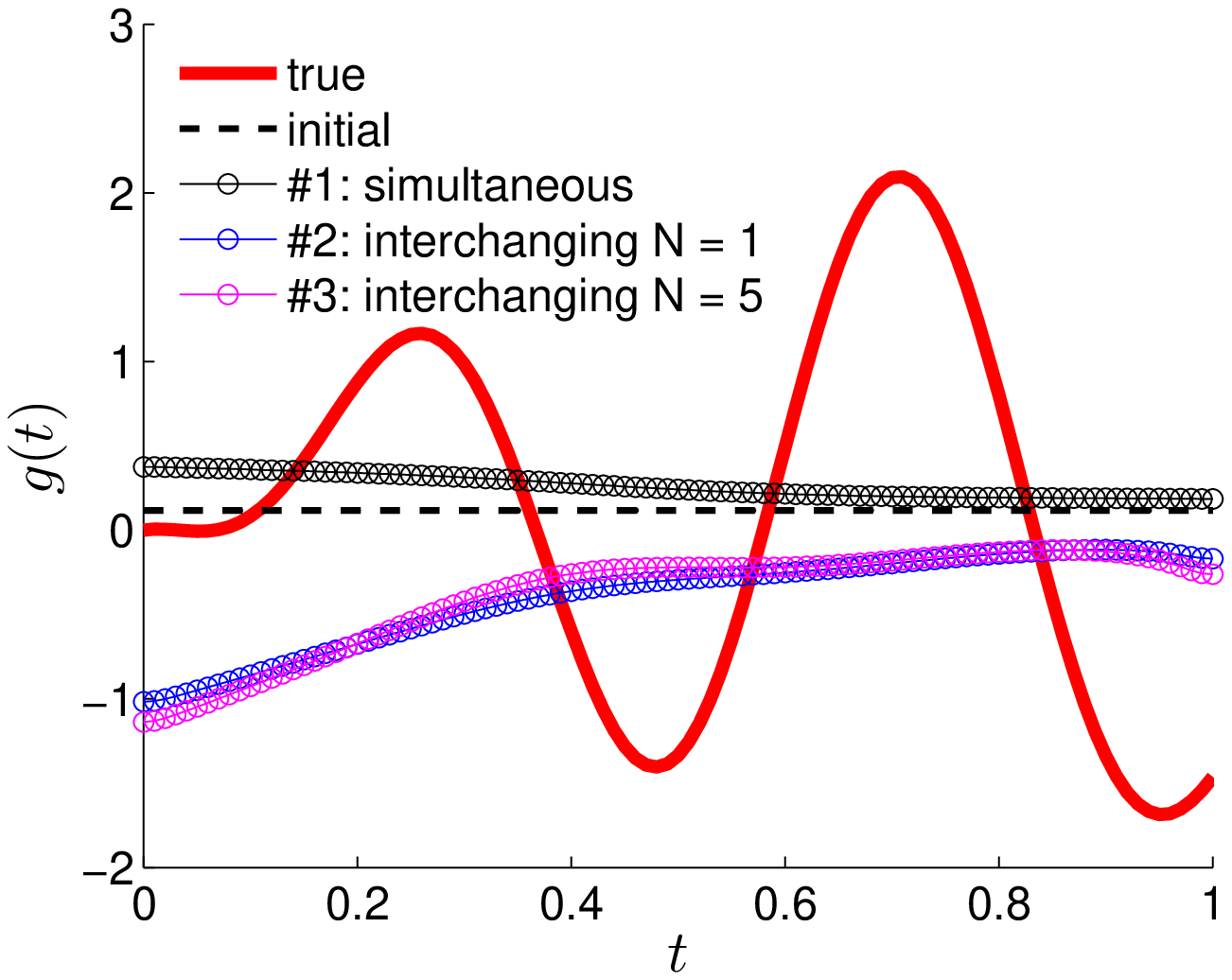}}}
  \mbox{
  \subfigure[model \#1: $w(x) \ \& \ \mu(t)$]{\includegraphics[width=0.33\textwidth]{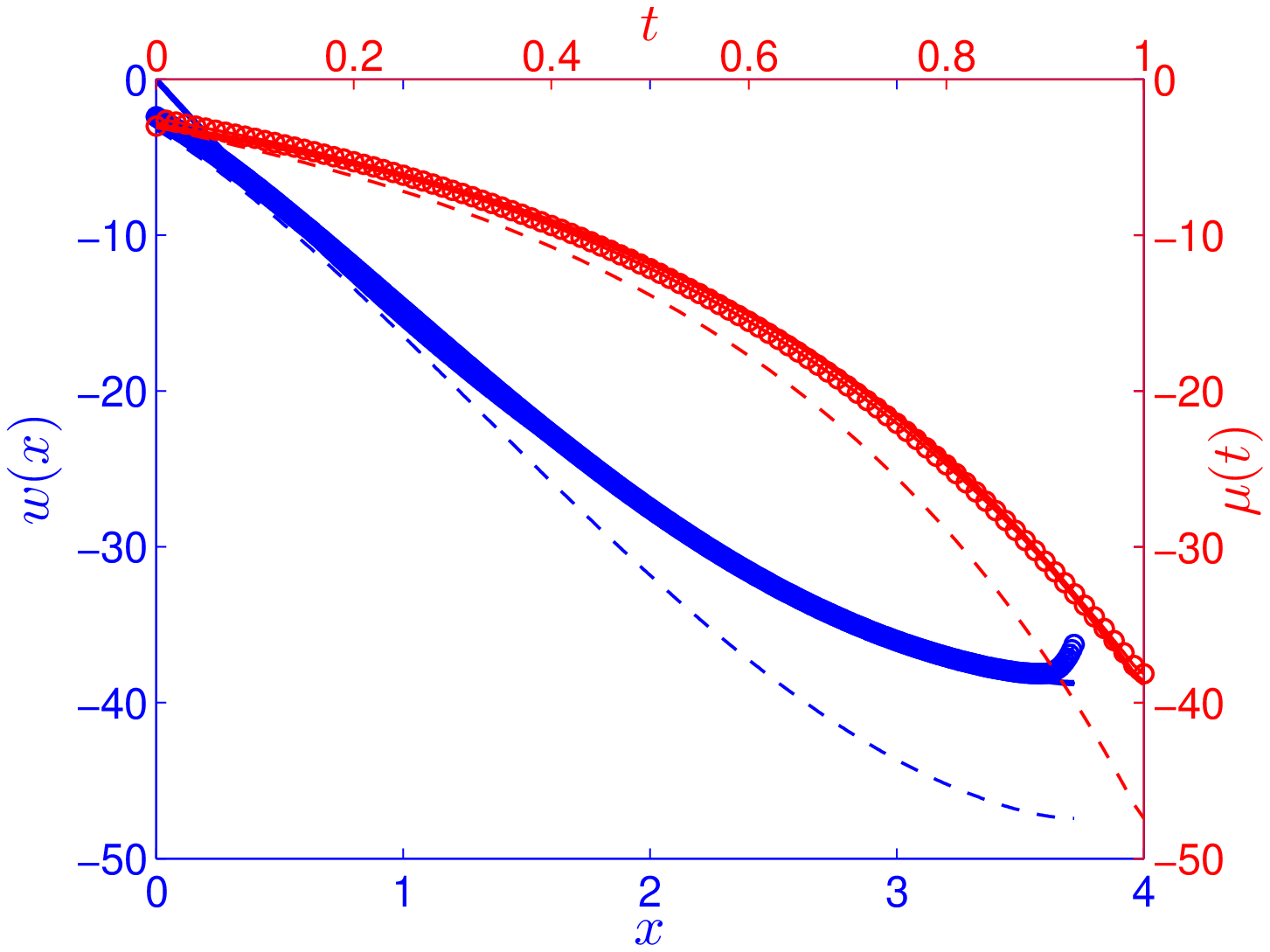}}
  \subfigure[model \#2: $w(x) \ \& \ \mu(t)$]{\includegraphics[width=0.33\textwidth]{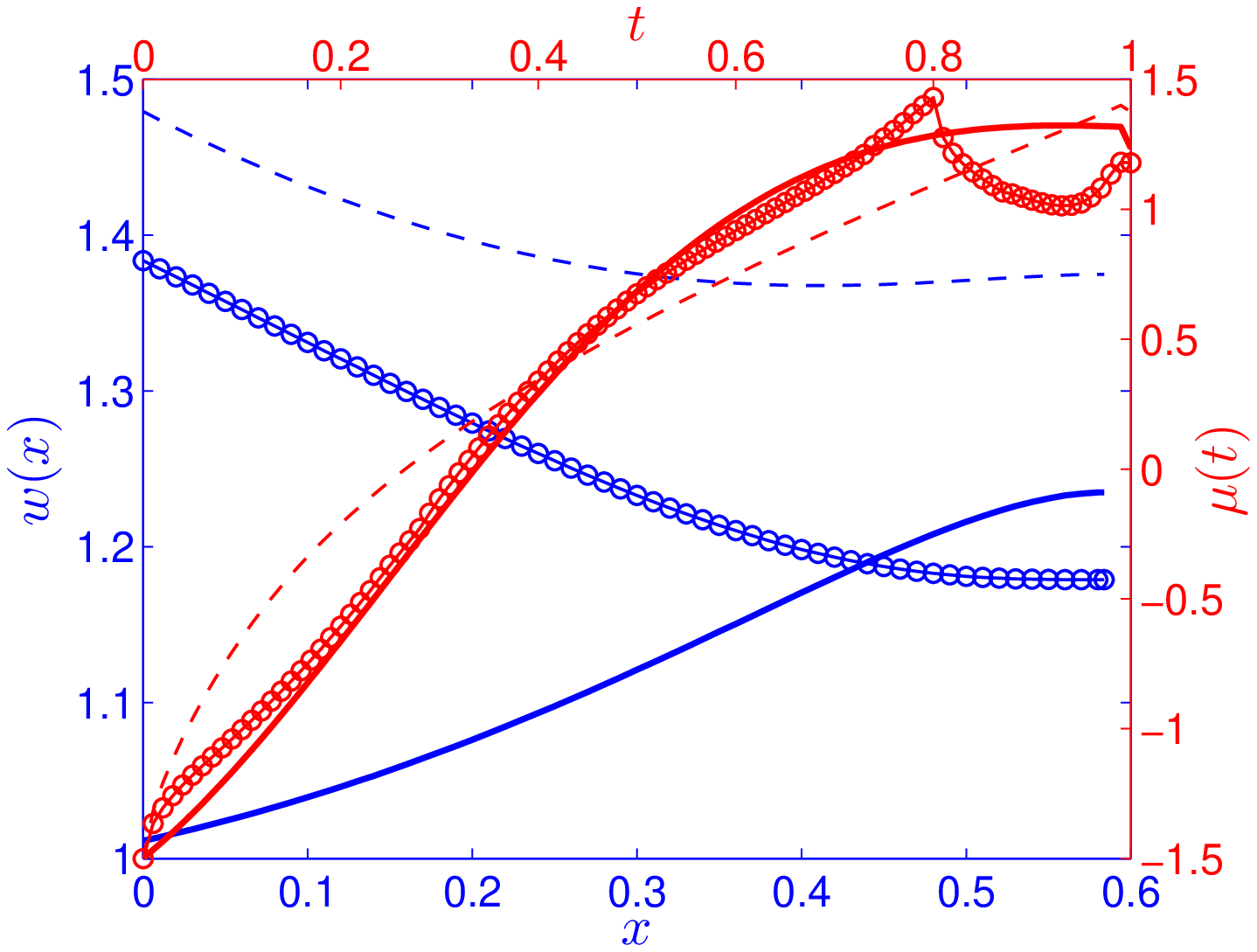}}
  \subfigure[model \#3: $w(x) \ \& \ \mu(t)$]{\includegraphics[width=0.33\textwidth]{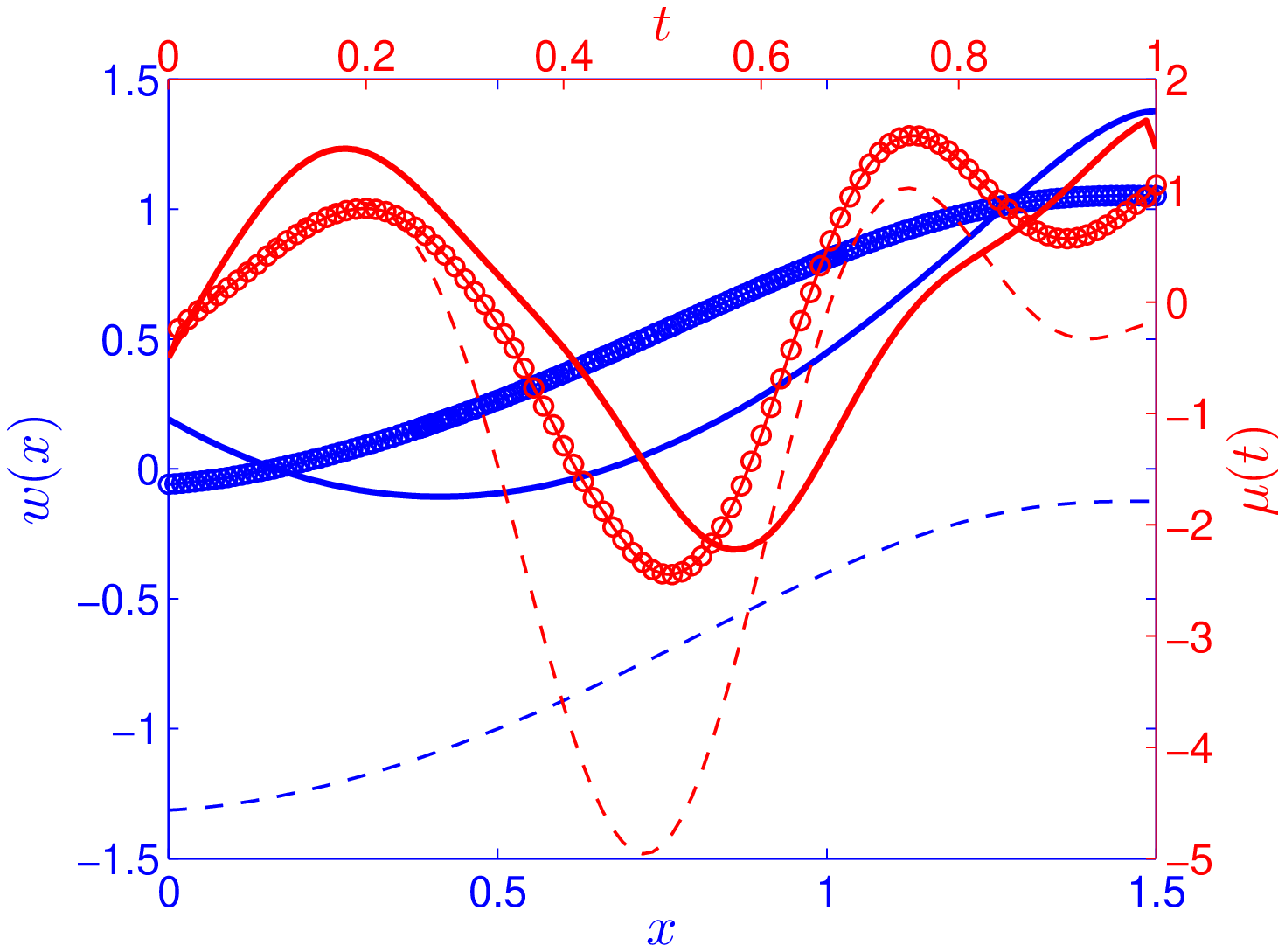}}}
  \end{center}
  \caption{Results of reconstructing (a,b,c) free boundary $s(t)$ and (d,e,f) left boundary heat flux $g(t)$
    for (a,d,g) model~\#1, (b,e,h) model~\#2, and (c,f,i) model~\#3. In (a-f) solid red and dashed black lines show
    respectively the shapes of true functions and their initial guesses; while black, blue and purple circles represent
    respectively optimal solutions for simultaneous reconstruction of $v = (s(t),g(t))$, and using interchanging order
    with $N = 1$ and $N = 5$. In (g-i) blue and red colors are used correspondingly for functions $w(x)$ and $\mu(t)$:
    dashed and solid lines represent correspondingly their values when $v = v_{\rm{ini}}$ and $v = v_{\rm{true}}$,
    while circles are used when $v = \bv$ obtained with simultaneous reconstruction (approach \#1).}
  \label{fig:opt_ell_sg}
\end{figure}

Figure~\ref{fig:opt_ell_sg} shows the results of identification both $s(t)$ and $g(t)$ for all three approaches:
simultaneous (black circles), and interchanging order for $N=1$ (blue circles) and $N=5$ (purple circles). We use
preconditioning procedure \eqref{eq:helm} for all three models supplied with $\ell^{\star}_{s,1} = 0.47$,
$\ell^{\star}_{s,2} = 0.2$, $\ell^{\star}_{s,3} = 0.52$ and $\ell^{\star}_{g,1} = 10^{-2}$, $\ell^{\star}_{g,2} = 0.6$,
$\ell^{\star}_{g,3} = 10^{-2}$ obtained previously by finding the best (minimal) values of $\mJ$ in the proximity of approximated
$\ell^*_s$ and $\ell^*_g$.

The results of identifying both $s(t)$ and $g(t)$ are consistent with our previous discussion on the complexity
of our models in particular and the complexity of the inverse Stefan problem in general. Hence, our conclusions
on the overall performance are two-fold. First, the quality of the obtained solution obviously depends
on the complexity of the model. As seen in Figure~\ref{fig:opt_ell_sg}(a,b,c) model~\#1 shows good convergence
for $s(t)$ for all three approaches used, while the results for models \#2 and \#3 are dependent on such
approaches. Interchanging gradients with $N=1$ and $N=5$ works well for model \#2 which is of moderate complexity,
but much better performance is shown by simultaneous gradient use for rather complicated model \#3. At the
same time, Figure~\ref{fig:opt_ell_sg}(d,e,f) shows that the performance in identifying $g(t)$ is poor
for all three models. This fact is consistent with the general statement that any gradient based approach
is sensitive to the choice of the optimization parameters: space and time discretization, initial guess,
smoothing parameter for preconditioning, step size in the control update procedure, and many other parameters
we do not consider in the current paper.

\begin{figure}[htb!]
  \begin{center}
  \mbox{
  \subfigure[model \#1: no preconditioning]{\includegraphics[width=0.33\textwidth]{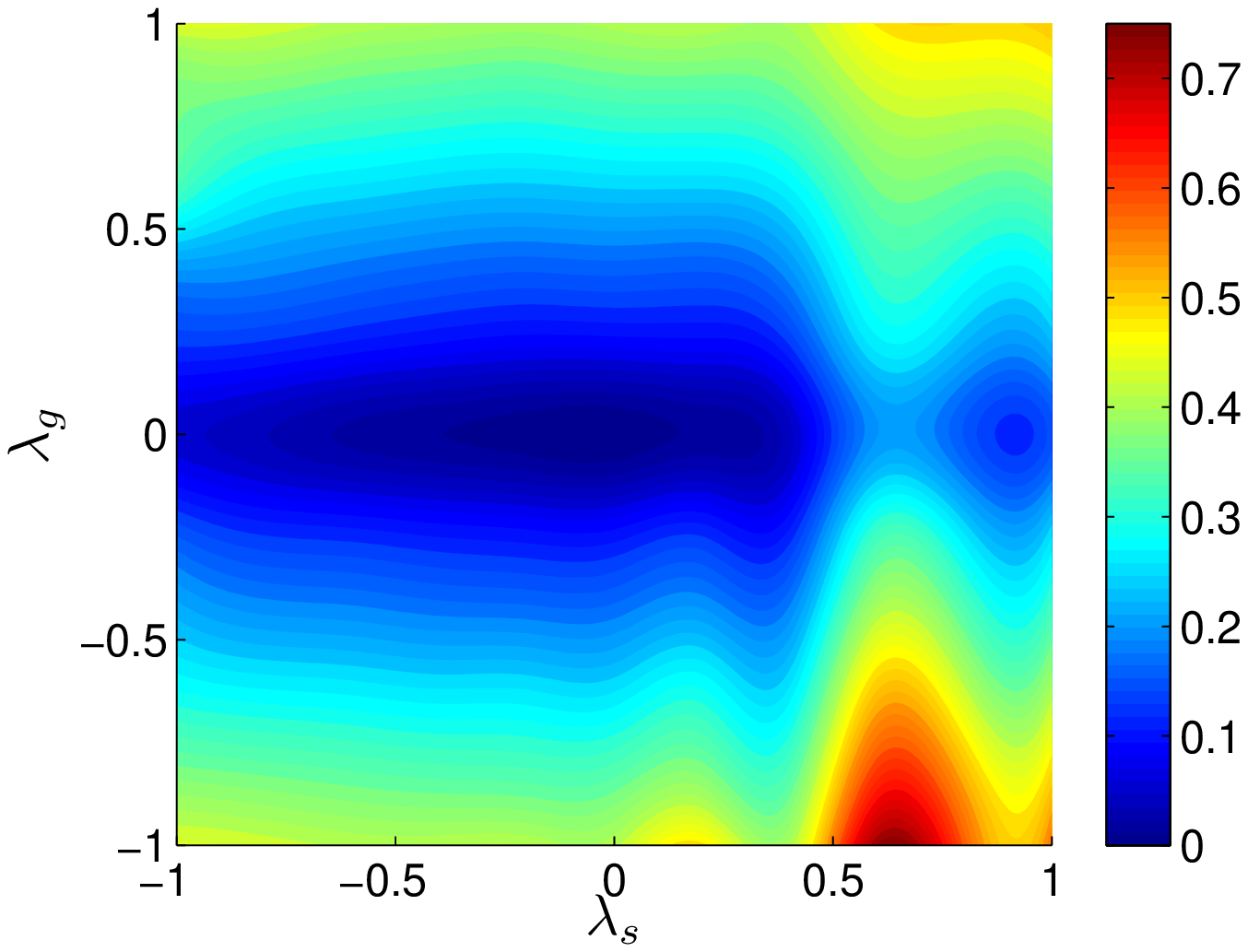}}
  \subfigure[model \#2: no preconditioning]{\includegraphics[width=0.33\textwidth]{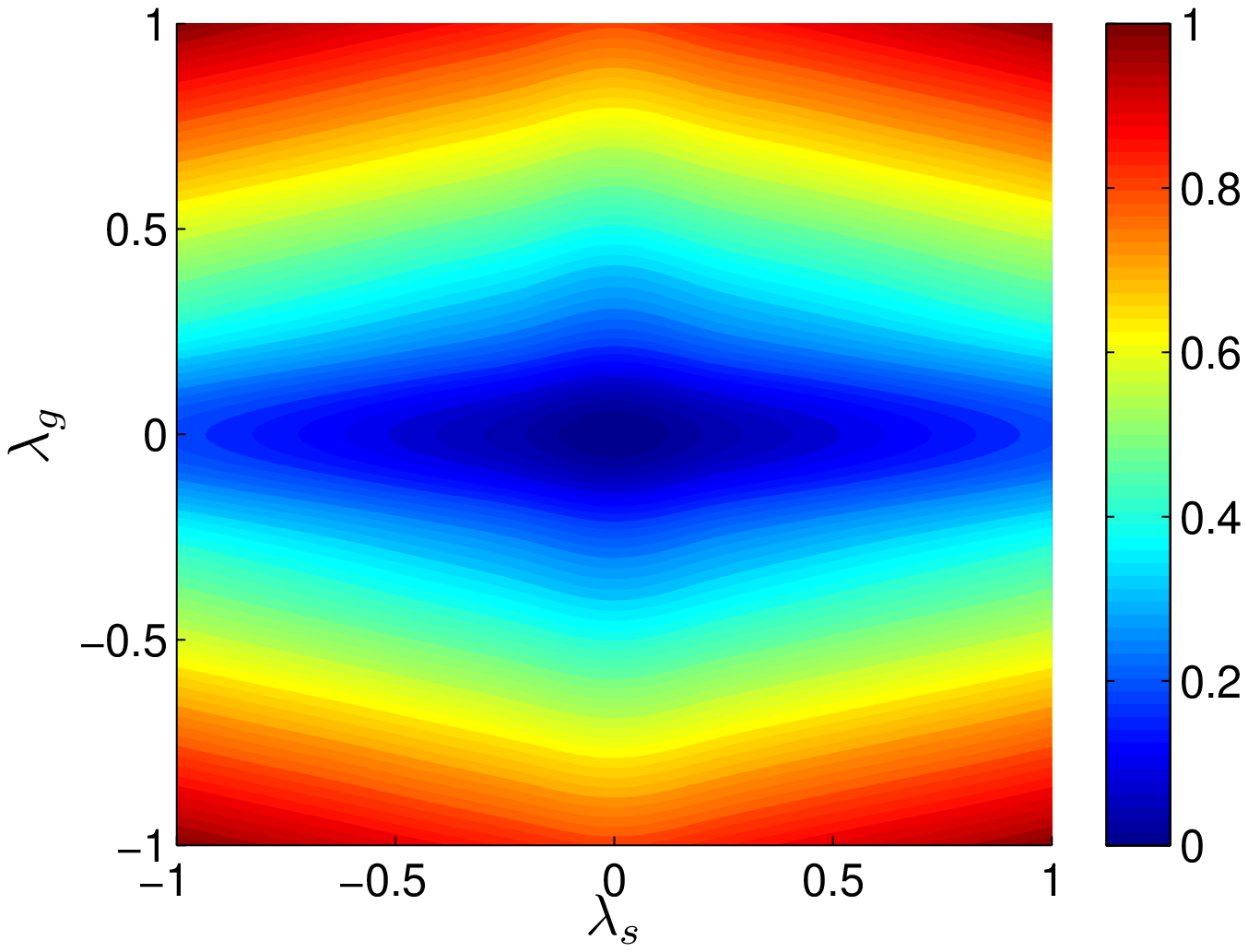}}
  \subfigure[model \#3: no preconditioning]{\includegraphics[width=0.33\textwidth]{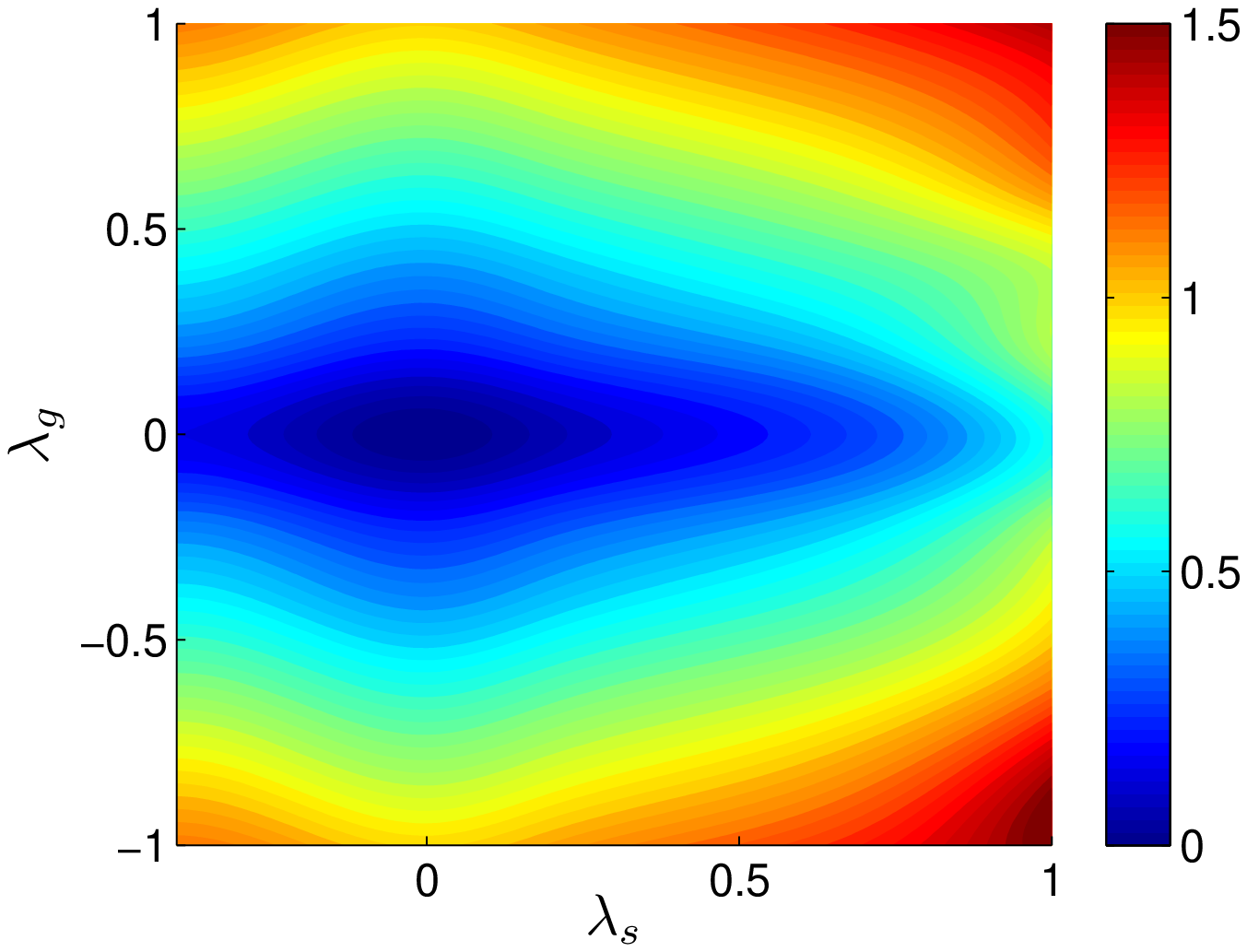}}}
  \mbox{
  \subfigure[model \#1: simultaneous]{\includegraphics[width=0.33\textwidth]{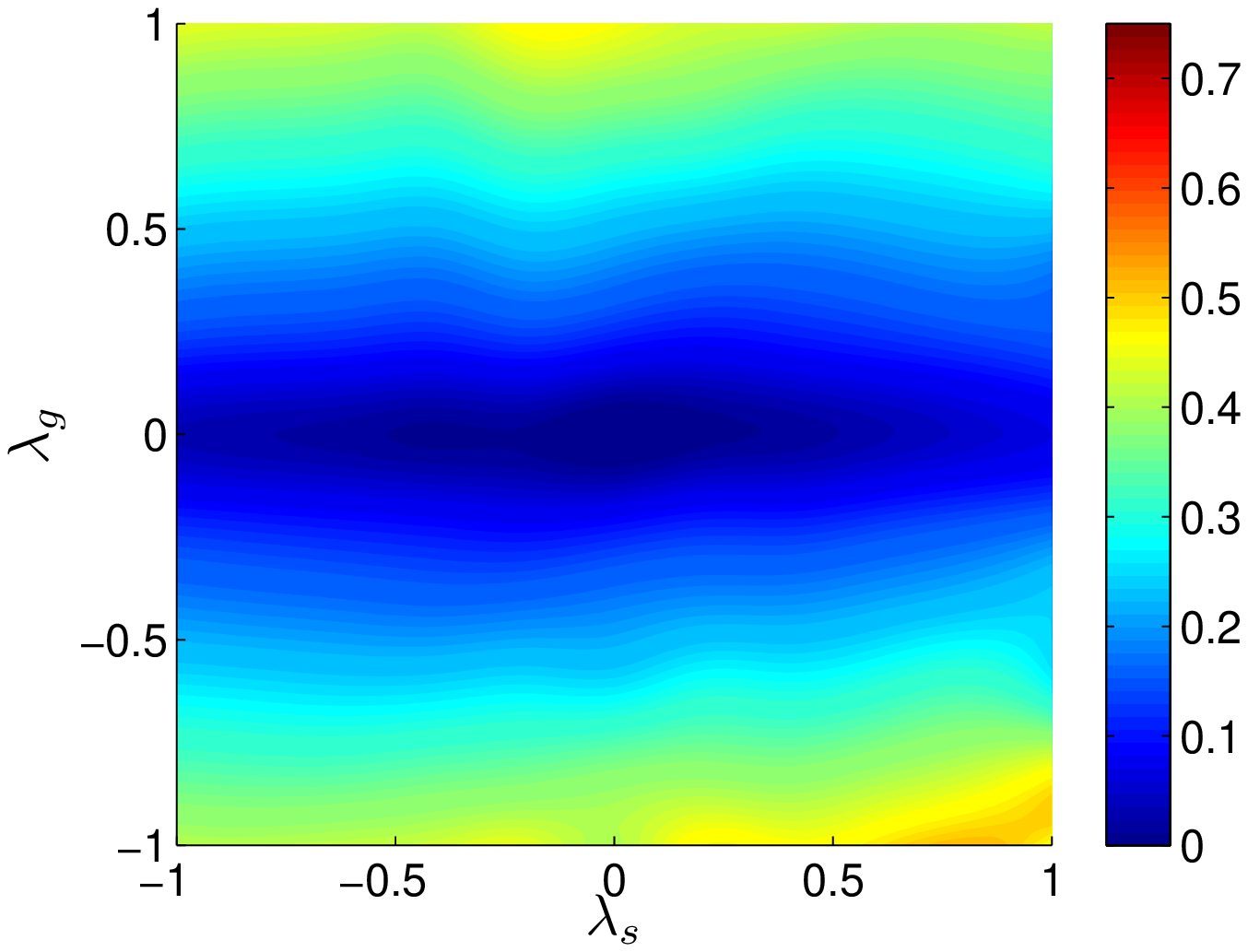}}
  \subfigure[model \#2: simultaneous]{\includegraphics[width=0.33\textwidth]{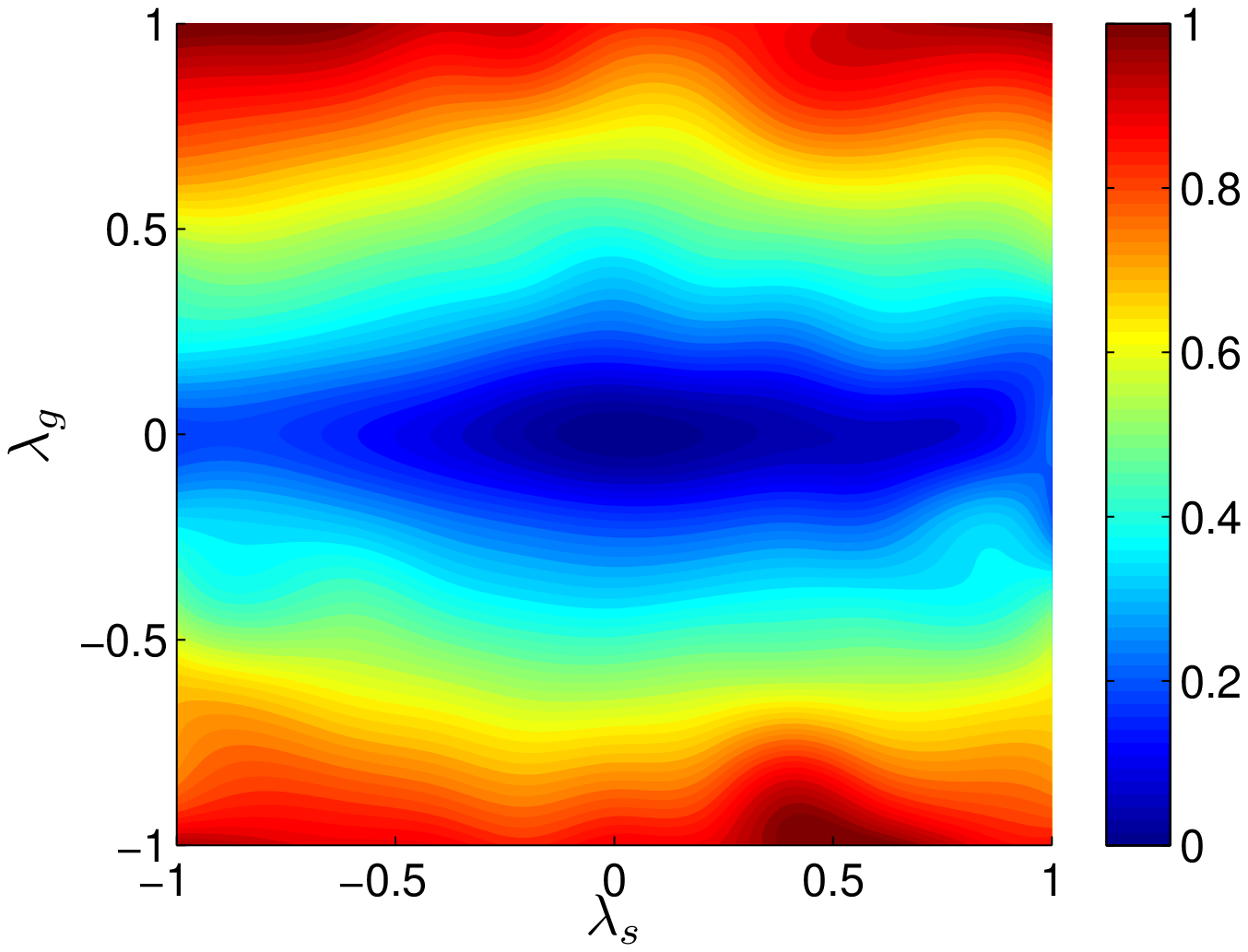}}
  \subfigure[model \#3: simultaneous]{\includegraphics[width=0.33\textwidth]{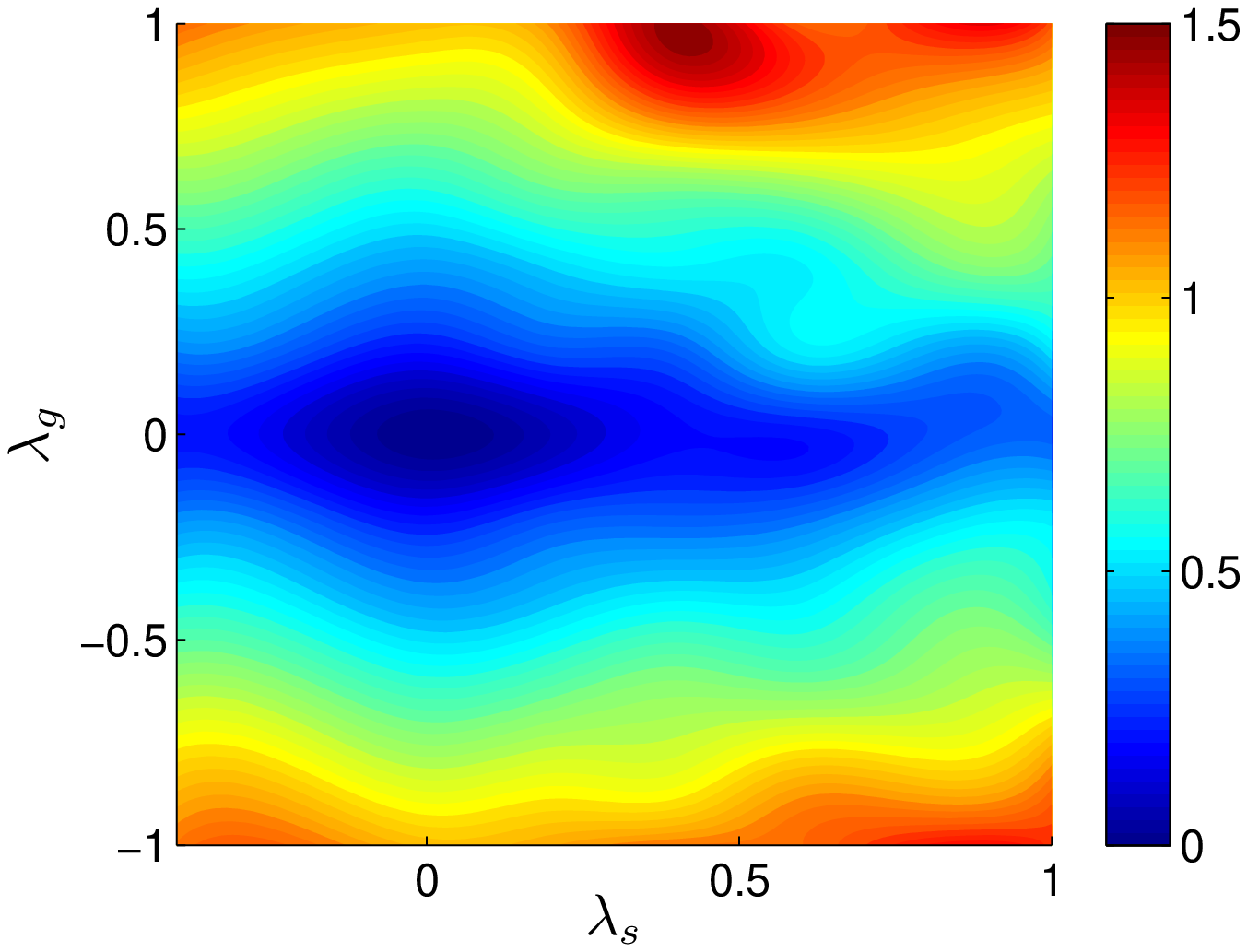}}}
  \mbox{
  \subfigure[model \#1: interchanging $N = 1$]{\includegraphics[width=0.33\textwidth]{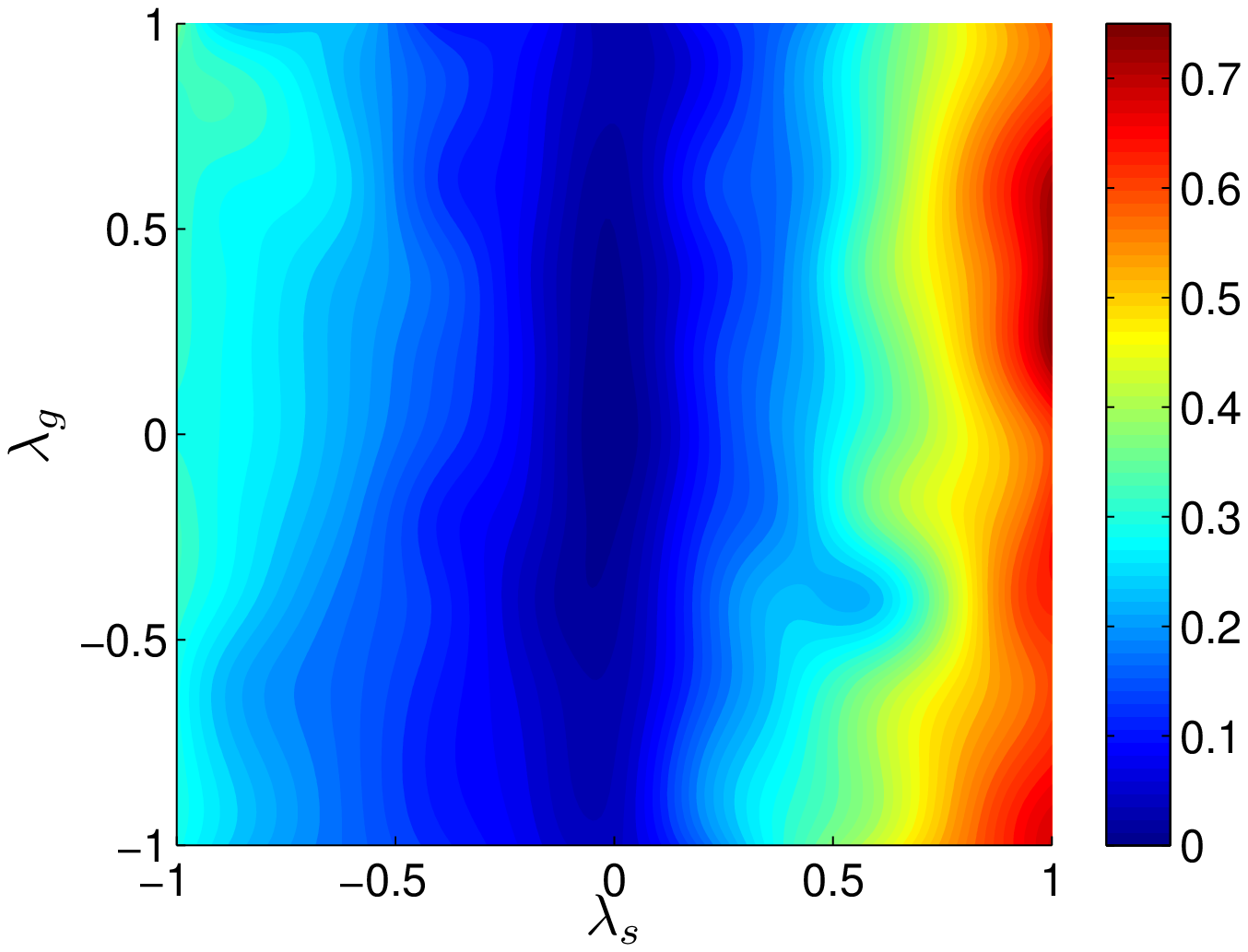}}
  \subfigure[model \#2: interchanging $N = 1$]{\includegraphics[width=0.33\textwidth]{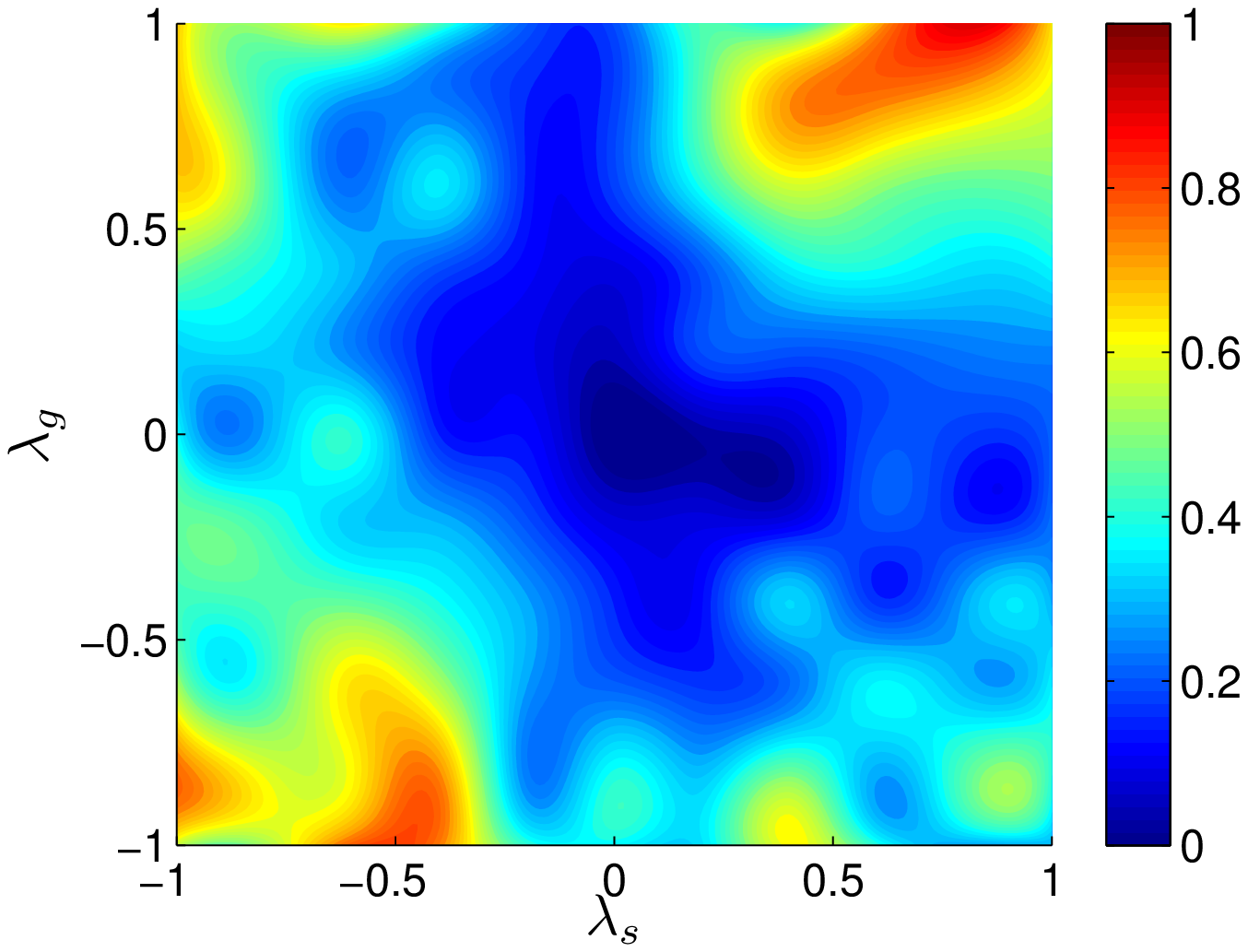}}
  \subfigure[model \#3: interchanging $N = 1$]{\includegraphics[width=0.33\textwidth]{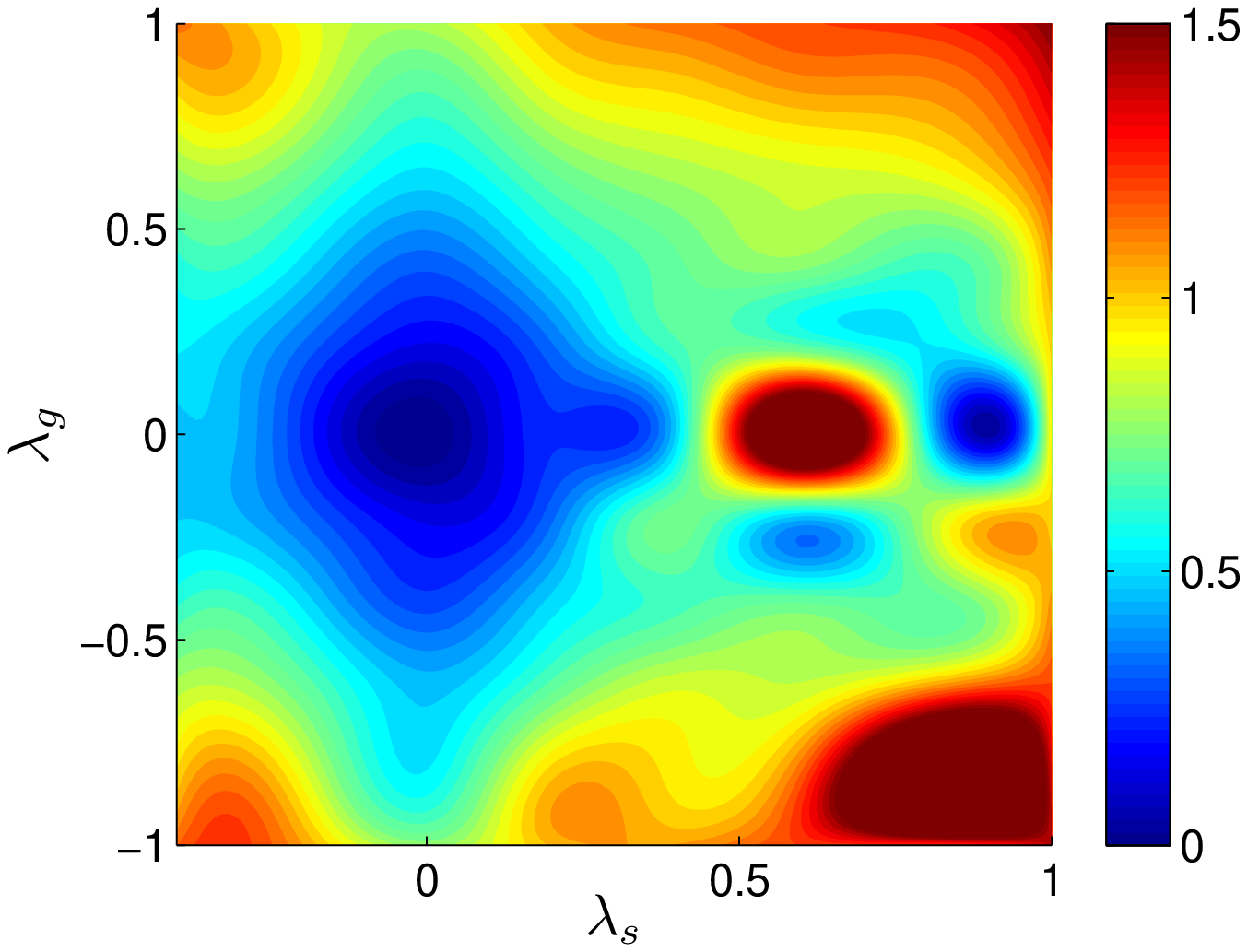}}}
  \end{center}
  \caption{Convergence analysis for both $s(t)$ and $g(t)$ performed for (a,d,g) model~\#1, (b,e,h) model~\#2, and
    (c,f,i) model~\#3 by evaluating the solution norms $\dfrac{\| \bs - s_{\rm{true}} \|_{L_2}}{\| s_{\rm{true}} \|_{L_2}} +
    \dfrac{\| \bg - g_{\rm{true}} \|_{L_2}}{\| g_{\rm{true}} \|_{L_2}}$ for three cases: (a,b,c) with no preconditioning
    by method~\#1, (d,e,f) with optimal preconditioning by method~\#1, and (g,h,i) with optimal preconditioning
    by method~\#2 ($N=1$).}
  \label{fig:conv_sol_sg}
\end{figure}

Finally, we validate our computational approach in terms of convergence of $s(t)$ and $g(t)$ to their respective
global optimal solutions by solving the same optimization problem for all three models using simultaneous identification
of $s(t)$ and $g(t)$ by method~\#1 and interchanging ($N=1$) method~\#2 starting with different initial guesses.
As we did it separately for $s(t)$ and $g(t)$, new initial guesses $s_{\rm{ini},\lambda_s}$ and $g_{\rm{ini},\lambda_g}$
are parameterized with respect to their proximity to their respective global minimizers $s_{\rm{true}}$ and $g_{\rm{true}}$
as shown by \eqref{eq:s_convex_comb} and \eqref{eq:g_convex_comb}. When using preconditioning procedure \eqref{eq:helm},
we apply smoothing parameters $\ell^*_{s,i}$ and $\ell^*_{g,i}, \, i = 1, 2, 3$.

The results of this convergence test for all three models are shown in Figure~\ref{fig:conv_sol_sg} for three cases:
\begin{itemize}
  \item with no preconditioning and simultaneous identification by method \#1 (a,b,c);
  \item with optimal preconditioning and simultaneous identification by method \#1 (d,e,f);
  \item with optimal preconditioning and interchanging ($N=1$) method \#2 (g,h,i).
\end{itemize}
For all cases color represents the norm $\dfrac{\| \bs - s_{\rm{true}} \|_{L_2}}{\| s_{\rm{true}} \|_{L_2}} +
\dfrac{\| \bg - g_{\rm{true}} \|_{L_2}}{\| g_{\rm{true}} \|_{L_2}}$. As seen in Figure~\ref{fig:conv_sol_sg}(a-f),
all three models show that the interval of convergence for $s(t)$ is much larger than that for $g(t)$ when
using simultaneous (method~\#1) reconstruction with and without preconditioning. We should note that adding
preconditioning improves slightly convergence for $s(t)$. For $g(t)$ optimal preconditioning is added
only for model~\#2 but without noticeable effect compared with $s(t)$. However, as seen in
Figure~\ref{fig:conv_sol_sg}(g-i), convergence for $g(t)$ could be improved by applying interchanging ($N=1$)
method \#2.

\section{Conclusions}
\label{sec:conclusions}
This paper presents computational analysis of the inverse Stefan type free boundary problem, where information
on the boundary heat flux is missing and must be
found along with the temperature and the free boundary. The motivation for this type of inverse problem arose in particular from the modeling of bioengineering problems on the laser ablation of biological tissues through a Stefan problem~\eqref{eq:pde-1}--\eqref{eq:pde-finaltemp}, where the free boundary $s(t)$ is the ablation depth at the moment $t$. We pursued the optimal control framework introduced in \cite{abdulla13,abdulla15}, where boundary heat flux and free boundary are components of the control vector,
and optimality criteria consist of the minimization of the quadratic declinations from the available measurements of the temperature distribution at the final moment, phase transition temperature on the free boundary, and  the final position of the free boundary. In recent papers~\cite{abdulla16,abdulla17}, the Fr\'echet differentiability and necessary optimality condition in Besov spaces
were established under minimal conditions on the data. In this paper we developed a gradient descent algorithm in Hilbert-Besov space $H^1\times H^1$ based on the formula for the Fr\'echet gradient which is an element of the dual space. By applying Riesz representation theorem, we implement {\it preconditioning} to calculate an equivalent form for the Fr\'echet gradient in $H^1\times H^1$ with increased regularity. Three model examples with various levels of complexity are considered. The following are the major outcomes:
\begin{itemize}
\item Gradient method with and without preconditioning is demonstrated to be an effective method for reconstruction of the local and global optimal control. Calibration of the preconditioning parameter demonstrates that there is an intermediate range of the parameter with best performance with respect to both cost functional and control criteria for the reconstruction of the free boundary. In general, preconditioning with optimal preconditioning parameter improves the convergence rate, but with the expense of increased computational time. 
\item Gradient method for the reconstruction of the free boundary is tested in the presence of additional measurements
on the position of the free boundary at some time instances with possible noise. Comparative analysis of alternative approaches when piecewise-linear interpolation of additional measurements is used as either the initial guess or as a regularization centroid of the Tikhonov regularization method. In the former, preconditioning has an advantage if the number of measurements are low, and has no improvement or even a negative effect otherwise. In the latter case, it is demonstrated that the Tikhonov regularization with optimal choice of the regularization parameter has a similar convergence effect as the original method without preconditioning, but with updated initial guess. These outcomes are consistent with sufficiently large Gaussian noise, up to $10\%$, added to the measurements. Hence, additional measurements of the free boundary have a regularizing effect on the reconstruction of the free boundary.
\item Gradient method is tested for identification of the free boundary and other control parameters. We developed alternative approaches such as simultaneous identification vs. identification in $N$-interchanging order or nested optimization, meaning that the identification algorithm switches between control parameters for every $N$ optimization iterations. We pursued two cases with $N=1$ and $N=5$ in all model examples. All three methods were accompanied with preconditioning with optimal choice of the parameters. Extensive comparative analysis demonstrates that the advantage of the methods are dependent on model complexity: all three methods worked well in the simplest model, nested optimization has an advantage in the model of moderate complexity, and simultaneous identification has  a clear advantage in the most complex model. 
\end{itemize}



\bibliographystyle{elsarticle-harv}
\bibliography{biblio_Abdulla,biblio_Bukshtynov,biblio_OPT}

\begin{thebibliography}{46}
\expandafter\ifx\csname natexlab\endcsname\relax\def\natexlab#1{#1}\fi
\expandafter\ifx\csname url\endcsname\relax
  \def\url#1{\texttt{#1}}\fi
\expandafter\ifx\csname urlprefix\endcsname\relax\def\urlprefix{URL }\fi

\bibitem[{Abdulla(2013)}]{abdulla13}
Abdulla, U.~G., 2013. On the optimal control of the free boundary problems for
  the second order parabolic equations. {I}. {W}ell-posedness and convergence
  of the method of lines. Inverse Problems and Imaging 7~(2), 307--340.

\bibitem[{Abdulla(2016)}]{abdulla15}
Abdulla, U.~G., 2016. On the optimal control of the free boundary problems for
  the second order parabolic equations. {II}. {C}onvergence of the method of
  finite differences. Inverse Problems and Imaging 10~(4), 869--898.

\bibitem[{Abdulla et~al.(2017)Abdulla, Cosgrove, and Goldfarb}]{abdulla17}
Abdulla, U.~G., Cosgrove, E., Goldfarb, J., 2017. On the {F}rechet
  differentability in optimal control of coefficients in parabolic free
  boundary problems. Evolution Equations and Control Theory 6~(4), 319--344.

\bibitem[{Abdulla and Goldfarb(2018)}]{abdulla16}
Abdulla, U.~G., Goldfarb, J., 2018. {F}rechet differentability in {B}esov
  spaces in the optimal control of parabolic free boundary problems. Journal of
  Inverse and Ill-posed Problems 26~(2).

\bibitem[{Abdulla and Poggi(2018)}]{abdullapoggi}
Abdulla, U.~G., Poggi, B., 2018. Optimal control of the multiphase stefan
  problem. Applied Mathematics and Optimization 77~(2).

\bibitem[{Baumeister(1980)}]{baumeister80}
Baumeister, J., 1980. Zur optimal {Steuerung} von frien {Randwertausgaben}.
  ZAMM 60, 335--339.

\bibitem[{Bell(1981)}]{bell81}
Bell, J.~B., 1981. The non-characteristic {Cauchy} problem for a class of
  equations with time dependence. {I.} problem in one space dimension. SIAM
  Journal on Mathematical Analysis 12~(5), 759--777.

\bibitem[{Berger(1977)}]{Berger77}
Berger, M.~S., 1977. Nonlinearity and Functional Analysis. Acad. Press, New
  York.

\bibitem[{Besov et~al.(1979{\natexlab{a}})Besov, Ilin, and Nikolskii}]{besov79}
Besov, O.~V., Ilin, V.~P., Nikolskii, S.~M., 1979{\natexlab{a}}. Integral
  Representations of Functions and Imbedding Theorems. Vol. Vol. 1. John Wiley
  \& Sons.

\bibitem[{Besov et~al.(1979{\natexlab{b}})Besov, Ilin, and
  Nikolskii}]{besov79a}
Besov, O.~V., Ilin, V.~P., Nikolskii, S.~M., 1979{\natexlab{b}}. Integral
  Representations of Functions and Imbedding Theorems. Vol. Vol. 2. John Wiley
  \& Sons.

\bibitem[{Budak and Vasileva(1972)}]{budak72}
Budak, B.~M., Vasileva, V.~N., 1972. On the solution of the inverse {Stefan}
  problem. Soviet Mathematics Doklady 13, 811--815.

\bibitem[{Budak and Vasileva(1973)}]{budak73}
Budak, B.~M., Vasileva, V.~N., 1973. On the solution of {Stefan}'s converse
  problem {II}. USSR Computational Mathematics and Mathematical Physics 13,
  97--110.

\bibitem[{Budak and Vasileva(1974)}]{budak74}
Budak, B.~M., Vasileva, V.~N., 1974. The solution of the inverse {Stefan}
  problem. USSR Computational Mathematics and Mathematical Physics 13~(1),
  130--151.

\bibitem[{Bukshtynov and Protas(2013)}]{Bukshtynov13}
Bukshtynov, V., Protas, B., 2013. Optimal reconstruction of material properties
  in complex multiphysics phenomena. Journal of Computational Physics 242,
  889--914.

\bibitem[{Bukshtynov et~al.(2011)Bukshtynov, Volkov, and Protas}]{Bukshtynov11}
Bukshtynov, V., Volkov, O., Protas, B., 2011. On optimal reconstruction of
  constitutive relations. Physica D: Nonlinear Phenomena 240~(16), 1228--1244.

\bibitem[{Cannon(1964)}]{cannon64}
Cannon, J.~R., 1964. A {Cauchy} problem for the heat equation. Annali di
  Matematica Pura Ed Applicata 66~(1), 155--165.

\bibitem[{Cannon and Jr.(1967)}]{cannon67}
Cannon, J.~R., Jr., J.~D., 1967. The {Cauchy} problem for the heat equation.
  SIAM Journal on Numerical Analysis 4~(3), 317--336.

\bibitem[{Carasso(1982)}]{carasso82}
Carasso, A., June 1982. Determining surface temperatures from interior
  observations. SIAM Journal on Applied Mathematics 42~(3), 558--574.

\bibitem[{Davis(2004)}]{umfpack}
Davis, T.~A., 2004. Algorithm 832: {UMFPACK V4.3} -- an unsymmetric-pattern
  multifrontal method. ACM Transactions on Mathematical Software (TOMS) 30~(2),
  196--199.

\bibitem[{Ewing(1979)}]{ewing79}
Ewing, R.~E., September 1979. The {Cauchy} problem for a linear parabolic
  equation. Journal of Mathematical Analysis and Applications 71~(1), 167--186.

\bibitem[{Ewing and Falk(1979)}]{ewing79a}
Ewing, R.~E., Falk, R., 1979. Numerical approximation of a {Cauchy} problem for
  a parabolic partial differential equations. Mathematics of Computation
  33~(148), 1125--1144.

\bibitem[{Fasano and Primicerio(1977)}]{fasano77}
Fasano, A., Primicerio, M., 1977. General free boundary problems for heat
  equations. Journal of Mathematical Analysis and Applications 57~(3),
  694--723.

\bibitem[{Gol'dman(1997)}]{goldman97}
Gol'dman, N.~L., 1997. Inverse {Stefan} Problems. Kluwer Academic Publishers
  Group, Dodrecht.

\bibitem[{Hecht(2012)}]{FreeFem12}
Hecht, F., 2012. New development in {FreeFem++}. Journal of Numerical
  Mathematics 20~(3-4), 251--265.

\bibitem[{Hoffman and Niezgodka(1981)}]{hoffman81}
Hoffman, K.~H., Niezgodka, M., 1981. Control of parabolic systems involving
  free boundaries. In: Proceedings of the International Conference on Free
  Boundary Problems.

\bibitem[{Hoffman and Sprekels(1982)}]{hoffman82}
Hoffman, K.~H., Sprekels, J., 1982. Real time control of free boundary in a
  two-phase {Stefan} problem. Numerical Functional Analysis and Optimization 5,
  47--76.

\bibitem[{Hoffman and Sprekels(1986)}]{hoffman86}
Hoffman, K.~H., Sprekels, J., 1986. On the identification of heat conductivity
  and latent heat conductivity as latent heat in a one-phase {Stefan} problem.
  Control and Cybernetics 15, 37--51.

\bibitem[{Jochum(1980{\natexlab{a}})}]{jochum80a}
Jochum, P., 1980{\natexlab{a}}. The inverse {Stefan} problem as a problem of
  nonlinear approximation theory. Journal of Approximation Theory 30, 37--51.

\bibitem[{Jochum(1980{\natexlab{b}})}]{jochum80}
Jochum, P., 1980{\natexlab{b}}. The numerical solution of the inverse {Stefan}
  problem. Numerical Mathematics 34, 411--429.

\bibitem[{Knabner(1983)}]{knabner83}
Knabner, P., 1983. Stability theorems for general free boundary problems of the
  {Stefan} type and applications. Applied Nonlinear Functional Analysis,
  Methoden und Verfahren der Mathematischen Physik 25, 95--116.

\bibitem[{Kufner et~al.(1977)Kufner, John, and Fu\v{c}ik}]{kufner77}
Kufner, A., John, O., Fu\v{c}ik, S., 1977. Function Spaces. Noordhoff
  International Publishing, Leyden, The Netherlands.

\bibitem[{Ladyzhenskaya et~al.(1968)Ladyzhenskaya, Solonnikov, and
  Uraltseva}]{ladyzhenskaya68}
Ladyzhenskaya, O.~A., Solonnikov, V.~A., Uraltseva, N.~N., 1968. Linear and
  Quasilinear Equations of Parabolic Type. Vol.~23 of Translations of
  Mathematical Monographs. American Mathematical Society, Providence, R. I.

\bibitem[{Lurye(1975)}]{lurye75}
Lurye, K.~A., 1975. Optimal Control in Problems of Mathematical Physics.
  Moscow. Nauka.

\bibitem[{Niezgodka(1979)}]{niezgodka79}
Niezgodka, M., 1979. Control of parabolic systems with free boundaries -
  application of inverse formulation. Control and Cybernetics 8, 213--225.

\bibitem[{Nikol'skii(1975)}]{nikolskii75}
Nikol'skii, S.~M., 1975. Approximation of Functions of Several Variables and
  Imbedding Theorems. Springer-Verlag, New York-Heidelberg.

\bibitem[{Nocedal and Wright(2006)}]{Nocedal06}
Nocedal, J., Wright, S.~J., 2006. Numerical Optimization, 2nd Edition.
  Springer, New York.

\bibitem[{Nochetto and C.Verdi(1987/88)}]{nochetto87}
Nochetto, R.~H., C.Verdi, 1987/88. The combined use of nonlinear {Chernoff}
  formula with a regularization procedure for two-phase {Stefan} problems.
  Numerical Functional Analysis and Optimization 9, 1177--1192.

\bibitem[{Press et~al.(2007)Press, Teukolsky, Vetterling, and
  Flannery}]{NumericalRecipes07}
Press, W.~H., Teukolsky, S.~A., Vetterling, W.~T., Flannery, B.~P., 2007.
  Numerical Recipes: The Art of Scientific Computing, 3rd Edition. Cambridge
  University Press.

\bibitem[{Primicerio(1982)}]{primicerio82}
Primicerio, M., 1982. The occurence of pathologies in some {Stefan}-like
  problems. In: Albrecht, J., Collatz, L., Hoffman, K.~H. (Eds.), Numerical
  Treatment of Free Boundary-Value problems. Vol.~58. ISNM, Birkhauser Verlag,
  Basel, pp. 233--244.

\bibitem[{Protas et~al.(2004)Protas, Bewley, and Hagen}]{ProtasBewleyHagen04}
Protas, B., Bewley, T., Hagen, G., 2004. A computational framework for the
  regularization of adjoint analysis in multiscale {PDE} systems. Journal of
  Computational Physics 195~(1), 49--89.

\bibitem[{Sagues(1982)}]{sagues82}
Sagues, C., 1982. Simulation and optimal control of free boundary. In:
  Albrecht, J., Collatz, L., Hoffman, K.~H. (Eds.), Numerical Treatment of Free
  Boundary-Value problems. Vol.~58. ISNM, Birkhauser Verlag, Basel, pp.
  270--287.

\bibitem[{Sherman(1971)}]{sherman71}
Sherman, B., 1971. General one-phase {Stefan} problems and free boundary
  problems for the heat equation with {Cauchy} data prescribed on the free
  boundary. SIAM J. Appl. Math. 20, 557--570.

\bibitem[{Solonnikov(1964)}]{solonnikov64}
Solonnikov, V.~A., 1964. A-priori estimates for solutions of second-order
  equations of parabolic type. Vol.~70 of Trudy Matematischeskogo instituta im.
  V. A. Steklova. Nauka, Moscow-Leningrad.

\bibitem[{Talenti and Vessella(1982)}]{talenti82}
Talenti, G., Vessella, S., June 1982. A note on an ill-posed problem for the
  heat equation. Journal of the Austrailian Mathematical Society 32~(3),
  358--368.

\bibitem[{Vasil'ev(1969)}]{vasilev69}
Vasil'ev, F.~P., 1969. The existence of a solution to a certain optimal
  {Stefan} problem. Computational Methods and Programming, 110--114.

\bibitem[{Yurii(1980)}]{yurii80}
Yurii, A.~D., 1980. On an optimal {Stefan} problem. Doklady Akademii nauk SSSR
  251, 1317--1321.

\end{thebibliography}

\end{document}